\documentclass[12pt,a4paper]{article}
\usepackage{amssymb,amsmath,latexsym,times,graphicx,standard}
\usepackage{pdfsync}

\textheight 
            675pt
\textwidth 
           460pt

\newcommand{\qed}{\hfill$\Box$\bigskip}

\newcommand{\eqn}[1]{(\ref{#1})}

\newcommand{\R}{\mathbb{R}}
\newcommand{\Z}{\mathbb{Z}}

\newcommand{\ve}{{\varepsilon}}

\newcommand{\dx}{\, \mathrm{d}x}

\newcommand{\dy}{\, \mathrm{d}y}
\newcommand{\dty}{\, \mathrm{d}{\tilde y}}

\newcommand{\dmu}{\, \mathrm{d}\mu}
\newcommand{\dtmu}{\, \mathrm{d}{\tilde \mu}}

\newcommand{\dkappa}{\, \mathrm{d}\kappa}

\newcommand{\ee}{\mathrm{e}}
\newcommand{\ii}{\mathrm{i}}

\renewcommand{\AA}{{\mathcal A}}
\newcommand{\BB}{{\mathcal B}}
\newcommand{\CC}{{\mathcal C}}
\newcommand{\DD}{{\mathcal D}}

\newcommand{\FF}{{\mathcal F}}
\newcommand{\GG}{{\mathcal G}}
\newcommand{\HH}{{\mathcal H}}

\newcommand{\LL}{{\mathcal L}}

\newcommand{\NN}{{\mathcal N}}
\newcommand{\RR}{{\mathcal R}}

\newcommand{\UU}{{\mathcal U}}
\newcommand{\XX}{{\mathcal X}}
\newcommand{\YY}{{\mathcal Y}}
\newcommand{\ZZ}{{\mathcal Z}}

\DeclareMathOperator{\sech}{sech}
\DeclareMathOperator{\cosech}{cosech}
\DeclareMathOperator{\supp}{supp}

\DeclareMathOperator{\re}{Re}
\DeclareMathOperator{\im}{Im}

\newtheorem{theorem}{Theorem}[section]
\newtheorem{lemma}[theorem]{Lemma}
\newtheorem{proposition}[theorem]{Proposition}
\newtheorem{corollary}[theorem]{Corollary}

\newcommand{\llangle}{\langle{\mskip-4mu}\langle}
\newcommand{\rrangle}{\rangle{\mskip-4mu}\rangle}

\newcounter{count}

\title{A dimension-breaking phenomenon for water waves with weak surface tension}

\author{M. D. Groves\thanks{Fachrichtung 6.1 - Mathematik, Universit\"{a}t des Saarlandes,
Postfach 151150, 66041 Saarbr\"{u}cken, Germany; 
Department of Mathematical Sciences, Loughborough
University, Loughborough, LE11 3TU, UK
} \and
S. M. Sun\footnote{Department of Mathematics, Virginia Polytechnic Institute and State University,
Blacksburg, VA 24061, USA}
\and E. Wahl\'en\footnote{Centre for Mathematical Sciences, Lund University, PO Box 118, 22100 Lund, Sweden}}
\date{}

\begin{document}
\maketitle

\begin{abstract}
\noindent

It is well known that the water-wave problem 
with weak surface tension has small-amplitude line solitary-wave solutions which to leading order 
are described by the nonlinear Schr\"odinger equation. 
The present paper contains an existence theory for three-dimensional periodically modulated solitary-wave solutions which have a solitary-wave profile in the direction of propagation and are periodic in the transverse direction;
they emanate from the line solitary waves in a dimension-breaking bifurcation.
In addition, it is shown that the line solitary waves are linearly unstable to long-wavelength transverse perturbations.
The key to these results is a formulation of the water wave problem as an evolutionary system in which the transverse
horizontal variable plays the role of time, a  careful study of the purely imaginary spectrum of the operator
obtained by linearising the evolutionary system at a line solitary wave, and an application of an
infinite-dimensional version of the classical Lyapunov centre theorem.
\end{abstract}

\section{Introduction} \label{Introduction}

In this article we consider the three-dimensional irrotational flow of a perfect fluid of unit density subject to the 
forces of gravity and surface tension. The fluid is bounded below by a rigid horizontal bottom 
$\{y=0\}$ and above by a free surface $\{y=h+\eta(x,z, t)\}$. Our primary interest is in travelling waves moving
without change of shape and with constant speed $c>0$ from left to right in the $x$-direction. In a moving frame
of reference in which such waves are stationary the equations of motion are
Laplace's equation
\begin{equation}
\phi_{xx}+\phi_{yy}+\phi_{zz}=0, \qquad\quad 0<y<1+\eta, \label{eq:Linear BVP}
\end{equation}
for the Eulerian velocity potential $\phi$ describing the flow, with boundary conditions
\begin{align}
\phi_y&=0, && y=0, \label{eq:Lower kinematic} \\
\phi_y&=\eta_t - \eta_x+\eta_x \phi_x+\eta_z \phi_z, && y=1+\eta, \label{eq:Upper kinematic}
\end{align}
and
\begin{equation}
\phi_t-\phi_x+\frac12 (\phi_x^{2}+\phi_y^{2}+\phi_z^2)+\alpha \eta  -\beta\! \left[\frac{\eta_x}{\sqrt{1+\eta_x^2+\eta_z^2}} \right]_x\!\!- \beta\! \left[\frac{\eta_z}{\sqrt{1+\eta_x^2+\eta_z^2}} \right]_z\!\!\!=0, \quad y=1+\eta. \label{eq:Bernoulli condition}
\end{equation}
Here we have introduced dimensionless variables, choosing $h$ as length scale and
$h/c$ as time scale; the parameters $\alpha$ and $\beta$ are defined in terms of the Froude
and Bond numbers $F=c/\sqrt{gh}$ and $\tau=\sigma/gh^2$ by the formulae
$\alpha=1/F^2$ and $\beta=\tau/F^2$, in which $g$ is the acceleration due to
gravity and $\sigma$ is the coefficient of surface tension.

\begin{figure}[ht]
\hspace{4cm}\includegraphics[width=8cm]{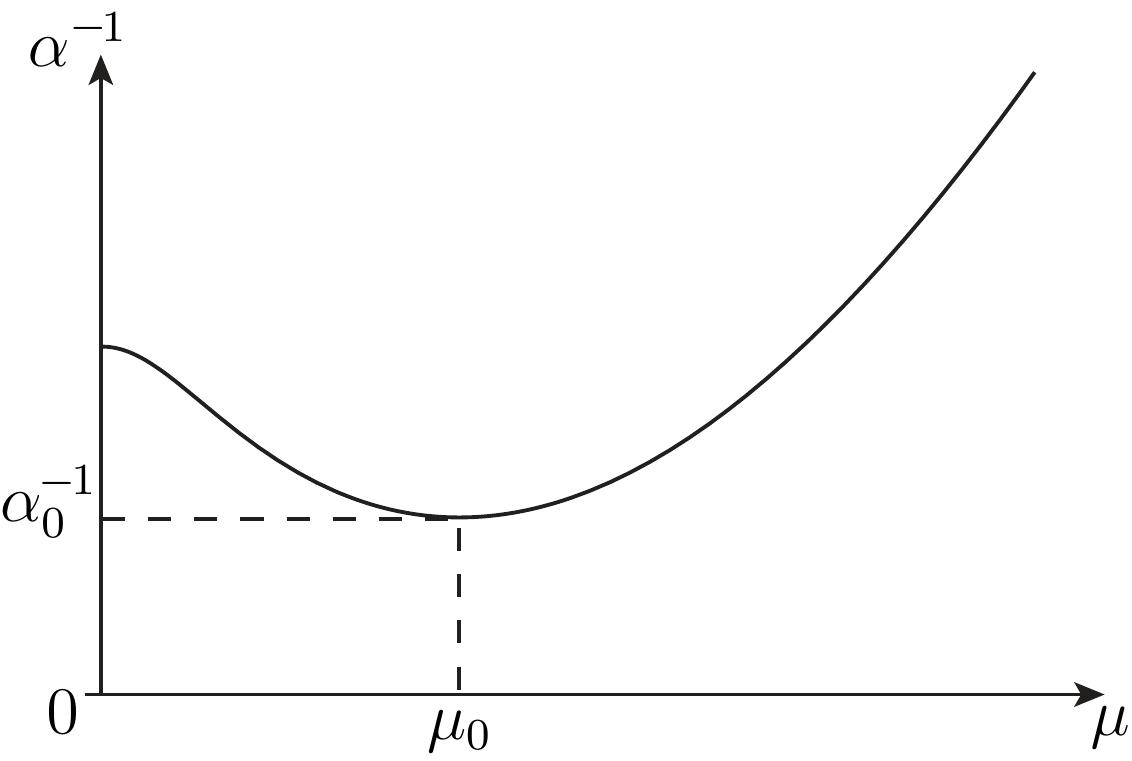}
{\it
\caption{The linear dispersion relation for a sinusoidal travelling wave train with wave number $\mu$ and $\tau=\tau_0$, where $\tau_0 \in (0,1/3)$.
\label{dispersion relation}}}
\end{figure}

Before describing the present contribution, let us briefly review the classical weakly nonlinear
theory with weak surface tension (that is, $0<\tau<1/3$)
for stationary wave packets consisting of a slowly varying envelope modulating a
periodic wave train. Figure \ref{dispersion relation}
shows the linear dispersion relation for a sinusoidal travelling wave train with wave number $\mu$;
for each fixed value $\tau_0 \in (0,1/3)$ of $\tau$ the dispersion curve has a unique
minimum at $(\mu,\alpha^{-1})=(\mu_0,\alpha_0^{-1})$. (The relationship between $\alpha_0$,
$\beta_0=\tau_0\alpha_0$ and $\mu_0$ can be expressed in the form
\[
\alpha_0=\frac{\mu_0^2}{2}\cosech^2 \mu_0 + \frac{\mu_0}{2} \coth \mu_0, \qquad
\beta_0=-\frac{1}{2}\cosech^2 \mu_0 +\frac{1}{2\mu_0} \coth \mu_0,
\]
which defines a curve in the $(\alpha,\beta)$-plane parameterised by $\mu \in (0,\infty)$.) Substituting
\[
\alpha=\alpha_0+\ve^2, \quad \beta=\beta_0, \quad \eta(x,z) = \frac{1}{2}\ve\zeta(\ve x, \ve z) \mathrm{e}^{\mathrm{i} \mu_0 x}+\frac{1}{2}\ve
\overline{\zeta(\ve x, \ve z)} \mathrm{e}^{-\mathrm{i} \mu_0 x}
\]
into equations \eqn{eq:Linear BVP}--\eqn{eq:Bernoulli condition}, one finds that
to leading order $\zeta$ satisfies the elliptic-elliptic Davey-Stewartson system
\begin{eqnarray*}
\zeta-A_1\zeta_{XX}-A_2\zeta_{ZZ} - A_3|\zeta|^2 \zeta +4 A_4\zeta \psi_X & = & 0, \\
- (1-\alpha_0^{-1})\psi_{XX}-\psi_{ZZ}-A_4(|\zeta|^2)_X & = & 0,
\end{eqnarray*}
where $X=\ve x$, $Z=\ve z$ (see Djordjevic \& Redekopp \cite{DjordjevicRedekopp77},
Ablowitz \& Segur \cite{AblowitzSegur79}, noting the misprint in equation
(2.24d), and
Sulem \& Sulem \cite[\S11.1.1]{SulemSulem});
formulae for the coefficients $A_1$, $A_3$ and $A_4$ are
given in the Appendix, while $A_2=\mu_0^{-1}\coth \mu_0$. In the
special case of $z$-independent waves these equations reduce to the cubic
nonlinear Schr\"{o}dinger equation
\[
\zeta-A_1 \partial_X^2 \zeta-A_5|\zeta|^2\zeta=0,
\]
where $A_5=A_3+4(1-\alpha_0^{-1})^{-1}A_4^2$.

\begin{figure}[ht]
\hspace{2cm}\includegraphics[width=12cm]{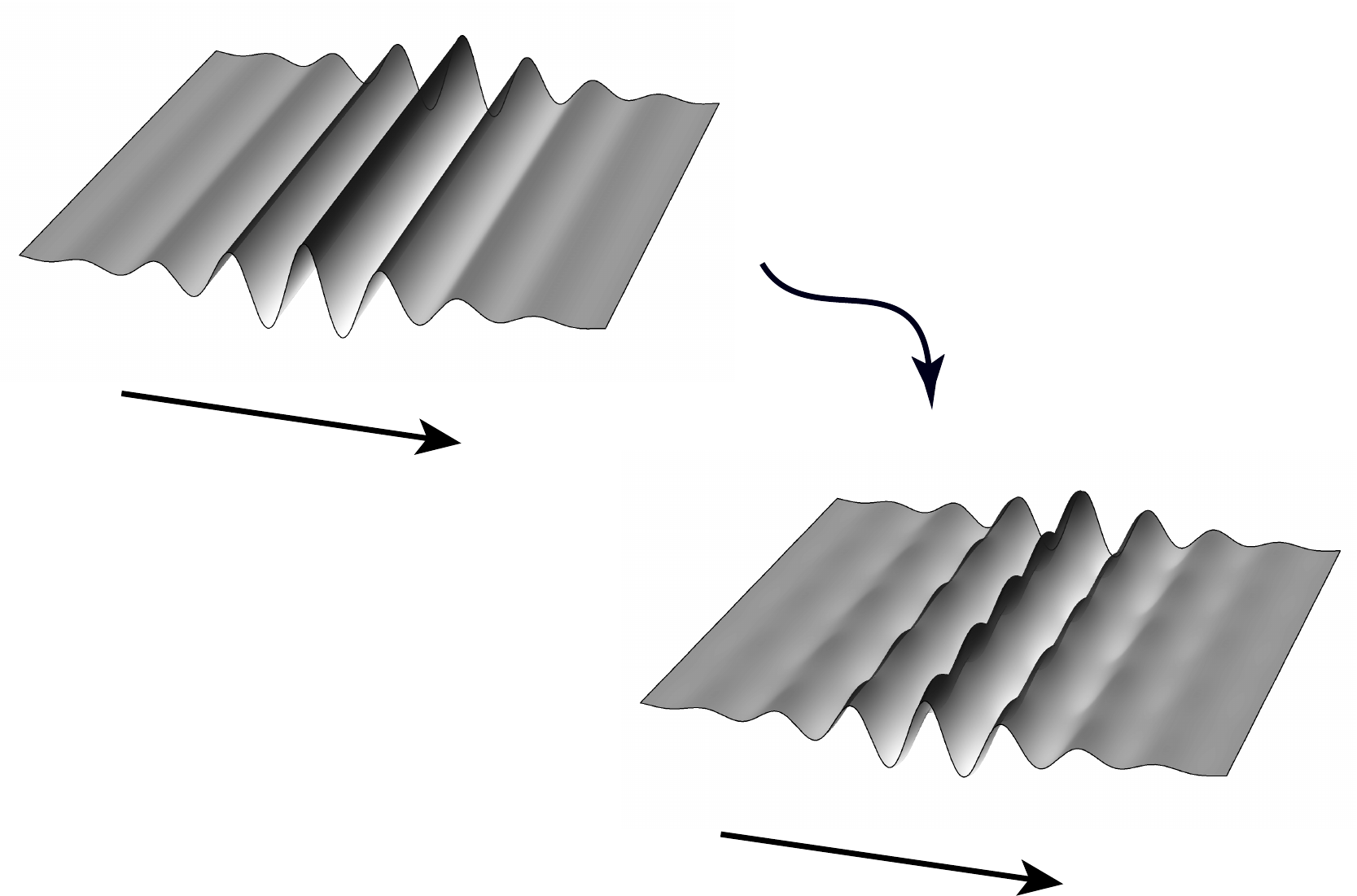}
{\it
\caption{The periodically modulated solitary wave on the right emerges from the line solitary wave on the left in a dimension breaking bifurcation.
\label{dbp}}}
\end{figure}

The water-wave problem \eqn{eq:Linear BVP}--\eqn{eq:Bernoulli condition} with weak surface
tension
has two families of $z$-independent,
small-amplitude, symmetric solitary-wave solutions $(\eta^\star_{\ve},\phi^\star_{\ve})$. We refer to
these waves, which were discovered by Iooss and Kirchg\"assner \cite{IoossKirchgaessner90}
(see also Dias \& Iooss \cite{DiasIooss93} and Iooss \& P\'{e}rou\`{e}me \cite{IoossPeroueme93}),
as {\em line solitary waves}; they are stationary solutions of \eqn{eq:Linear BVP}--\eqn{eq:Bernoulli condition} with
\[
\alpha=\alpha_0+\ve^2, \quad \beta=\beta_0, \quad
\eta_\ve^\star(x)=\ve \zeta^\star(\ve x)\cos(\mu_0 x)+O(\ve^2),
\]
where the function
\[
\zeta^\star(X)=
\pm \left(\frac{2}{A_5}\right)^{1/2}\sech\left(\frac{X}{A_1^{1/2}}\right)
\]
is a symmetric solution of the cubic nonlinear Schr\"odinger equation (and of course a symmetric
line solitary-wave solution of the Davey-Stewartson system).
The line solitary waves are wave packets consisting of a decaying envelope which modulates an underlying periodic wave train with wave number $\mu_0$
(see Figure \ref{dbp}).

In this paper we again consider weak surface tension ($0<\tau<1/3$) and
prove the existence of three-dimensional \emph{periodically modulated solitary waves} which have a solitary-wave profile in the $x$-direction and are periodic in the transverse direction $z$. These waves emanate from the Iooss-Kirchg\"{a}ssner line solitary waves in a \emph{dimension-breaking bifurcation} --- a phenomenon in which a spatially inhomogeneous solution of a partial differential equation emerges  from a solution which is homogeneous in one or more spatial dimensions (Haragus \& Kirchg\"{a}ssner \cite{HaragusKirchgaessner95}) --- and are sketched
in Figure \ref{dbp}. Periodically modulated solitary waves have also been computed numerically,
albeit for water of infinite depth, by Milewski \& Zhang \cite{MilewskiWang14}.

Our existence proof is based upon \emph{spatial dynamics}, a framework for studying stationary boundary-value problems by treating them as evolution equations in which an unbounded spatial variable plays the role of time. The method was pioneered by Kirchg\"assner \cite{Kirchgaessner82} and used extensively in the study of two-dimensional travelling water waves (see e.g.\ Groves \cite{Groves04} and the references therein). It was extended to the three-dimensional water-wave
setting by Groves \& Mielke \cite{GrovesMielke01} and has been used to construct a wide variety of three-dimensional gravity-capillary waves (see Groves \& Haragus \cite{GrovesHaragus03}).  
We use the method by formulating the hydrodynamic problem for travelling waves
as a reversible evolutionary equation
\begin{equation}
u_z=f(u), \label{Spatial dynamics formulation 1 in intro}
\end{equation}
where $u$ belongs to a space $\XX$ of
functions which converge to zero as $x \rightarrow \pm\infty$ and
$f: \DD(f) \subseteq \XX \rightarrow \XX$ is an unbounded, densely defined
vector field. Observe that
the nonzero equilibria of equation
\eqn{Spatial dynamics formulation 1 in intro} are
precisely the line solitary-wave solutions
of the water-wave problem; of particular interest in this respect are the
equilibria $u_\ve^\star$ corresponding to the
Iooss-Kirchg\"{a}ssner line solitary waves.
A solution to equation \eqn{Spatial dynamics formulation 1 in intro}
of the form
\begin{equation}
u(x,z)=u_\ve^\star(x)+u^\prime(x,z), \label{Basic Ansatz}
\end{equation}
where $u^\prime$ has small amplitude and is periodic in $z$, corresponds to a periodically modulated
solitary wave which is generated by a dimension-breaking bifurcation from
a Iooss-Kirchg\"{a}ssner line solitary wave. Accordingly, we construct a
family of small-amplitude periodic solutions to the equation
\begin{equation}
u^\prime_z = f(u_\ve^\star+u^\prime) \label{Apply LI to this}
\end{equation}
for each sufficiently small value of $\ve$. Explicit formulae for equations \eqn{Spatial dynamics formulation 1 in intro} and
\eqn{Apply LI to this} (together with the definitions of suitable function spaces in which
to study them) are presented in Section \ref{Spatial dynamics} below.

Small-amplitude periodic solutions of reversible evolutionary equations
are classically found using the Lyapunov centre theorem,
which asserts the existence of a family of these solutions
with frequency near $\omega_0$ provided that $\pm \ii \omega_0$ are nonresonant
imaginary eigenvalues of the corresponding linearised system.
The theorem remains true for infinite-dimensional
systems with finite-dimensional linear central subspaces under
additional hypotheses on the decay rate of the resolvent operator along the imaginary axis
(Devaney \cite{Devaney76b}). 
In the present setting the central subspace is however
infinite-dimensional due to the presence of essential spectrum at the origin
(a feature typical of spatial dynamics formulations for problems in unbounded
domains).
It was pointed out by Iooss \cite{Iooss99} that this difficulty
does not automatically rule out an application of the reversible
Lyapunov centre theorem. Denoting the linear and nonlinear parts of
the reversible vector field in question by respectively $L$ and $N$,
one finds upon examining the proof of the Lyapunov centre theorem that
it is not required that $L$ is invertible on the whole of phase space,
rather merely on the range of $N$ (`the Iooss condition at the origin').
We thus arrive at
the following generalisation of the reversible Lyapunov centre theorem,
which we refer to as the \emph{Lyapunov-Iooss theorem} (see also Bagri \& Groves
\cite{BagriGroves14} for an application to doubly periodic surface waves).\pagebreak

\begin{theorem} \label{LI theorem}
Consider the differential equation
\begin{equation}
v_\tau=L(v)+N(v), \label{Basic DE}
\end{equation}
in which $v(\tau)$ belongs to a Banach space $\XX$. Suppose that $\YY$, $\ZZ$ are further Banach spaces
with the properties that
\begin{list}{(\roman{count})}{\usecounter{count}}
\item
$\ZZ$ is continuously embedded in $\YY$ and continuously and densely embedded in $\XX$,
\item
$L: \ZZ \subseteq \XX \rightarrow \XX$ is a closed linear operator,
\item
there is an open neighbourhood $\UU$ of the origin in $\YY$
such that $L \in \LL(\YY,\XX)$ and $N \in C^3_{\mathrm{b,u}}(\UU, \XX)$ (and hence $N \in C^3_{\mathrm{b,u}}(\UU \cap \ZZ, \XX)$)
with $N(0)=0$, $\mathrm{d}N[0]=0$.
\end{list}
Suppose further that
\begin{list}{(\roman{count})}{\usecounter{count}}
\setcounter{count}{3}
\item equation \eqn{Basic DE} is reversible: there exists an involution $S \in \LL(\XX,\XX)$
with $SLv=-LSv$ and $SN(v)=-N(Sv)$ for all $v \in \UU$,
\end{list}
and that the following spectral hypotheses are satisfied.
\begin{list}{(\roman{count})}{\usecounter{count}}
\setcounter{count}{4}
\item $\pm \ii \omega_0$ are nonzero simple eigenvalues of $L$;
\item $\ii n \omega_0 \in \rho(L)$ for $n \in \Z\!\setminus\!\{-1,0,1\}$;
\item $\|(L-\ii n \omega_0 I )^{-1}\|_{\XX \rightarrow \XX} = o(1)$ and $\|(L-\ii n \omega_0 I )^{-1}\|_{\XX \rightarrow \ZZ} = O(1)$ as $n \rightarrow \infty$;
\item For each $v^\dag \in \UU$ the equation
$$Lv=-N(v^\dag)$$
has a unique solution $v \in \YY$ and the mapping $v^\dag \mapsto u$ belongs to $C^3_{\mathrm{b,u}}(\UU, \YY)$.
\end{list}
Under these hypotheses there exist an open neighbourhood $\NN$ of the origin in $\R$ and a\linebreak
continuously differentiable branch $\{(v(s),\omega(s))\}_{s \in S}$ of reversible, $2\pi/\omega(s)$-periodic solutions in
$C^1_\mathrm{per}({\mathbb R}, \YY \oplus \XX) \cap C_\mathrm{per}({\mathbb R}, \YY \oplus \ZZ)$
to \eqn{Basic DE} with amplitude $O(|s|)$.
Here the direct sum refers to the decomposition of a function into its mode $0$ and higher-order Fourier components,
the subscript `per' indicates a $2\pi/\omega(s)$-periodic function and $\omega(s)=\omega_0+O(|s|^2)$.
\end{theorem}

Our main task is now to study the purely imaginary spectrum of the linearisation $L$ of
the vector field on the right-hand side of \eqn{Apply LI to this} (see Sections
\ref{sec:derivation} and \ref{sec:analysis}). Motivated by the
weakly nonlinear scaling for modulated wave packets, we use the decomposition
$\eta=\eta_1+\eta_2$,
where the supports of the Fourier transforms of $\eta_1$ and $\eta_2$ lie respectively
in $S=[-\mu_0-\delta,-\mu_0+\delta] \cup [\mu_0-\delta,\mu_0+\delta]$ and $\R \setminus S$,
and write
$$\eta_1(x) = \frac{1}{2}\ve \zeta(\ve x) \ee^{\ii \mu_0 x} + \frac{1}{2}\ve \overline{\zeta(\ve x)} \ee^{-\ii \mu_0 x}.$$
We find, after a lengthy calculation,
that the pair $\pm \ii \ve k$, $k>0$, are eigenvalues of $L$ if and only if the operator
\begin{align*}
&\BB_{\ve, k}(\zeta, \psi)
:=\\
&
\begin{pmatrix}
A_2^{-1}\zeta-A_2^{-1}A_1\zeta_{xx}+k^2 \zeta-A_2^{-1}(A_3+A_5)\zeta^{\star 2}\zeta
-A_2^{-1}A_3\zeta^{\star 2}\overline{\zeta}
+4A_2^{-1}A_4\zeta^\star \psi_x-\RR_{\ve, k}(\zeta)\\[2mm]
-(1-\alpha_0^{-1})\psi_{xx}+ k^2 \psi-2A_4\re(\zeta^\star \zeta)_x
\end{pmatrix}
\end{align*}
has a zero eigenvalue, where $\RR_{\ve,k}(\zeta)$ is a remainder term
which is $O(\ve)$ in an appropriate sense; note that
$\BB_{0, k}(\zeta, \psi)$ is precisely the operator obtained
by linearising the Davey-Stewartson system at its line solitary wave
and taking the Fourier transform with respect to $z$.
A comparison with
well-studied operators in mathematical physics
shows that $\BB_{0,0}=\BB_{0,k}-k^2 I$ has precisely
one simple negative eigenvalue $-k_0^2$, so that $\BB_{0,k}$ has a
simple zero eigenvalue precisely when $k=k_0$, and we deduce from a spectral perturbation argument
that $\BB_{\ve,k}$ is invertible for $k \in [k_\mathrm{min},\infty) \setminus \{k_\ve\}$, where
$k_\mathrm{min}$ is a fixed small, positive number and $|k_0-k_\ve| = O(\ve)$, while
$\BB_{\ve,k_\ve}$ has a simple zero eigenvalue;
the corresponding eigenvector is invariant under
the symmetry $\tilde{R}:(\zeta(x), \psi(x)) \mapsto (\overline{\zeta(-x)},-\psi(-x))$,
thus corresponding to a linear water wave whose surface profile is symmetric in $x$. 

A similar argument is used to verify the Iooss condition at the origin. Applying the
above reduction to the equation
$$Lv=-N(v^\dag),$$
where $N$ is the nonlinear part of the vector field on the right-hand side of \eqn{Apply LI to this}, we obtain the reduced equation
$$\CC_\ve(\zeta) = \zeta^\dag,$$
where $\CC_0$ is the operator obtained by linearising the cubic nonlinear Schr\"{o}dinger equation at its solitary wave and $\CC_\ve-\CC_0$ is $O(\varepsilon)$ in an appropriate
sense. The operator $\CC_0$ is invertible in the class
of functions which are invariant under the symmetry $\zeta(x) \mapsto \overline{\zeta(-x)}$,
and it follows that the reduced equation is solvable in this class. We also show that
$\BB_{\ve,k}|_{\mathrm{Fix}\, \tilde{R}}$ has no zero eigenvalue in the interval $(0,k_\mathrm{min})$.
Working in function spaces corresponding to symmetric surface waves, we therefore obtain
the following existence result for periodically modulated solitary waves from
Theorem \ref{LI theorem}.

\begin{theorem}
There exist an open neighbourhood
${\mathcal N}_\epsilon$ of the origin in ${\mathbb R}$
and a family of periodically modulated solitary waves
$\{(\eta_s(x,z),\phi_s(x,y,z))\}_{s \in {\mathcal N}_\epsilon}$ which
emerges from the Iooss-Kirchg\"{a}ssner line solitary wave in
a dimension-breaking bifurcation. Here
\[\eta_s(x,z)=\eta_\ve^\star(x)+\eta^\star_s(x,z),\]
in which $\eta_s^\star(\cdot\,,\cdot)$ has amplitude $O(|s|)$ and is even in
both arguments and periodic in its second
with frequency $\varepsilon k_\varepsilon + O(|s|^2)$;
the positive number $k_\varepsilon$ satisfies $|k_\varepsilon-k_0| =O(\epsilon)$.
\end{theorem}

The spectral information necessary for the Lyapunov-Iooss theorem can also be used to study the stability properties of line solitary waves. The spatial dynamics formulation
\eqn{Spatial dynamics formulation 1 in intro} can be generalised to time-dependent waves; the result is an equation of the form
\begin{equation}
u_z = d u_t + f(u), \label{Spatial dynamics formulation 2 in intro}
\end{equation}
where $d: \DD(f) \subset \XX \rightarrow \XX$ is a linear operator (see Section
\ref{Spatial dynamics}).
We say that a line solitary wave $u^\star$ (that is, a solution of
\eqn{Spatial dynamics formulation 2 in intro} which is independent of $x$ and $t$)
is `transversely linearly unstable' if the linear equation
$$
u_z = d u_t + \mathrm{d}f[u^\star](u)
$$
has a solution of the form $\ee^{\sigma t} u_\sigma(z)$, where
$\re \sigma>0$ and $u_\sigma$ is bounded (see Section \ref{Spatial dynamics} for
a precise definition). The following theorem, which is
due to Godey \cite{Godey14}, shows that the existence of a pair of simple purely
imaginary eigenvalues of $L=\mathrm{d}f[u^\star]$ implies the transverse
linear instability of $u^\star$.

\begin{theorem} \label{thm:G theorem}
Consider the differential equation
\begin{equation}
v_\tau = d_1 v_t  + d_2 v_{tt}+ Lv , \label{Linearised multi-DE}
\end{equation}
in which $v(\tau,t)$ belongs to a Banach space $\XX$. Suppose that
$\ZZ$ is a further Banach space with the properties that
\begin{list}{(\roman{count})}{\usecounter{count}}
\item
$\ZZ$ is continuously and densely embedded in $\XX$,
\item
$L: \ZZ \subseteq \XX \rightarrow \XX$ and $d_j:\ZZ \subseteq \XX \rightarrow \XX$ are closed linear operators,
\item
the equation is reversible: there exists an involution $S \in \LL(\XX,\XX)$ with $LSv=-SLv$ and
$d_jSv=-Sd_jv$ for all $v \in \ZZ$,
\item
$L$ has a pair $\pm \ii \omega_0$ of isolated purely imaginary eigenvalues with odd algebraic multiplicity.
\end{list}
Under these hypotheses equation \eqn{Linearised multi-DE}
has a solution of the form $\ee^{\sigma t} v_\sigma(\tau)$, where
$v_\sigma \in C^1({\mathbb R}, \XX) \cap C(\R, \ZZ)$ is periodic,  for each sufficiently small positive
value of $\sigma$; its period tends to $2\pi/\omega_0$ as $\sigma \rightarrow 0$.
\end{theorem}

This result applies to the Iooss-Kirchg\"{a}ssner line solitary waves and 
confirms a prediction made by Ablowitz \& Segur \cite[p.\ 704]{AblowitzSegur79}
on the basis of the Davey-Stewartson approximation; note that it is not necessary
to restrict to symmetric surface waves to apply Theorem \ref{thm:G theorem}.

\begin{theorem}
The Iooss-Kirchg\"{a}ssner line solitary wave $(\eta^\star_{\ve},\phi^\star_{\ve})$
is transversely linearly unstable with respect to periodic perturbations with sufficiently
large period.
\end{theorem}

The stability properties of line solitary waves with weak surface tension with respect to two- dimensional perturbations were studied by Buffoni
\cite{Buffoni05,Buffoni09} and Groves \& Wahl\'{e}n \cite{GrovesWahlen10}.
They proved the conditional energetic stability of the set of minimisers of the energy subject to the constraint of fixed momentum. The precise relation between the set of minimisers and the Iooss- Kirchg\"{a}ssner waves still remains unclear, however.

The case of strong surface tension (that is, $\tau>1/3$) is rather different. The linear dispersion curve for a sinusoidal travelling wave train has its minimum at wavenumber $\mu=0$ (and $\alpha=1$), and the appropriate weakly nonlinear approximation is the
Kadomtsev-Petviashvili (KP-I) equation,
which reduces to the Korteweg-deVries
equation in the special case of $z$-independent waves. This
equation admits explicit line and periodically modulated solitary-wave solutions
(Tajiri \& Murakami \cite{TajiriMurakami90}); the line solitary waves are symmetric waves of
depression with monotonically and exponentially decaying tails.
An existence theory for line solitary waves in the parameter regime
$\alpha=1+\ve^2$ and $\beta=\beta_0$, where $\beta_0>1/3$, was given
by Kirchg\"{a}ssner \cite{Kirchgaessner88} and
Amick \& Kirchg\"{a}ssner \cite{AmickKirchgaessner89}, while their dimension breaking
to periodically modulated solitary waves was established by Groves, Haragus \& Sun
\cite{GrovesHaragusSun02} (see also Groves, Haragus \& Sun \cite{GrovesHaragusSun01}
and Pego \& Sun \cite{PegoSun04} for a discussion of their linear transverse instability
and Rousset and Tzvetkov \cite{RoussetTzvetkov11} for a proof of their nonlinear 
instability with respect to three-dimensional perturbations). The present treatment of the weak
surface tension case follows the strategy of Groves, Haragus \& Sun for strong
surface tension, and we indeed refer to their paper for the derivation of our spatial
dynamics formulation and some \emph{a priori} linear estimates. Our mathematics,
however, differs from theirs in several respects: we use simpler function spaces for the spatial dynamics
formulation and work directly with the variable $\eta$, an approach which leads to conceptually
simpler calculations and illuminates the connection with the weakly nonlinear model more clearly.

\section{Dimension-breaking phenomena}

\subsection{Spatial dynamics} \label{Spatial dynamics}

In this section we formulate equations \eqn{eq:Linear BVP}--\eqn{eq:Bernoulli condition} as an evolutionary equation in which the
unbounded horizontal direction $z$ plays the role of time.
The first step is to use the change of variable
\[
\phi(x,y,z,t)=\Phi(x,y^\prime,z,t), \qquad y=y^\prime(1+\eta(x,z,t))
\]
to map the variable fluid domain $\{(x,y,z): x, z \in \R, 0 < y < 1+ \eta(x,z,t)\}$ into the fixed strip $\{(x,y^\prime,z) \in \R \times (0,1) \times \R\}$;
we drop the primes for notational convenience.
Defining
\begin{align*}
\omega&=-\int_0^1 \left(\Phi_z-\frac{y\eta_z \Phi_y}{1+\eta}\right)y\Phi_y \dy+\frac{\beta \eta_z}{\sqrt{1+\eta_x^2+\eta_z^2}},\\
\xi&=\left(\Phi_z-\frac{y\eta_z \Phi_y}{1+\eta}\right)(1+\eta),
\end{align*}
one finds that equations \eqn{eq:Linear BVP}--\eqn{eq:Bernoulli condition} can be formulated as
\begin{equation}
u_z=du_t+f(u) \label{zt ES}
\end{equation}
for $u=(\eta, \omega, \Phi, \xi)$
with boundary conditions 
\begin{equation}
 \Phi_y=y \eta_t + B(u) \quad \text{on } y=0,1.
 \label{zt ES BC}
\end{equation}
Here
$$du=(0,\Phi|_{y=1},0,0),$$
$f(u)=(f_1(u), f_2(u), f_3(u), f_4(u))$ with
\begin{align*}
f_1(u)&=W\left(\frac{1+\eta_x^2}{\beta^2-W^2}\right)^{1/2},\\
f_2(u)&=\frac{W}{(1+\eta)^2}\left(\frac{1+\eta_x^2}{\beta^2-W^2}\right)^{1/2} \int_0^1 y\Phi_y\xi \dy -\left[\eta_x \left(\frac{\beta^2-W^2 }{1+\eta_x^2}\right)^{1/2} \right]_x+\alpha \eta -\Phi_x\big|_{y=1}\\
&\qquad\mbox{} + \int_0^1 \left\{ \frac{\xi^2-\Phi_y^{2}}{2(1+\eta)^2}+\frac12 \left(\Phi_x-\frac{y\eta_x \Phi_y}{1+\eta}\right)^2 +\left[\left(\Phi_x-\frac{y\eta_x \Phi_y}{1+\eta} \right)y\Phi_y\right]_x\right.\\
&\left. \qquad \qquad \qquad \mbox{} + \left(\Phi_x-\frac{y\eta_x \Phi_y}{1+\eta}\right)\frac{y\eta_x \Phi_y}{1+\eta}\right\} \dy,\\
f_3(u)&=\frac{\xi}{1+\eta}+\frac{Wy\Phi_y}{1+\eta} \left(\frac{1+\eta_x^2}{\beta^2-W^2}\right)^{1/2},\\
f_4(u)&=-\frac{\Phi_{yy}}{1+\eta}-\left[(1+\eta)\left(\Phi_x-\frac{y\eta_x \Phi_y}{1+\eta}\right)\right]_x +\left[\left(\Phi_x-\frac{y\eta_x \Phi_y}{1+\eta}\right)y\eta_x\right]_y\\
&\qquad \mbox{} +\frac{W (y\xi)_y}{1+\eta}\left(\frac{1+\eta_x^2}{\beta^2-W^2}\right)^{1/2},
\end{align*}
and
\[
B(u)=-y\eta_x+y\eta_x \Phi_x+\frac{\eta\Phi_y}{1+\eta}
-\frac{y^2 \eta_x^2 \Phi_y}{1+\eta}
+\frac{Wy\xi}{1+\eta} \left(\frac{1+\eta_x^2}{\beta^2-W^2}\right)^{1/2},
\]
where
$$W=\omega + \frac{1}{1+\eta}\int_0^1 y \Phi_y \xi \dy;$$
note the helpful identity
\begin{equation}
f_4(u) = -[(1+\eta)\Phi_x - y \eta_x\Phi_y]_x + [-\Phi_y+y\eta_x+B(\eta,\omega,\Phi,\xi)]_y.
\label{Helpful identity}
\end{equation}
In keeping with our study of the Iooss-Kirchg\"{a}ssner line solitary waves we
henceforth work in the parameter regime
$$\alpha=\alpha_0+\ve^2,\qquad \beta=\beta_0.$$

To identify an appropriate functional-analytic setting for these equations,
let us first specialise to stationary solutions, which satisfy
\begin{equation}
\label{z ES}
u_z=f(u)
\end{equation}
with boundary conditions 
\begin{equation}
\label{z ES BC}
 \Phi_y= B(u) \quad \text{on } y=0,1.
\end{equation}
Define
\[
X^s=H^{s+1}(\R)\times H^{s}(\R)\times H^{s+1}(\Sigma)\times H^{s}(\Sigma), \qquad s \geq 0,
\]
where $\Sigma$ is the strip $\R \times (0,1)$,
and let $M$ be an open neighbourhood of the origin in $X^1$ contained in the set 
$$
\{u \in X^1 : |W(x)|<\beta, \eta(x)>-1 \text{ for all } x\in \R\}
$$
(note that $M$ is a manifold domain of $X^0$).
The results presented by Bagri \& Groves \cite[Proposition 2.1]{BagriGroves14}) show that $B, f: M \rightarrow H^1(\Sigma)$
are analytic mappings. Equations \eqn{z ES}, \eqn{z ES BC} therefore constitute an evolutionary equation in the infinite-dimensional
phase space $X^0$ with
$$\DD(f) = \{u \in M: \Phi_y=B(u) \mbox{ on }y=0,1\}.$$
Furthermore, $f$ is the Hamiltonian vector field for the Hamiltonian system $(X^0, \Omega, H)$, where the position-independent 
symplectic form $\Omega$ on $X^0$ is given by
\[
\Omega(u_1, u_2)=\int_\R (\omega_2\eta_1-\eta_2\omega_1)\dx +
\int_\Sigma (\xi_2\Phi_1-\Phi_2\xi_1)\dy\dx
\]
and $H: M\to \R$ is the analytic function defined by
\begin{align*}
H(u)&=\int_{\Sigma}\left\{ (1+\eta)\Phi_x-y\eta_x \Phi_y+\frac{\xi^2-\Phi_y^{2}}{2(1+\eta)}-\frac{(1+\eta)}{2}\left(\Phi_x-\frac{y\eta_x\Phi_y}{1+\eta}\right)^2\right\}\dy\dx\\
&\qquad \mbox{}+\int_{\R} \left\{-\frac12 (\alpha_0+\ve^2) \eta^2 +\beta_0 -(\beta_0^2-W^2)^{1/2} (1+\eta_x^2)^{1/2}\right\} \dx
\end{align*}
(see Groves, Haragus \& Sun \cite[\S 2]{GrovesHaragusSun02}).
This Hamiltonian system is {\em reversible}, that is, Hamilton's equations are invariant under the transformation
\[
z\mapsto -z, \quad u\mapsto S(u),
\]
where the {\em reverser} $S: X^s\to X^s$ is defined by
$S(\eta, \omega, \Phi, \xi)=(\eta, -\omega, \Phi, -\xi)$.

Equation \eqn{z ES} is also invariant under the reflection
$R: X^s\to X^s$ given by
\[
R(\eta(x), \omega(x), \Phi(x,y), \xi(x,y))=(\eta(-x), \omega(-x), -\Phi(-x,y), -\xi(-x,y)),
\]
and one may seek solutions which are invariant under this symmetry by replacing $X^s$ by
\begin{align*}
X_\mathrm{r}^s & := X^s \cap \mathrm{Fix}\, R \\
& =
H_\mathrm{e}^{s+1}(\R)\times H_\mathrm{e}^s(\R)\times H_\mathrm{o}^{s+1}(\Sigma)\times H_\mathrm{o}^{s}(\Sigma),
\end{align*}
where
\begin{align*}
H_\mathrm{e}^s(\R)&=\{w\in H^s(\R): w(x)=w(-x) \text{ for all } x \in \R\},\\
H_\mathrm{o}^s(\Sigma)&=\{w\in H^s(\Sigma): w(x,y)=-w(-x,y) \text{ for all } (x,y)\in \Sigma\},
\end{align*}
(with corresponding definitions for $H_\mathrm{o}^s(\R)$ and $H_\mathrm{e}^s(\Sigma)$)
and $M$ by $M_\mathrm{r}:=M \cap \mathrm{Fix}\, R$.
It is also possible to extend $M_\mathrm{r}$
by replacing $X_\mathrm{r}^1$ in its definition with the extended function space $X_\star^1$, where
\[
X^s_\star=H_\mathrm{e}^{s+1}(\R)\times H_\mathrm{e}^s(\R)\times H_{\star,\mathrm{o}}^{s+1}(\Sigma)\times H_\mathrm{o}^{s}(\Sigma)
\]
and
$$H_{\star,\mathrm{o}}^{s+1}(\Sigma)=\{w\in L_\text{loc}^2(\Sigma): w_x, w_y \in H^s(\Sigma), w(x,y)=-w(-x,y) \text{ for all } (x,y)\in \Sigma\}
$$
(a Banach space with norm
$\|u\|_{\star, s+1}:=(\|u_x\|_s^2+\|u_y\|_s^2)^{1/2}$); note the relationship $M_\mathrm{r}=M_\star \cap X_\mathrm{r}^1$ between $M_\mathrm{r}$
and its extension $M_\star$.
This feature allows one to consider solutions to \eqn{z ES}, \eqn{z ES BC} whose $\Phi$-component is unbounded;
in particular line solitary waves fall into this category.

The line solitary waves found by Iooss \& Kirchg\"{a}ssner are equilibrium solutions of
equation \eqn{z ES}, \eqn{z ES BC} of the form $u^\star=(\eta^\star, 0, \Phi^\star,0) \in M_\star$, where $\eta^\star: \R \rightarrow \R$ and $\Phi^\star: \overline{\Sigma} \rightarrow \R$
are smooth functions whose derivatives are all $O(\ve)$ (see the Appendix), and we seek solutions of the form
\[
u=u^\star+u^\prime, \qquad u^\prime \in M^\prime,
\]
where $M^\prime$ is an open neighbourhood of the origin in $X^1$ chosen small enough so that $u^\star+M^\prime \subseteq M$;
note that we can here replace $M$ with $M_\mathrm{r}$ or $M_\star$.
Substituting this \emph{Ansatz} into \eqn{z ES}, \eqn{z ES BC} and again dropping
the prime for notational convenience, we find that
\begin{equation}
u_z=f(u^\star+u) \label{z ES perturbed}
\end{equation}
with boundary conditions
\begin{equation}
\Phi_y=B_\mathrm{l}(\eta, \Phi)+B_\mathrm{nl}(\eta, \omega, \Phi, \xi) \quad \text{on } y=0,1, \label{z ES BC perturbed}
\end{equation}
where
\begin{align*}
B_\mathrm{l}(\eta,\Phi)&=
\mathrm{d}B[\eta^\star,0,\Phi^\star,0](\eta,\omega,\Phi,\xi) \\
&=y(-\eta_x+\eta_x^\star \Phi_x+\Phi_x^\star \eta_x)\\
&\qquad\mbox{}+\frac{\eta^\star \Phi_y}{1+\eta^\star}
+\frac{\Phi_y^\star \eta }{(1+\eta^\star)^2}
+\frac{y^2(\eta_x^\star)^2\Phi_y^\star \eta}{(1+\eta^\star)^2}
-\frac{y^2(\eta_x^\star)^2 \Phi_y}{1+\eta^\star}
-\frac{2y^2\eta_x^\star \Phi_y^\star \eta_x}{1+\eta^\star}
\end{align*}
and
$$
B_\mathrm{nl}(\eta,\omega,\Phi,\xi)=B(\eta^\star+\eta, \omega, \Phi^\star+\Phi, \xi)
-B(\eta^\star,0, \Phi^\star,0)
-B_\mathrm{l}(\eta, \Phi).
$$
A small-amplitude periodic solution of \eqn{z ES perturbed}, \eqn{z ES BC perturbed} corresponds to a periodically modulated solitary wave
in the form of a three-dimensional perturbation of a line solitary wave.

Let us now return to \eqn{zt ES}, \eqn{zt ES BC}. Observing that $d \in \LL(X_0,X_0)$, one finds that these equations
constitute a reversible evolutionary equation with phase space $H^1((-T,T),X_0)$;
the domain of its vector field is
$$\{u \in H^2((-T,T),X_0) \cap H^1((-T,T), M): \Phi_y=y\eta_t+B(u) \mbox{ on }y=0,1\}$$
and its reverser is given by the pointwise extension of $S: X_0 \rightarrow X_0$ to $H^1((-T,T),X_0)$.
Seeking solutions of the form $u=u^\star+u^\prime$, where $u^\prime$ has pointwise values in
$M^\prime$, we find that
$$u_z = d u_t + f(u^\star+u)$$
with boundary conditions
$$
\Phi_y=y\eta_t+B_\mathrm{l}(\eta, \Phi)+B_\mathrm{nl}(\eta, \omega, \Phi, \xi) \quad \text{on } y=0,1,
$$
where we have again dropped the primes. We say that the solution $u^\star$ of \eqn{zt ES}, \eqn{zt ES BC}
is \emph{transversely linearly unstable} if the linear equation
\begin{equation}
u_z = d u_t + \mathrm{d}f[u^\star](u) \label{TS}
\end{equation}
with linear boundary conditions
\begin{equation}
\Phi_y=y\eta_t+B_\mathrm{l}(\eta, \Phi) \quad \text{on } y=0,1 \label{TS BC}
\end{equation}
has a solution of the form $\ee^{\sigma t} u_\sigma(z)$, with
$u_\sigma \in C_\mathrm{b}^1({\mathbb R}, X^0) \cap C_\mathrm{b}(\R, X^1)$ and
$\re \sigma>0$.

\subsection{Boundary conditions} \label{Boundary conditions}

Equations \eqn{z ES perturbed}, \eqn{z ES BC perturbed} cannot be studied using standard methods
for evolutionary equations because of the nonlinear boundary conditions. This difficulty is handled using a change of variable
which leads to an equivalent problem in a linear space. For $(\eta, \omega, \Phi, \xi)$ in $M$,
$M_\mathrm{r}$ or $M_\star$
define
$Q(\eta, \omega, \Phi, \xi)= (\eta, \omega,  \Gamma, \xi)$, 
where
\[
\Gamma=\Phi+\Theta_y
\]
and $\Theta \in H^3(\Sigma)$ is the unique solution of the boundary-value problem
\begin{align*}
-\Delta \Theta + B_\mathrm{l}(0, \Theta_y)&=B_\mathrm{nl}(\eta, \omega, \Phi, \xi) &&\text{in } \Sigma, \\
\Theta &=0 &&\text{on } y=0, 1
\end{align*}
(the linear operator on the left-hand side of this boundary-value problem is uniformly strongly elliptic with
smooth coefficients); note in particular that
\begin{equation}
-\Phi_y +B_\mathrm{l}(\eta,\Phi) + B_\mathrm{nl}(\eta,\omega,\Phi,\xi) = - \Gamma_y + B_\mathrm{l}(\eta,\Gamma) - \Theta_{xx}.
\label{Key feature of linear diffeo}
\end{equation}
The following result is established by the method used by  Groves, Haragus \& Sun \cite[Lemma 2.1]{GrovesHaragusSun02}.

\begin{lemma}
\hspace{1cm}
\begin{list}{(\roman{count})}{\usecounter{count}}
\item
The mapping $Q$ is a near-identity, analytic diffeomorphism from a 
neighbourhood $M$ of the origin in $X^1$ onto a neighbourhood $\tilde M$ 
of the origin in $X^1$.
\item
For each $u\in M$ the operator $\mathrm{d}Q[u]$
also defines an isomorphism
$\widehat{\mathrm{d}Q}[u]: X^0\to X^0$, and
the operators $\widehat{\mathrm{d}Q}[u]$, $\widehat{\mathrm{d}Q}[u]^{-1} \in  \LL(X^0, X^0)$ depend 
analytically on $u\in M$.
\item[(iii)] Statements (i) and (ii) also hold when $M$,  $\tilde{M}$ are replaced by
$M_\mathrm{r}$, $\tilde{M}_\mathrm{r}$ or $M_\star$, $\tilde{M}_\star$, where
$M_\mathrm{r}=M \cap \mathrm{Fix}\, R=M_\star \cap X_\mathrm{r}^1$ and
$\tilde{M}_\mathrm{r}=\tilde{M} \cap \mathrm{Fix}\, R=\tilde{M}_\star \cap X_\mathrm{r}^1$.
\end{list}
\end{lemma}

The above change of variable transforms \eqn{z ES perturbed} into the equation
\begin{equation}
v_z=L v+N(v) \label{z ES final}
\end{equation}
with linear boundary conditions
\begin{equation}
\Gamma_y=B_\mathrm{l}(\eta, \Gamma) \quad \text{on } y=0,1 \label{z ES BC final}
\end{equation}
for the variable $v=Q(u)$, in which $L:=\mathrm{d}\tilde{f}[0]=\mathrm{d}f[u^\star]$ and $N:=\tilde{f}-L$
are the linear and nonlinear parts of the transformed vector field\[
\tilde{f}(v):=\widehat{\mathrm{d}Q}[Q^{-1}(v)](f(u^\star+Q^{-1}(v))).
\]
Here
$$
\DD(L)=Y_1:=\left\{(\eta, \omega, \Gamma, \xi)\in X^1
: \Gamma_y=B_\mathrm{l}(\eta, \Gamma) \mbox{ on } y=0, 1\right\}$$
and
$\DD(N)$ is the neighbourhood $U := \tilde M \cap Y^1$
of the origin in $Y^1$, and we may replace $\tilde{M}$, $Y^1$, $U$
by $\tilde{M}_\mathrm{r}$, $Y_\mathrm{r}^1$, $U_\mathrm{r}$
or $\tilde{M}_\star$, $Y_\star^1$, $U_\star$ in the usual fashion.
The linear operator $L$ is given by the explicit formula
\[
L \begin{pmatrix} \eta \\ \omega \\ \Gamma \\ \xi \end{pmatrix}
=\begin{pmatrix}
\displaystyle \frac{\omega}{\beta_0}+h_1(\omega,\xi)\\
(\alpha_0+\ve^2)\eta-\Gamma_x|_{y=1}-\beta_0\eta_{xx}+h_2(\eta,\Gamma)\\
\xi+H_1(\omega,\xi)\\
-\Gamma_{xx}-\Gamma_{yy}+H_2(\eta,\Gamma),
\end{pmatrix},
\]
where
\begin{align*}
h_1(\omega,\xi)&=\frac{(1+(\eta_x^\star)^2)^{1/2}}{\beta_0}\left(\omega+\frac1{1+\eta^\star}\int_0^1 y\Phi_y^\star\xi \dy\right)-\frac{\omega}{\beta_0}, \\
h_2(\eta,\Gamma)&=\beta_0 \eta_{xx}-\beta_0 \left[\frac{\eta_x}{(1+(\eta_x^\star)^2)^{3/2}}\right]_x\\
&\hspace{-1cm} \mbox{}+
\int_0^1\!\!\left\{ \Phi_x^\star \Gamma_x-\frac{\Phi_y^\star\Gamma_y}{(1+\eta^\star)^2}+\frac{(\Phi_y^\star)^2\eta}{(1+\eta^\star)^3}
-\frac{y^2 (\eta_x^\star)^2\Phi_y^\star\Gamma_y}{(1+\eta^\star)^2}
-\frac{y^2 \eta_x^\star(\Phi_y^\star)^2\eta_x}{(1+\eta^\star)^2}
+\frac{y^2 \eta_x^\star(\Phi_y^\star)^2\eta}{(1+\eta^\star)^3}\right.
\\
&\hspace{-1cm} \qquad \qquad \left.\mbox{}
+\left[
y\Phi_y^\star\Gamma_x+y\Phi_x^\star\Gamma_y-
\frac{2y^2\eta_x^\star\Phi_y^\star\Gamma_y}{1+\eta^\star}
-\frac{y^2(\Phi_y^\star)^2\eta_x}{1+\eta^\star}
+\frac{y^2(\Phi_y^\star)^2\eta_x^\star\eta}{(1+\eta^\star)^2}\right]_x \right\} \dy,\\
H_1(\omega,\xi)&=-\frac{\eta^\star \xi}{1+\eta^\star}
+\frac{(1+(\eta_x^\star)^2)^{1/2}}{\beta_0(1+\eta^\star)}
\left(
\omega+\frac1{1+\eta^\star}\int_0^1 y\Phi_y^\star\xi \dy
\right)y\Phi_y^\star,
\\
H_2(\eta,\Gamma)&=(F_1(\eta,\Gamma))_x+(F_3(\eta,\Gamma))_y
\end{align*}
and
\begin{align*}
F_1(\eta,\Gamma)&=-\eta^\star\Gamma_x-\Phi_x^\star\eta +y\Phi_y^\star\eta_x+y\eta_x^\star \Gamma_y,\\
F_3(\eta,\Gamma)&=y\eta_x + B_\mathrm{l}(\eta,\Gamma).
\end{align*}
Furthermore, the change of variable preserves the reversibility of the evolutionary equation.

A similar issue arises when seeking solutions of \eqn{TS}, \eqn{TS BC} which depend analytically
upon $t$, which we regard as a parameter. In this setting the equations constitute an evolutionary system with
phase space $\tilde{X}^0$, the domain of whose vector field is the subset of functions in $\tilde{X}^1$ which satisfy \eqn{TS BC};
here the tilde denotes the analytic dependence upon $t$.
Standard methods cannot be applied because of the dependence of the boundary conditions upon the parameter, and
this problem is also resolved using a change of variables.
Define
\[
\Gamma=\Phi+\Theta_{yt},
\]
where $\Theta$ is the unique solution of the boundary-value problem
\begin{align*}
-\Delta \Theta + B_\mathrm{l}(0, \Theta_y)&=y\eta &&\text{in } \Sigma, \\
\Theta &=0 &&\text{on } y=0, 1.
\end{align*}
The change of variable $u=(\eta,\omega,\Phi,\xi) \mapsto (\eta,\omega,\Gamma,\xi)=v$ is a linear isomorphism
$\tilde{X}^0 \rightarrow \tilde{X}^0$ and
$\tilde{X}^1 \rightarrow \tilde{X}^1$. Note in particular that
$$
-\Phi_y +y \eta_t + B_\mathrm{l}(\eta,\Phi) = - \Gamma_y + B_\mathrm{l}(\eta,\Gamma) - \Theta_{xxt},
$$
so that this change of variable transforms \eqn{TS}, \eqn{TS BC} into
\begin{equation}
v_z = d_1 v_t + d_2 v_{tt} + Lv \label{TS final}
\end{equation}
with boundary conditions
\begin{equation}
\Gamma_y=B_\mathrm{l}(\eta, \Gamma) \quad \text{on } y=0,1, \label{TS BC final}
\end{equation}
where $d_1, d_2 \in \LL(\tilde{X}^0,\tilde{X}^0)$ are given by
$$
 d_1\begin{pmatrix} \eta \\ \omega \\ \Gamma \\ \xi \end{pmatrix}
=\begin{pmatrix}
0 \\
\Gamma|_{y=1}+\Theta_{xy}|_{y=1} - h_2(0,\Theta_y) \\
\tilde{\Theta}_y \\
-\eta -  (F_1(0,\Theta_y))_x \end{pmatrix}, \qquad
 d_2\begin{pmatrix} \eta \\ \omega \\ \Gamma \\ \xi \end{pmatrix}
=\begin{pmatrix}
0 \\ - \Theta_y|_{y=1} \\ 0 \\ 0
\end{pmatrix}
$$
and $\tilde{\Theta}$ is the unique solution of the boundary-value problem
\begin{align*}
-\Delta \tilde{\Theta} + B_\mathrm{l}(0, \tilde{\Theta}_y)&=y\left(\frac{\omega}{\beta_0} + h_1(\omega,\xi)\right) &&\text{in } \Sigma, \\
\tilde{\Theta} &=0 &&\text{on } y=0, 1.
\end{align*}

\subsection{Existence results}

We construct periodically modulated solitary waves by applying the Lyapunov-Iooss theorem
(Theorem \ref{LI theorem}) to \eqn{z ES final}, \eqn{z ES BC final}, taking $\XX=X_\mathrm{r}^0$, $\YY=Y_\star^1$, $\ZZ=Y_\mathrm{r}^1$ and $\UU=U_\star$; similarly, 
we demonstrate the transverse linear instability of the Iooss-Kirchg\"{a}ssner solitary waves
by applying Theorem \ref{thm:G theorem}
to \eqn{TS final}, \eqn{TS BC final}, taking
$\XX=X^0$ and $\ZZ=Y^1$, to construct a solution
of the form $\ee^{\sigma t} v_\sigma(z)$,
where $v_\sigma \in C^1(\R,X^0) \cap C(\R, Y^1)$ is
periodic. (In both cases we of course take $\tau=z$ and
$S(\eta,\omega,\Gamma,\xi) = (\eta,-\omega,\Gamma,-\xi)$.)
The spectral hypotheses are verified by studying the resolvent equations
\[
(L-\ii \lambda I)u=u^\dag, \qquad \lambda \in {\mathbb R},
\]
for $L$, that is
\begin{align}
\label{eq:spectral 1}
\frac{\omega}{\beta_0}+h_1(\omega,\xi) &=\ii \lambda \eta+\eta^\dag,\\
\label{eq:spectral 2}
(\alpha_0+\ve^2)\eta-\Gamma_x|_{y=1}-\beta_0\eta_{xx}+h_2(\eta,\Gamma) &=\ii \lambda \omega+\omega^\dag,\\
\label{eq:spectral 3}
\xi+H_1(\omega,\xi)&=\ii \lambda \Gamma+\Gamma^\dag,\\
\label{eq:spectral 4}
-\Gamma_{xx}-\Gamma_{yy}+(F_1(\eta,\Gamma))_x + (F_3(\eta,\Gamma))_y &=\ii \lambda \xi+\xi^\dag
\end{align}
with
\begin{align}
\Gamma_y&=0 && \text{on } y=0, \label{eq:spectral BC1}\\
\Gamma_y & = -\eta_x + F_3(\eta,\Gamma) && \text{on } y=1; \label{eq:spectral BC2}
\end{align}
since $L$ is real and anticommutes with the reverser $S$ it suffices to examine
non-negative values of $\lambda$, real values of $\eta$, $\Gamma$, $\omega^\dag$, $\xi^\dag$
and imaginary values of $\omega$, $\xi$, $\eta^\dag$, $\Gamma^\dag$.

The arguments given in Section \ref{Introduction} suggest that the support of
the Fourier transform $\hat{\eta}=\FF[\eta]$ of $\eta$ is concentrated near wavenumbers $\mu=\pm \mu_0$;
we therefore decompose it into the sum
of a function $\hat{\eta}_1$ with spectrum near $\mu=\pm \mu_0$
and a function $\hat{\eta}_2$ whose support is bounded away from these points. To this end
choose $\delta \in (0,\mu_0/3)$, let $\chi_0$ be the characteristic function of the interval $[-\delta,\delta]$
and define
$$
\eta_1=\chi(D)\eta:=\FF^{-1}[\chi \hat \eta], \qquad \eta_2=(1-\chi(D))\eta:=\FF^{-1}[(1-\chi)\hat \eta],
$$
where $\chi(\mu)=\chi_0(\mu-\mu_0)+\chi_0(\mu+\mu_0)$. In Section
\ref{sec:derivation} we reduce \eqn{eq:spectral 1}--\eqn{eq:spectral BC2} to a linear equation
of Davey-Stewartson type by determining $\eta_2$,
$\omega$, $\Gamma$ and $\xi$ as functions of $\eta_1$ and writing
$$\eta_1(x) = \frac{1}{2}\ve \zeta(\ve x) \ee^{\ii \mu_0 x} + \frac{1}{2}\ve \overline{\zeta(\ve x)} \ee^{-\ii \mu_0 x},$$
for some $\zeta \in \chi_0(\ve D)L^2(\R)$,
so that $\eta_1 \in \chi(D)L^2(\R)$. The results of these calculations are summarised in the
following theorem.

\begin{theorem}
\label{thm:reduction}
$ $
\begin{list}{(\roman{count})}{\usecounter{count}}
\item
There exists a constant $k_\mathrm{max}>0$ such that equations \eqn{eq:spectral 1}--\eqn{eq:spectral BC2} have a unique
solution $u \in Y^1$ for each $\lambda > \ve k_\mathrm{max}$ and each
$u^\dag \in X^0$.
\item
Suppose that $\lambda=\ve k$, where $0< k \le k_\mathrm{max}$. 
There exist $\zeta^\dag_{\ve, k}\in \LL(X^0, \chi_0(\ve D)L^2(\R))$ and an injection $\check u_{\ve, k} \in \LL(\chi_0(\ve D)L^2(\R)\times X^0, Y^1)$ such that
$u\in Y^1$ solves the resolvent equations
$$(L-\mathrm{i}\lambda I)u = u^\dag$$
if and only if $u=\check u_{\ve, k}(\zeta, u^\dag)$ for some solution
of the reduced equation 
\begin{align*}
\BB_{\ve, k}(\zeta, \psi)=(\zeta_{\ve, k}^\dag(u^\dag),0).
\end{align*}
The operator $\BB_{\ve, k}: \DD_\BB\subseteq W\to W$ is given by
\begin{align*}
&\BB_{\ve, k}(\zeta, \psi)
=\\
&
\begin{pmatrix}
A_2^{-1}\zeta-A_2^{-1}A_1\zeta_{xx}+k^2 \zeta-A_2^{-1}(A_3+A_5)\zeta^{\star 2}\zeta
-A_2^{-1}A_3\zeta^{\star 2}\overline{\zeta}
+4A_2^{-1}A_4\zeta^\star \psi_x-\RR_{\ve, k}(\zeta)\\[2mm]
-(1-\alpha_0^{-1})\psi_{xx}+ k^2 \psi-2A_4\re(\zeta^\star \zeta)_x
\end{pmatrix},
\end{align*}
where $W=L^2(\R) \times L^2(\R)$,
$\DD_\BB=H^2(\R) \times H^2(\R)$
and $\RR_{\ve, k} \in \LL(H^2(\R),L^2(\R))$ satisfies the estimate
\[
 \|\RR_{\ve, k}(\zeta)\|_0\le c\ve \|\zeta\|_1;
\]
each solution $(\zeta,\psi)$ of the reduced equation satisfies $\zeta \in \chi_0(\ve D)H^2(\R)$.
In particular, $L-\ii\ve k I$ is (semi-)Fredholm if $\BB_{\ve, k}$ is (semi-)Fredholm,
and the kernels of $L-\ii \ve k I$ and $\BB_{\ve, k}$ have the same dimension.

Furthermore, the reduction preserves the invariance of the resolvent equations under the
reflection $R$; its action on $W$ is given by $\tilde{R}:(\zeta(x), \psi(x)) \mapsto (\overline{\zeta(-x)},-\psi(-x))$.
\item
Define $Z_\mathrm{r}^0=H_\mathrm{e}^1(\R) \times L_\mathrm{e}^2(\R) \times H_\mathrm{o}^1(\Sigma) \times H_\mathrm{e}^{1}(\Sigma) \times \mathring{H}_\mathrm{o}^{1}(\Sigma)$, where
$$
\mathring{H}_\mathrm{o}^{1}(\Sigma)
=\{u \in H_\mathrm{o}^1(\Sigma): u|_{y=0}=u|_{y=1}=0\},$$
and
$$
H_\mathrm{c}^s(\R)= \{u \in H^s(\R): u(-x) = \overline{u(x)} \text{ for all } x \in \R\}, \qquad s \geq 0.
$$
There exist $\zeta^\dag_{\ve, 0}\in \LL(Z_\mathrm{r}^0, \chi_0(\ve D)L_\mathrm{c}^2(\R))$ and an injection $\check u_{\ve, 0} \in \LL(\chi_0(\ve D)L_\mathrm{c}^2(\R)\times Z_\mathrm{r}^0, Y_\star^1)$ such that
$u\in Y_\star^1$ solves
$$Lu=u^\dag,$$
where $\xi^\dag=-(\xi_1^\dag)_x-(\xi_3^\dag)_y$,
if and only if $u=\check u_{\ve, 0}(\zeta, \eta^\dag,\omega^\dag,\Gamma^\dag,\xi_1^\dag,\xi_3^\dag)$ for some
solution of the reduced equation
$$\CC_\ve \zeta = \zeta^\dag_{\ve,0}(\eta^\dag,\omega^\dag,\Gamma^\dag,\xi_1^\dag,\xi_3^\dag).$$
The operator $\CC_\ve: H_\mathrm{c}^2(\R) \subseteq L_\mathrm{c}^2(\R) \to L_\mathrm{c}^2(\R)$ is given by
$$
\CC_\ve \zeta =
A_2^{-1}\zeta-A_2^{-1}A_1\zeta_{xx}-2A_2^{-1}A_5\zeta^{\star 2}\zeta
-A_2^{-1}A_5\zeta^{\star 2}\overline{\zeta}
-\RR_{\ve,0}(\zeta)
$$
and $\RR_{\ve,0} \in \LL(H_\mathrm{c}^2(\R),L_\mathrm{c}^2(\R))$ satisfies the estimate
\[
 \|\RR_{\ve,0}(\zeta)\|_0\le c\ve \|\zeta\|_1;
\]
each solution of the reduced equation lies in $\chi_0(\ve D)H_\mathrm{c}^2(\R)$.
\end{list}
\end{theorem}

Part (i) of Theorem \ref{thm:reduction} implies that $\ii\lambda \in \rho(L)$ for $\lambda > \ve k_\mathrm{max}$,
while parts (ii) and (iii) show in particular that $\ii\ve k \in \rho(L)$ if $\BB_{\ve,k}$ (for $k > 0$)
or $\CC_\ve$ (for $k=0$) is invertible. We study the invertibility of these operators in
Section \ref{sec:analysis}, establishing the following results, whose corollaries relate them to
the operator $L$.

\begin{lemma}
For each sufficiently small $k_\mathrm{min}>0$ there
exists a unique number $k_\ve \in [k_\mathrm{min}, k_\mathrm{max}]$ with the following properties.
\begin{list}{(\roman{count})}{\usecounter{count}}
\item The operator $\BB_{\ve, k}: \DD_\BB \subseteq W \to W$ is an isomorphism for 
each $k\in [k_\mathrm{min}, k_\mathrm{max}]\setminus \{k_\ve\}$.
\item
The operator $\BB_{\ve, k_\ve}: \DD_\BB \subseteq W \to W$ is Fredholm with index 
$0$ and has a one-dimensional kernel which lies in $\mathrm{Fix}\, \tilde{R}$.
\end{list}
\end{lemma}
\begin{corollary}
The imaginary number $\ii\ve k$ belongs to
$\rho(L)$ for each $k \in [k_\mathrm{min},\infty) \setminus \{k_\ve\}$,
while $\ii\ve k_\ve$ is a simple
eigenvalue of $L$ whose eigenspace lies in $\mathrm{Fix}\, R$.
\end{corollary}
{\bf Proof.} It remains only to show that the eigenvalue $\ii \ve k_\ve$ is algebraically simple.
Observe that $\Omega(Lu_1,u_2)=-\Omega(u_1,Lu_2)$ and in particular that
$\Omega((L-\ii\ve k_\ve)u_1,u_2)=-\Omega(u_1,(L+\ii\ve k_\ve)u_2)$
for $u_1, u_2 \in \DD(L)$. It follows that $\Omega(f,\bar{u}_\ve)=0$ is a necessary
condition for $f \in X_\mathrm{r}^0$ to lie in the range of $L-\ii\ve k_\ve I$,
where $u_\varepsilon$ is an eigenvector corresponding to the eigenvalue $k_\ve$.
Using the formulae
\begin{align*}
\omega&=\frac{\ii \ve k_\ve \beta_0 \eta_\ve}{(1+(\eta_x^\star)^2)^{1/2}}
-\ii\ve k_\ve \int_0^1 y\Phi_y^\star\Gamma_\ve\dy + \frac{\ii\ve k_\ve \eta_\ve}{1+\eta^\star}\int_0^1 y\Phi_y^\star\dy, \\
\xi&=(1+\eta^\star)\ii \ve k_\ve\Gamma_\ve-\ii\ve k_\ve \eta_\ve
y\Phi_y^\star
\end{align*}
(see equations \eqn{eq:omega}, \eqn{eq:xi}), we find that
\[
\Omega(u_\ve, \bar{u}_\ve)=-2\ii \ve k_\ve \beta_0  \int_\R |\eta_\ve|^2\dx 
-2\ii \ve k_\ve \int_\Sigma |\Gamma_\ve|^2 \dy\dx+O(\ve^2(\|\eta_\ve\|_0^2+\|\Gamma_\ve\|_0^2))\ne 0,
\]
so that $u_\ve$ does not lie in the range of $L-\ii \ve k_\ve I$.\qed
\begin{lemma}  \label{lem:even analysis}
\hspace{1cm}
\begin{list}{(\roman{count})}{\usecounter{count}}
\item
The operator $\CC_\ve: H_\mathrm{c}^2(\R) \subseteq L_\mathrm{c}^2(\R) \to L_\mathrm{c}^2(\R)$ is an isomorphism.
\item
There exists $k_\mathrm{min}>0$ such that 
$\BB_{\ve,k}|_{\mathrm{Fix}\, \tilde{R}}$ is semi-Fredholm and injective for $k \in (0,k_\mathrm{min})$.
\end{list}
\end{lemma}
\begin{corollary}
The operator $L-\ii \ve k I:Y_\mathrm{r}^1 \to X_\mathrm{r}^0$ is an isomorphism for each $k\in (0,k_\mathrm{min})$.
\end{corollary}
{\bf Proof.} Because $k \mapsto \BB_{\ve,k}$ is a
continuous mapping $(0,k_\mathrm{max}) \mapsto \LL(\DD_\BB \cap \mathrm{Fix}\, \tilde{R}, \mathrm{Fix}\, \tilde{R})$
(see Section \ref{Analysis step 3})
and the Fredholm index of $\BB_{\ve,k}|_{\mathrm{Fix}\, \tilde{R}}$ is $0$ for $k \in [k_\mathrm{min},k_\mathrm{max}]$,
the same is true for $k \in (0,k_\mathrm{min})$. Using Lemma \ref{lem:even analysis}, we find that
$\BB_{\ve,k}|_{\mathrm{Fix}\, \tilde{R}}$ is invertible for $k \in (0,k_\mathrm{min})$,
and the stated result follows from Theorem \ref{thm:reduction}(ii).\qed

The resolvent decay estimate is a consequence of the next lemma, whose
proof is given by
Groves, Haragus \& Sun \cite[Lemma 3.25]{GrovesHaragusSun02}.
\begin{lemma}
There exists a constant $\lambda^\star>0$ such that the solution $u \in Y^1$
of equations \eqn{eq:spectral 1}--\eqn{eq:spectral BC2}
satisfies the estimate 
$$
\|u\|_{X^1}^2+ \lambda^2 \|u\|_{X^0}^2 \le c\|u^\dag\|_{X^0}^2
$$
for each $\lambda > \lambda^\star$ and each $u^\dag \in X_\mathrm{r}^0$.
\end{lemma}

It remains to establish the Iooss condition at the origin. For this purpose we first
record the following consequence of Theorem \ref{thm:reduction}(iii) and Lemma \ref{lem:even analysis}(i). 

\begin{proposition}
\label{prop:second spectral result}
For each $\eta^\dag \in H_\mathrm{e}^1(\R)$, $\omega^\dag \in L_\mathrm{e}^2(\R)$, $\Gamma^\dag \in H_\mathrm{o}^{s}(\Sigma)$
and 
\[
\xi^\dag=-(\xi_1^\dag)_x-(\xi_3^\dag)_y,
\]
where $\xi_1^\dag\in H_\mathrm{e}^{1}(\Sigma)$ and $\xi_3^\dag\in \mathring{H}_\mathrm{o}^{1}(\Sigma)$,
the equation
$$Lu = u^\dag$$
has a unique solution $u \in Y_\star^1$ which satisfies the estimate
\[
\|u\|_{Y_\star^1}\le c(\|\eta^\dag\|_{1}+\|\omega^\dag\|_0+\|\Gamma^\dag\|_{1}
+\|\xi_1^\dag\|_{1}+\|\xi_3^\dag\|_{1}).
\]
\end{proposition}
\begin{lemma} \label{Spectral result at origin}
For each $v^\dag \in U_\star$ the equation
\[
L v=-N(v^\dag)
\]
has a unique solution $v\in Y_\star^1$ which satisfies the estimate
\[
 \|v\|_{Y_\star^1}\le c\|v^\dag\|_{Y_\star^1}.
\]
\end{lemma}
{\bf Proof.} Choose $v^\dag \in U_\star$ and write $u^\dag=Q^{-1}(v^\dag)$.

Recall that
\[
g(v)=\widehat{\mathrm{d}Q}[u](f(u^\star+u)), \qquad u=Q^{-1}(v)
\]
and that $Q$, $\widehat{\mathrm{d}Q}[u]$ do not alter the fourth component of their arguments.
The fourth component of $g(v^\dag)$ is therefore given by
\begin{align*}
g_4(v^\dag) & =
f_4(u^\star+u^\dag) \\
&=-[(1+\eta^\star+\eta^\dag)(\Phi_x^\star+\Gamma_x^\dag)-y(\eta_x^\star+\eta_x^\dag)(\Phi_y^\star+\Gamma_y^\dag)]_x\\
&\qquad \qquad \mbox{}+
[-\Phi_y^\star-\Gamma_y^\dag+B(\eta^\star+\eta^\dag, \omega^\dag, \Phi^\star+\Gamma^\dag, \xi^\dag)+y(\eta_x^\star+\eta_x^\dag)]_y\\
&=-[(1+\eta^\star+\eta^\dag)\Gamma_x^\dag+\Phi_x^\star \eta^\dag-y(\eta_x^\star+\eta_x^\dag)\Gamma_y^\dag-y\Phi_y^\star\eta_x^\dag]_x\\
&\qquad \qquad \mbox{}+
[-\Gamma_y^\dag+B_\text{l}(\eta^\dag, \Gamma^\dag)+B_\mathrm{nl}(\eta^\dag,\omega^\dag,\Gamma^\dag,\xi^\dag)+y\eta_x^\dag]_y\\
&= -[(1+\eta^\star+\eta^\dag)\Gamma_x^\dag+\Phi_x^\star \eta^\dag-y(\eta_x^\star+\eta_x^\dag)\Gamma_y^\dag-y\Phi_y^\star\eta_x^\dag]_x\\
&\qquad \qquad \mbox{}+
[-\Gamma_y^\dag-\Theta_{xx}^\dag+B_\mathrm{l}(\eta^\dag,\Gamma^\dag)+y\eta^\dag_x]_y\\
&= -[\Gamma_x^\dag+(\eta^\star+\eta^\dag)\Gamma_x^\dag+\Phi_x^\star \eta^\dag-y(\eta_x^\star+\eta_x^\dag)\Gamma_y^\dag-y\Phi_y^\star\eta_x^\dag]_x\\
&\qquad \qquad \mbox{}+
[-\Gamma_y^\dag+B_\mathrm{l}(\eta^\dag, \Gamma^\dag)+y\eta^\dag_x]_y,
\end{align*}
in which the first line follows from \eqn{Helpful identity}, the second from the fact that $f_4(u^\star)=0$ and the definitions of $B_\mathrm{l}$, $B_\mathrm{nl}$,
the third from \eqn{Key feature of linear diffeo} and the fourth from the identity $\Gamma_{xx}^\dag = \Phi_{xx}^\dag + \Theta_{xxy}^\dag$.
According to the definition $N(v)=g(v)-Lv$ the fourth component of $N(v^\dag)$ is
\begin{align*}
N_{4}(v^\dag) &=
f_4(u^\star+u^\dag)-L_{4}v^\dag\\
&=
-[F_1(\eta^\dag,\Gamma^\dag)+(\eta^\star+\eta^\dag)\Gamma_x^\dag+\Phi_x^\star \eta^\dag-y(\eta_x^\star+\eta_x^\dag)\Gamma_y^\dag-y\Phi_y^\star\eta_x^\dag]_x\\
&\qquad \qquad \mbox{}+
[-F_3(\eta^\dag,\Gamma^\dag)+B_\mathrm{l}(\eta^\dag, \Gamma^\dag)+y\eta^\dag_x]_y\\
&= - [F_1(\eta^\dag,\Gamma^\dag)+(\eta^\star+\eta^\dag)\Gamma_x^\dag+\Phi_x^\star \eta^\dag-y(\eta_x^\star+\eta_x^\dag)\Gamma_y^\dag-y\Phi_y^\star\eta_x^\dag]_x.
\end{align*}

Applying Proposition \ref{prop:second spectral result} with
$$\xi_1^\dag = - \big(F_1(\eta^\dag,\Gamma^\dag)+(\eta^\star+\eta^\dag)\Gamma_x^\dag+\Phi_x^\star \eta^\dag-y(\eta_x^\star+\eta_x^\dag)\Gamma_y^\dag-y\Phi_y^\star\eta_x^\dag\big), \quad \xi_3^\dag=0,$$
one finds that the equation
$$Lv=-N(v^\dag)$$
has a unique solution $v \in Y_\star^1$ which satisfies
$$\|v\|_{Y_\star^1} \leq c(\|N_{1}(v^\dag)\|_1 + \|N_{2}(v^\dag)\|_0 + \|N_{3}(v^\dag)\|_1 + \|\xi_1^\dag\|_1).$$
The assertion follows by combining this estimate with
\[
\|\xi_1^\dag\|_1 \ \leq\  c(\|\eta^\dag\|_2+\|\nabla \Gamma^\dag\|_1 + \|\nabla \Gamma^\dag\|_1) \ \leq\  c\|v^\dag\|_{Y_\star^1}
\]
and
$$\|N(v^\dag)\|_{X^0} \leq c \|v^\dag\|_{Y_\star^1}.\eqno{\Box}$$

\section{Derivation of the reduced equation}
\label{sec:derivation}

\subsection{Reduction to a single pseudodifferential equation} \label{Derivation step 1}
We begin by eliminating $\omega$, $\Gamma$ and $\xi$ from \eqn{eq:spectral 1}--\eqn{eq:spectral BC2} and deriving a single, equivalent equation for $\eta$. Equations \eqn{eq:spectral 1} and \eqn{eq:spectral 3} clearly yield the explicit formulae
\begin{align}
\label{eq:omega}
\omega&=\frac{\beta_0}{(1+(\eta_x^\star)^2)^{1/2}}(\eta^\dag+\ii \lambda \eta)
-\frac{1}{1+\eta^\star}\int_0^1 y\Phi_y^\star\xi \dy,\\
\label{eq:xi}
\xi&=(1+\eta^\star)(\Gamma^\dag+\ii \lambda \Gamma)-y\Phi_y^\star(\eta^\dag+\ii \lambda \eta).
\end{align}
for $\omega$ and $\xi$, and substituting these formulae into \eqn{eq:spectral 4}--\eqn{eq:spectral BC2}, we obtain the  boundary-value problem
\begin{align}
-\hat \Gamma_{yy}+q^2\hat \Gamma &=\hat F^\dag(u^\dag)-\ii \mu \hat F_1(\eta,\Gamma)-\ii \lambda \hat F_2(\eta,\Gamma) -(\hat F_3(\eta,\Gamma))_y, && 0<y<1, \label{eq:Phi BC1} \\
\hat \Gamma_y&=0 && \text{on } y=0,  \label{eq:Phi BC2} \\
\hat \Gamma_y &=-\ii \mu \hat \eta+\hat F_3(\eta,\Gamma) && \text{on } y=1 \label{eq:Phi BC3}
\end{align}
for the Fourier transform $\hat{\Gamma}(\mu,y)$ of $\Gamma(x,y)$, where
$q=\sqrt{\mu^2+\lambda^2}$,
\begin{align*}
F_2(\eta,\Gamma) &=-\ii \lambda \eta^\star \Gamma+\ii \lambda y \Phi_y^\star \eta,\\
F^\dag(u^\dag) &=\xi^\dag+\ii \lambda(1+\eta^\star)\Gamma^\dag -\ii\lambda y\Phi_y^\star\eta^\dag
\end{align*}
and $\xi^\dag = - \xi_{1x}^\dag - \xi_{3y}^\dag$ for $\lambda=0$.
On the other hand, substituting \eqn{eq:omega}, \eqn{eq:xi} into 
\eqn{eq:spectral 2} yields
\begin{align*}
(\alpha_0+\ve^2)\eta&-\Gamma_x|_{y=1}-\beta_0\eta_{xx}+h_2(\eta,\Gamma)\\
&=\omega^\dag+\ii \lambda 
\frac{\beta_0}{(1+(\eta_x^\star)^2)^{1/2}}(\eta^\dag+\ii \lambda \eta)\\
&\qquad \mbox{}
-\frac{\ii \lambda}{1+\eta^\star}\int_0^1 y\Phi_y^\star
[(1+\eta^\star)(\Gamma^\dag+\ii \lambda \Gamma)-y\Phi_y^\star(\eta^\dag+\ii \lambda \eta)]
 \dy,
\end{align*}
which we write as
\begin{equation}
g_\ve(D, \lambda)\eta=\NN(\eta, \Gamma, u^\dag), \label{eq:Red eqn 1}
\end{equation}
where
$$
g_\ve(\mu,\lambda)=\alpha_0+\ve^2 + \beta_0 q^2 - \frac{\mu^2}{q^2}q\coth q
$$
and
\begin{align*}
\NN(&\eta, \Gamma, u^\dag)\\
&=\, \omega^\dag+\ii \lambda \beta_0 \frac{\eta^\dag}{(1+(\eta_x^\star)^2)^{1/2}}
+\Gamma_x|_{y=1}-\FF^{-1}\left[\frac{\mu^2}{q^2} q \coth q \,\hat \eta\right]
\\
&\qquad \mbox{}-\beta_0\left(\eta_{xx}-\left[\frac{\eta_x}{(1+(\eta_x^\star)^2)^{3/2}}\right]_x\right)
-\lambda^2\beta_0 \left(\frac1{(1+(\eta_x^\star)^2)^{1/2}}-1\right)\eta\\
&\qquad
\mbox{}-\frac{\ii\lambda}{1+\eta^\star}\int_0^1 [-\ii \lambda y^2(\Phi_y^\star)^2\eta - y^2(\Phi_y^\star)^2\eta^\dag
+\ii \lambda  y\Phi_y^\star \Gamma
+y\Phi_y^\star\Gamma^\dag
+\ii \lambda  y\Phi_y^\star\eta^\star \Gamma
+y\Phi_y^\star\eta^\star \Gamma^\dag]
\dy\\
&\qquad\mbox{}-\int_0^1\left\{ \Phi_x^\star \Gamma_x-\frac{\Phi_y^\star\Gamma_y}{(1+\eta^\star)^2}+\frac{(\Phi_y^\star)^2\eta}{(1+\eta^\star)^3}
-\frac{y^2 (\eta_x^\star)^2\Phi_y^\star\Gamma_y}{(1+\eta^\star)^2}
-\frac{y^2 \eta_x^\star (\Phi_y^\star)^2\eta_x}{(1+\eta^\star)^2}
+\frac{y^2 \eta_x^\star(\Phi_y^\star)^2\eta}{(1+\eta^\star)^3}\right.
\\
&\qquad \qquad \qquad\left.
+\left[
y\Phi_y^\star\Gamma_x+y\Phi_x^\star\Gamma_y-
\frac{2y^2\eta_x^\star\Phi_y^\star\Gamma_y}{1+\eta^\star}
-\frac{y^2(\Phi_y^\star)^2\eta_x}{1+\eta^\star}
+\frac{y^2(\Phi_y^\star)^2\eta_x^\star\eta}{(1+\eta^\star)^2}\right]_x \right\} \dy.
\end{align*}

The next step is to express $\Gamma$ as a function of $\eta$. To this end we formulate the boundary-value problem \eqn{eq:Phi BC1}--\eqn{eq:Phi BC3} as the integral equation
\begin{equation}
\Gamma = \GG_1(F_1(\eta,\Gamma),F_2(\eta,\Gamma),F_3(\eta,\Gamma),\eta)+\GG_2(F^\dag(u^\dag)), \label{eq:Phi integral}
\end{equation}
where
\begin{align*}
\GG_1(P_1,P_2,P_3,p)
&=\FF^{-1}\left[\int_0^1\{G(y, \tilde y)(-\ii \mu \hat P_1-\ii \lambda \hat P_2)+G_{\tilde y}(y,\tilde y) \hat P_3\} \dty
-\ii \mu G(y,1)\hat p\right], \\
\GG_2(P)&=\FF^{-1}\left[\int_0^1 G(y, \tilde y)\hat P \dty\right]
\end{align*}
and
\[
G(y, \tilde y)=\begin{cases}
\displaystyle{\frac{\cosh(qy)\cosh(q(1-\tilde y))}{q \sinh q}}, & 0\le y\le \tilde y \le 1,\vspace{5mm}\\
\displaystyle{\frac{\cosh(q\tilde y)\cosh(q(1-y))}{q \sinh q}}, &0\le \tilde y\le y\le 1.
\end{cases}
\]
The following proposition records the mapping properties of the integral operators
appearing on the right hand side of \eqn{eq:Phi integral}; it is proved using estimates for
$G(y,\tilde{y})$ similar to those given by Groves, Haragus \& Sun \cite[Proposition 3.2]{GrovesHaragusSun02}.

\begin{proposition} \label{prop:Green mapping}
$ $
\begin{list}{(\roman{count})}{\usecounter{count}}
\item
The mapping
$(P_1,P_2,P_3,p) \mapsto \GG_1(P_1,P_2,P_3,p)$
defines a linear function
$$
\begin{cases}
H^1(\Sigma)
\times H^1(\Sigma) \times H^1(\Sigma)
\times H^2(\R) \rightarrow H^2(\Sigma), & \lambda>0, \\
H_\mathrm{e}^1(\Sigma)
\times H_\mathrm{o}^1(\Sigma) \times H_\mathrm{o}^1(\Sigma)
\times H_\mathrm{e}^2(\R) \rightarrow
H_{\star,\mathrm{o}}^2(\Sigma), & \lambda=0,
\end{cases}
$$
which satisfies the estimate
\begin{align*}
\|\nabla \GG_1&\|_1 + \lambda \|\GG_1\|_1 + \lambda^2 \|\GG_1\|_0 \\
&\leq c (\|P_1\|_1 +\|P_2\|_1+ \|P_3\|_1 + \lambda(\|P_2\|_0+\|P_2\|_0+\|P_3\|_0)+(1+\lambda^2)\|p\|_0+\|p_{xx}\|_0).
\end{align*}
\item
Suppose that $\lambda>0$. The mapping $P \mapsto \GG_2(P)$ defines a linear function
$L^2(\Sigma) \rightarrow H^2(\Sigma)$ which satisfies the estimate
$$\|\nabla \GG_2\|_1 + \lambda \|\GG_2\|_1 + \lambda^2 \|\GG_2\|_0 \leq c(1+\lambda^{-1})\|P\|_0.$$
\item
Suppose that $\lambda=0$ and $P=-P_{1x} - P_{3y}$, where $P_1 \in H^1_\mathrm{e}(\Sigma)$
and $P_3 \in \mathring{H}^1_\mathrm{o}(\Sigma)$.
The mapping $(P_1,P_3) \mapsto \GG_2(P)$ defines a linear function
$H_\mathrm{e}^1(\Sigma) \times H_\mathrm{o}^1(\Sigma) \rightarrow H_{\star,\mathrm{o}}^2(\Sigma)$ which satisfies the estimate
$$\|\nabla \GG_2\|_1 \leq c(\|P_1\|_1 + \|P_3\|_1).$$
\end{list}
\end{proposition}
\begin{theorem}
\label{thm:Phi}
$ $
\begin{list}{(\roman{count})}{\usecounter{count}}
\item
Suppose that $\lambda>0$.
For each sufficiently small value of $\ve$ and each
$\eta\in H^{2}(\R)$,
$F^\dag \in L^2(\Sigma)$, equation \eqn{eq:Phi integral} has a unique solution
$\tilde{\Gamma}=\tilde{\Gamma}(\eta, u^\dag) \in H^{2}(\Sigma)$ 
which satisfies the estimate
$$
\|\nabla \tilde{\Gamma}\|_1 + \lambda \|\tilde{\Gamma}\|_1 + \lambda^2 \|\tilde{\Gamma}\|_0 \leq c ((1+\lambda^2)\|\eta\|_0 + \|\eta_{xx}\|_0 + (1+\lambda^{-1}) \|F^\dag\|_0).
$$
\item
Suppose that $\lambda=0$.
For each sufficiently small value of $\ve$ and each
$\eta\in H_\mathrm{e}^{2}(\R)$,
$\xi_1^\dag \in H^1_\mathrm{e}(\Sigma)$,
$\xi_3^\dag \in \mathring{H}^1_\mathrm{o}(\Sigma)$,
equation \eqn{eq:Phi integral} has a unique solution
$\tilde{\Gamma}=\tilde{\Gamma}(\eta, \xi_1^\dag,F_2^\dag) \in H_{\star,\mathrm{o}}^{2}(\Sigma)$ 
which satisfies the estimate
$$
\|\nabla \tilde{\Gamma}\|_1 \leq c (\|\eta\|_2+\|\xi_1^\dag\|_1+\|\xi_3^\dag\|_1).
$$
\end{list}
\end{theorem}
{\bf Proof.} It follows from the estimates
\begin{align*}
\|F_1\|_1,\ \|F_3\|_1 &\leq c\ve(\|\eta\|_1+\|\eta_x\|_1 + \|\nabla \Gamma\|_1) \\
&\leq c\ve(\|\eta\|_0+\|\eta_{xx}\|_0 + \|\nabla \Gamma\|_1), \\
\\
\lambda \|F_1\|_0,\ \lambda\|F_3\|_0 &\leq c\ve(\lambda\|\eta\|_0+\lambda\|\eta_x\|_0 + \lambda\|\nabla \Gamma\|_0) \\
&\leq c\ve((1+\lambda^2)\|\eta\|_0+\|\eta_{xx}\|_0 + \lambda\|\nabla \Gamma\|_0), \\
\\
\|F_2\|_1 &\leq c\ve(\lambda \|\eta\|_1 + \lambda \|\Gamma\|_1) \\
&\leq c\ve((1+\lambda^2)\|\eta\|_0+\|\eta_{xx}\|_0 + \lambda\|\Gamma\|_1), \\
\\
\lambda \|F_2\|_0 &\leq c(\lambda^2 \|\eta\|_0 + \lambda^2 \|\Gamma\|_0)
\end{align*}
and Proposition \ref{prop:Green mapping}(i) that
\begin{align*}
&\|\nabla\GG_1(\eta,\Gamma)\|_1 + \lambda \|\GG_1(\eta,\Gamma)\|_1 + \lambda^2 \|\GG_1(\eta,\Gamma)\|_0 \\
&\qquad\leq c(\ve(\|\nabla \Gamma\|_1 + \lambda \|\Gamma\|_1 + \lambda^2 \|\Gamma\|_0) + (1+\lambda^2)\|\eta\|_0 + \|\eta_{xx}\|_0),
\end{align*}
in which $\GG_1(\eta,\Gamma)$ is used an abbreviation for
$\GG_1(F_1(\eta,\Gamma),F_2(\eta,\Gamma),F_3(\eta,\Gamma),\eta)$. Using this estimate together with Proposition \ref{prop:Green mapping}(ii) or
(iii) as appropriate and equipping $H^2(\Sigma)$ with the norm $\Gamma \mapsto \|\nabla \Gamma\|_1 + \lambda \|\Gamma\|_1 + \lambda^2 \|\Phi\|_0$, one finds that
\eqn{eq:Phi integral} admits a unique solution which satisfies the stated estimate.\qed

Substituting $\Gamma=\tilde{\Gamma}(\eta, u^\dag)$ into \eqn{eq:Red eqn 1},
we obtain a single equation for $\eta$, namely
\begin{equation}
\label{eq:Red eqn 2}
g_\ve(D, \lambda)\eta= \tilde{\NN}(\eta, u^\dag),
\end{equation}
where
\[
\tilde{\NN}(\eta, u^\dag)=\NN(\eta, \tilde{\Gamma}(\eta,u^\dag),u^\dag).
\]
Finally, we record an estimate  for $g_\ve(\mu,\lambda)$
which is obtained by elementary calculus and use it to solve \eqn{eq:Red eqn 2}
for large values of $\lambda$, thus proving Theorem \ref{thm:reduction}(i).

\begin{proposition} \label{prop:g}
For each fixed $\lambda^\star>0$ the quantity $g_\ve(\mu,\lambda)$, which is even in both its arguments,
satisfies the estimates
\[
c((\mu-\mu_0)^2+\lambda^2+\ve^2) \le g_\ve (\mu, \lambda)\le \frac{1}{c}((\mu-\mu_0)^2+\lambda^2+\ve^2), \qquad (\mu,\lambda) \in S
\]
and
$$
c(1+q^2+\ve^2) \le g_\ve (\mu, \lambda)\le \frac{1}{c}(1+q^2+\ve^2), \qquad (\mu,\lambda) \in [0,\infty)^2 \setminus S,
$$
where $S=[\mu_0-\delta,\mu_0+\delta] \times [0,\lambda^\star]$.
\end{proposition}
\begin{lemma}
Choose $\lambda^\star>0$. For each sufficiently small value of $\ve$, each
$\lambda>\lambda^\star$ and each $u^\dag\in X^0$,
equation \eqn{eq:Red eqn 2} has a unique solution $\eta=\eta(u^\dag) \in H^2(\R)$ which satisfies the estimate
$$(1+\lambda^2)\|\eta\|_0 + \|\eta_{xx}\|_0 \leq c_\lambda \|u^\dag\|_{X^0}.$$
\end{lemma}
{\bf Proof.} Write \eqn{eq:Red eqn 2} as
$$g_\ve(D,\lambda)\eta - \tilde\NN(\eta,0) = \tilde\NN(0,u^\dag)$$
and recall the estimates $g_\ve(\mu,\lambda) \geq c (1+\mu^2+\lambda^2)$ (see
Proposition \ref{prop:g}) and
$$\|\tilde{\NN}(\eta,0)\|_0 \leq c\ve((1+\lambda^2)\|\eta\|_0+\|\eta_{xx}\|_0).$$
Equipping $H^2(\R)$ with the norm $\eta \mapsto (1+\lambda^2)\|\eta\|_0 + \|\eta_{xx}\|_0$, one finds that the operator
$$g_\ve(D,\lambda) - \tilde{\NN}(\cdot,0):L^2_\mathrm{e}(\R) \rightarrow H^2_\mathrm{e}(\R)$$
is invertible for each sufficiently small value of $\ve$,
and the result follows from this fact and the estimate
$$
\|\tilde{\NN}(0,u^\dag)\|_0 \leq c_\lambda \|u^\dag\|_{X^0}.\eqno{\Box}
$$

In view of the previous lemma we now fix $\lambda^\star>0$ and henceforth suppose that $\lambda\leq\lambda^\star$.
We write \eqn{eq:Red eqn 2} as
\begin{equation}
\label{eq:eta1}
g_\ve(\mu,\lambda)\hat \eta_1 =
\chi \FF[\tilde{\NN}(\eta_1+\eta_2,u^\dag)]
\end{equation}
and
\begin{equation}
\label{eq:eta2}
g_\ve(\mu,\lambda)\hat \eta_2 =
(1-\chi) \FF[\tilde{\NN}(\eta_1+\eta_2,u^\dag)],
\end{equation}
where $\hat{\eta}_1 = \chi \hat{\eta}$ and $\hat{\eta}_2 = (1-\chi) \hat{\eta}$. 

\begin{lemma} \label{lem:invert}
The linear operator
$$g_\ve(D,\lambda) - (1-\chi(D))\tilde \NN(\cdot,0): (1-\chi(D))H^2(\R) \rightarrow (1-\chi(D))L^2(\R)$$ is invertible
for each sufficiently small value of $\ve$.
\end{lemma}
{\bf Proof.} This result follows from the estimates $g_\ve(\mu,\lambda) \geq c (1+q^2)$ for $||\mu|-\mu_0| > \delta$ (see
Proposition \ref{prop:g}) and
$$\|\tilde\NN(\eta,0)\|_0 \leq c\ve\|\eta\|_2.\eqno{\Box}$$
\begin{corollary}
$ $
\begin{list}{(\roman{count})}{\usecounter{count}}
\item
Suppose that $\lambda>0$.
For each sufficiently small value of $\ve$ and each $\eta_1 \in \chi(D)L^2(\R)$, $u^\dag\in X^0$,
equation \eqn{eq:eta2} has a unique solution $\eta_2=\check \eta_2(\eta_1, u^\dag)\in H^2(\R)$ which satisfies the estimate
$$
\|\check \eta_2\|_{2}
\le
c\ve\|\eta_1\|_0 + c_\lambda \|u^\dag\|_{X^0}.
$$
\item
Suppose that $\lambda=0$.
For each sufficiently small value of $\ve$ and each $\eta_1 \in \chi(D)L^2_\mathrm{e}(\R)$, $u^\dag\in X_\mathrm{r}^0$,
equation \eqn{eq:eta2} has a unique solution $\eta_2=\check \eta_2(\eta_1, u^\dag)\in H_\mathrm{e}^2(\R)$ which satisfies the estimate
$$
\|\check \eta_2\|_{2}
\le
c(\ve\|\eta_1\|_0 + \|\omega^\dag\|_0 + \|\xi_1^\dag\|_1 + \|\xi_3^\dag\|_1).
$$
\end{list}
\end{corollary}
{\bf Proof.}
Write \eqn{eq:eta2} as
$$g_\ve(D,\lambda)(\eta_2) - (1-\chi(D))\tilde \NN(\eta_2,0)=(1-\chi(D))\big(\tilde{\NN}(\eta_1,0)+\tilde{\NN}(0,u^\dag)\big)$$
and use Lemma \ref{lem:invert} and the estimates
$$\|\tilde\NN(\eta_1,0)\|_0 \leq c\ve\|\eta_1\|_0, \qquad
\|\tilde\NN(0,u^\dag)\|_0 \leq \begin{cases}
c_\lambda \|u^\dag\|_{X^0}, & \lambda>0, \\
c(\|\omega^\dag\|_0 + \|\xi_1^\dag\|_1 + \|\xi_3^\dag\|_1), & \lambda=0.
\end{cases}
$$
\qed

Substituting $\check \eta_2=\check \eta_2(\eta_1,u^\dag)$ into \eqn{eq:eta1},
we obtain a single equation for $\eta_1$, namely
\begin{equation}
\label{eq:Red eqn 3}
g_\ve(\mu,\lambda)\hat \eta_1 =
\chi \FF [\check{\NN}(\eta_1,u^\dag)],
\end{equation}
where
\begin{align*}
\check{\NN}(\eta_1,u^\dag)&=\tilde{\NN}(\eta_1+\check \eta_2(\eta_1,u^\dag), u^\dag)\\
&=\NN(\eta_1+\check{\eta}_2(\eta_1,u^\dag), \check \Gamma(\eta_1, u^\dag), u^\dag)
\end{align*}
and
$$
\check \Gamma(\eta_1, u^\dag)=\tilde{\Gamma}(\eta_1+\check \eta_2(\eta_1, u^\dag), u^\dag).
$$

\subsection{Solution of the reduced equation for intermediate spectral values}
\label{Derivation step 2}
The next step is to derive more precise estimates for the functions introduced in the previous section. We
expand $\eta^\star$ and $\Phi^\star$ as
$$
\eta^\star(x)=\eta^\star_1(x)+\eta^\star_2(x)+\eta^\star_\mathrm{r}(x),\qquad
\Phi^\star(x,y)=\Phi^\star_1(x,y)+\Phi^\star_2(x,y)+ \Phi^\star_\mathrm{r}(x,y)
$$
and $\tilde{\Gamma}=\tilde{\Gamma}(\eta, u^\dag)$ as
\[
\tilde{\Gamma}(\eta,u^\dag)=\tilde{\Gamma}^{(1)}(\eta)+\tilde{\Gamma}^{(2)}(\eta)+ \tilde{\Gamma}^{(3)}(\eta)+ \tilde{\Gamma}^\mathrm{r}(\eta)
+\tilde \Gamma^\dag(u^\dag).
\] 
Here $\eta^\star_1$, $\eta^\star_2$, $\eta^\star_\mathrm{r}$ and
$\Phi^\star_1$, $\Phi^\star_2$, $\Phi^\star_\mathrm{r}$ are given in the Appendix
(equations \eqn{Defn eta1star}--\eqn{Defn Phi2star}),
\begin{align}
\label{eq:Phi1}
\tilde{\Gamma}^{(1)}(\eta)&=\GG_1(0,0,\eta),\\
\label{eq:Phi2}
\tilde{\Gamma}^{(2)}(\eta)&=\GG_1(F_1^{(2)}(\eta),F_2^{(2)}(\eta),F_3^{(2)}(\eta),0), \\
\label{eq:Phi3}
\tilde{\Gamma}^{(3)}(\eta)&=\GG_1(F_1^{(3)}(\eta),F_2^{(3)}(\eta),F_3^{(3)}(\eta),0)
\end{align}
with
\begin{align*}
F_1^{(2)}(\eta)&=-\eta^\star_1\tilde{\Gamma}_x^{(1)}(\eta)-\Phi_{1x}^\star\eta +y\Phi_{1y}^\star\eta_x+y\eta_{1x}^\star \tilde{\Gamma}_y^{(1)}(\eta), \\
F_2^{(2)}(\eta)&=-\ii \lambda \eta^\star_1 \tilde{\Gamma}^{(1)}(\eta)+\ii \lambda y \Phi_{1y}^\star \eta, \\
F_3^{(2)}(\eta)&=
y\eta_{1x}^\star\tilde{\Gamma}_x^{(1)}(\eta)+y\Phi_{1x}^\star\eta_x
+\eta^\star_1 \tilde{\Gamma}_y^{(1)}(\eta)+\Phi^\star_{1y} \eta, \\
\\
F_1^{(3)}(\eta)&=-\eta^\star_1\tilde{\Gamma}_x^{(2)}(\eta)-\eta^\star_2\tilde{\Gamma}_x^{(1)}(\eta)
-\Phi_{2x}^\star\eta+y\Phi_{2y}^\star\eta_x+y\eta_{1x}^\star\tilde  \Phi_y^{(2)}(\eta)+y\eta_{2x}^\star \tilde{\Gamma}_y^{(1)}(\eta),
\\
F_2^{(3)}(\eta)&=-\ii \lambda \eta_1^\star\tilde{\Gamma}^{(2)}(\eta)
-\ii \lambda \eta^\star_2 \tilde{\Gamma}^{(1)}(\eta)+\ii \lambda y \Phi_{2y}^\star \eta,
\\
F_3^{(3)}(\eta)&=
y\eta_{1x}^\star \tilde{\Gamma}_x^{(2)}(\eta)
+
y\eta_{2x}^\star\tilde{\Gamma}_x^{(1)}(\eta)
+y\Phi_{2x}^\star\eta_x
+\eta_1^\star\tilde{\Gamma}_y^{(2)}(\eta)
+\eta^\star_2 \tilde{\Gamma}_y^{(1)}(\eta)
+\Phi^\star_{2y} \eta
\\
&\qquad\mbox{}
-(\eta_1^\star)^2\tilde{\Gamma}_y^{(1)}(\eta)
-2\Phi^\star_{1y}\eta_1^\star \eta
-y^2(\eta_{1x}^\star)^2\tilde{\Gamma}_y^{(1)}(\eta)-2y^2\eta_{1x}^\star\Phi_{1y}^\star\eta_x,
\end{align*}
and $\tilde{\Gamma}^\mathrm{r}(\eta)$, $\tilde \Gamma^\dag(u^\dag)$ are determined implicitly by
\begin{equation}
\label{eq:Phir}
\tilde{\Gamma}^\mathrm{r}(\eta)=\GG_1(F_1^\mathrm{r}(\eta),F_2^\mathrm{r}(\eta),F_3^\mathrm{r}(\eta),0)
\end{equation}
with
$$
F_j^\mathrm{r}(\eta)=F_j(\eta,\tilde{\Gamma}^{(1)}(\eta)+\tilde{\Gamma}^{(2)}(\eta)+\tilde{\Gamma}^{(3)}(\eta)+\tilde{\Gamma}^\mathrm{r}(\eta))-F_j^{(2)}(\eta)-F_j^{(3)}(\eta), \qquad j=1,2,3,
$$
and
$$
\tilde \Gamma^\dag=\GG_1(F_1(0, \tilde \Gamma^\dag),F_2(0, \tilde \Gamma^\dag),F_3(0, \tilde \Gamma^\dag),0)+\GG_2(F^\dag(u^\dag)).
$$

\begin{lemma} \label{lemma:tildePhi estimates}
The functions $\tilde{\Gamma}^{(1)}$, $\tilde{\Gamma}^{(2)}$, $\tilde{\Gamma}^{(3)}$, $\tilde{\Gamma}^\mathrm{r}$ 
and $\tilde{\Gamma}^\dag$
satisfy the estimates
\begin{align*}
\|\nabla \tilde{\Gamma}^{(j)}(\eta)\|_{1} + \lambda \|\tilde{\Gamma}^{(j)}(\eta)\|_{1} + \lambda^2 \|\tilde{\Gamma}^{(j)}(\eta)\|_{0} &\le c\ve^{j-1} \|\eta\|_2,\qquad j=1,2,3,\\
\|\nabla \tilde{\Gamma}^\mathrm{r}(\eta)\|_{1} + \lambda \|\tilde{\Gamma}^\mathrm{r}(\eta)\|_{1} + \lambda^2 \|\tilde{\Gamma}^\mathrm{r}(\eta)\|_{0} &\le c\ve^3 \|\eta\|_2, \\
\|\nabla \tilde{\Gamma}^{\dag}\|_1+\lambda \| \tilde{\Gamma}^{\dag}\|_1 +\lambda^2 \| \tilde{\Gamma}^{\dag}\|_0 
& \le 
\begin{cases}
c_\lambda \|u^\dag\|_{X^0}, & \lambda>0, \\
c(\|\omega^\dag\|_0+\|\xi_1^\dag\|_1 + \|\xi_3^\dag\|_1), & \lambda=0.
\end{cases}
\end{align*}
\end{lemma}

We correspondingly write
\[
\tilde{\NN}(\eta, u^\dag)=\tilde{\NN}^{(2)}(\eta)+\tilde{\NN}^{(3)}(\eta)+\tilde{\NN}^\mathrm{r}(\eta)+\tilde{\NN}^\dag(u^\dag),
\]
where
\begin{align}
\tilde{\NN}^{(2)}(\eta)&=\tilde{\Gamma}_x^{(2)}(\eta)\Big|_{y=1}-\ii\lambda\int_0^1 y\Phi_{1y}^\star
\ii \lambda \tilde{\Gamma}^{(1)}(\eta)\dy \nonumber \\
&\qquad \mbox{}-\int_0^1\left\{ \Phi_{1x}^\star \tilde{\Gamma}_x^{(1)}(\eta)-\Phi_{1y}^\star \tilde{\Gamma}_y^{(1)}(\eta)+\left[
y\Phi_{1y}^\star\tilde{\Gamma}_x^{(1)}(\eta)+y\Phi_{1x}^\star\tilde{\Gamma}_y^{(1)}(\eta)\right]_x \right\} \dy,
\label{Eqn for tildeNN2}\\
\nonumber \\
\tilde{\NN}^{(3)}(\eta)&=\tilde{\Gamma}_x^{(3)}(\eta)\Big|_{y=1}
-\frac{3\beta_0}{2} \left[(\eta_{1x}^\star)^ 2 \eta_x\right]_x
+\frac{\lambda^2 \beta_0}{2} \eta_{1x}^{\star 2}\eta \nonumber \\
&
\qquad \mbox{}-\ii\lambda\int_0^1\left\{-\ii \lambda y^2(\Phi_{1y}^\star)^2\eta
 + \ii \lambda y\Phi_{1y}^\star  \tilde{\Gamma}^{(2)}(\eta)
 - \ii \lambda y\Phi_{2y}^\star \tilde{\Gamma}^{(1)}(\eta)
 \right\}\dy \nonumber \\
&\qquad \mbox{}-\int_0^1\Bigl\{ \Phi_{1x}^\star  \tilde{\Gamma}_x^{(2)}(\eta)-\Phi_{1y}^\star \tilde{\Gamma}_y^{(2)}(\eta)
+2\Phi_{1y}^\star \eta^\star_1 \tilde{\Gamma}_y^{(1)} (\eta)
+(\Phi_{1y}^\star)^2\eta \} \dy \nonumber \\
&\qquad \mbox{}-\int_0^1
\left[
y\Phi_{1y}^\star \tilde{\Gamma}_x^{(2)}(\eta)+y\Phi_{1x}^\star \tilde{\Gamma}_y^{(2)}(\eta)
-2y^2\eta_{1x}^\star\Phi_{1y}^\star \tilde{\Gamma}_y^{(1)}(\eta)
-y^2(\Phi_{1y}^\star)^2 \eta_x
\right]_x  \dy \nonumber \\
&\qquad-\int_0^1\left\{ \Phi_{2x}^\star \tilde{\Gamma}_x^{(1)}(\eta)-\Phi_{2y}^\star \tilde{\Gamma}_y^{(1)}(\eta)+\left[
y\Phi_{2y}^\star\tilde{\Gamma}_x^{(1)}(\eta)+y\Phi_{2x}^\star\tilde{\Gamma}_y^{(1)}(\eta)\right]_x \right\} \dy,\
\label{Eqn for tildeNN3} \\
\nonumber \\
\tilde \NN(\eta)&=\tilde \NN(\eta,0)-\tilde \NN^{(2)}(\eta)-\tilde \NN^{(3)}(\eta), \nonumber \\
\nonumber \\
\tilde{\NN}^\dag(u^\dag)&= \tilde\NN(0,u^\dag). \nonumber
\end{align}

\begin{lemma} \label{lemma:tildeNestimates}
The functions $\tilde \NN^{(2)}$, $\tilde \NN^{(3)}$, $\tilde \NN^\mathrm{r}$ and
$\tilde \NN^{\dag}$ satisfy the estimates
\begin{align*}
\| \tilde{\NN}^{(j)}(\eta)\|_0 &\le c\ve^{j-1} \|\eta\|_2,\qquad j=2,3, \\
\| \tilde{\NN}^\mathrm{r}(\eta)\|_0 &\le c\ve^3 \|\eta\|_2, \\
\| \tilde{\NN}^{\dag}(u^\dag)\|_0
& \le \begin{cases}
c_\lambda \|u^\dag\|_{X^0}, & \lambda>0, \\
c(\|\omega^\dag\|_0+\|\xi_1^\dag\|_1 + \|\xi_3^\dag\|_1), & \lambda=0.
\end{cases}
\end{align*}
\end{lemma}

Similarly, we expand $\check \eta_2=\check \eta_2 (\eta_1, u^\dag)$ as
$$
\check \eta_2 (\eta_1, u^\dag)
= \check{\eta}_2^{(2)}(\eta_1) + \check{\eta}_2^\mathrm{r}(\eta_1)+\check \eta_2^\dag(u^\dag),$$
where
\[
\FF[ \check \eta_2^{(2)}]=
g_\ve(\mu,\lambda)^{-1}(1-\chi)\FF[
\tilde{\NN}^{(2)}(\eta_1)]
\]
and
$\check \eta_2^\mathrm{r}$, $\check \eta_2^\dag$ are determined implicitly by
\begin{align*}
\FF[\check \eta_2^\mathrm{r}] 
&=g_\ve(\mu,\lambda)^{-1}(1-\chi)\FF[
\tilde{\NN}^{(3)}(\eta_1)+\tilde{\NN}^\mathrm{r}(\eta_1)+
\tilde{\NN}(\check{\eta}_2^{(2)}+\check \eta_2^\mathrm{r},0)], \\
\FF[\check \eta_2^\dag] 
&=g_\ve(\mu,\lambda)^{-1}(1-\chi)\FF[\tilde{\NN}(\check \eta_2^\dag,u^\dag)],
\end{align*}
and $\check \Gamma = \check{\Gamma}(\eta_1,u^\dag)$ as
\[
\check \Gamma(\eta_1,u^\dag)=\check \Gamma^{(1)}(\eta_1)+\check \Gamma^{(2)}(\eta_1)+\check \Gamma^{(3)}(\eta_1)+\check \Gamma^{\mathrm{r}}(\eta_1)+
\check \Gamma^{\dag}(u^\dag),
\]
where
\begin{align*}
\check \Gamma^{(1)}&=\tilde{\Gamma}^{(1)}(\eta_1),\\
\check \Gamma^{(2)}&=\tilde{\Gamma}^{(2)}(\eta_1)+\tilde{\Gamma}^{(1)}(\check \eta_2^{(2)}(\eta_1)),\\
\check \Gamma^{(3)}&=
\tilde{\Gamma}^{(3)}(\eta_1)
+\tilde{\Gamma}^{(2)}(\check \eta_2^{(2)}(\eta_1))
+\tilde{\Gamma}^{(1)}(\check \eta_2^{\mathrm{r}}(\eta_1)),\\
\check \Gamma^{\mathrm{r}}&=
\tilde{\Gamma}^{\mathrm{r}}(\eta_1+\check \eta_2^{(2)}(\eta_1)+\check \eta_2^{\mathrm{r}}(\eta_1))
+\tilde{\Gamma}^{(3)}(\check \eta_2^{(2)}(\eta_1)+\check \eta_2^{\mathrm{r}}(\eta_1))
+\tilde{\Gamma}^{(2)}(\check \eta_2^\mathrm{r}(\eta_1)),\\
\check \Gamma^{\dag}&=\tilde{\Gamma}(\check\eta_2^\dag(u^\dag), u^\dag).
\end{align*}

\begin{proposition} \label{prop:checkestimates}
$ $
\begin{list}{(\roman{count})}{\usecounter{count}}
\item
The functions $\check \eta_2^{(2)}$, $\check \eta_2^{\mathrm{r}}$ and $\check \eta_2^{\dag}$ 
satisfy the estimates
$$
\| \check \eta_2^{(2)} \|_2 \le c\ve \|\eta_1\|_{0}, \qquad
\|\check \eta_2^\mathrm{r}\|_2 \le c\ve^2 \|\eta_1\|_0, \qquad
\| \check \eta_2^{\dag} \|_2 \le \begin{cases}
c_\lambda \|u^\dag\|_{X^0}, & \lambda>0, \\
c(\|\omega^\dag\|_0+\|\xi_1^\dag\|_1 + \|\xi_3^\dag\|_1), & \lambda=0.
\end{cases}$$
\item
The functions $\check \Gamma^{(1)}$, $\check \Gamma^{(2)}$, $\check \Gamma^{(3)}$, $\check \Gamma^\mathrm{r}$ and $\check \Gamma^{\dag}$ satisfy the estimates
\begin{align*}
\|\nabla\check \Gamma^{(j)}\|_1+ \lambda \|\check \Gamma^{(j)}\|_0 &\le c\ve^{j-1} \|\eta_1\|_{0},\quad j=1,2,3,\\
\|\nabla\check \Gamma^\mathrm{r}\|_1+\lambda\|\check \Gamma^\mathrm{r}\|_{0} &\le c\ve^3\|\eta_1\|_{0},\\
\|\nabla\check \Gamma^\dag\|_1+\lambda\| \check \Gamma^{\dag}\|_0 &\le  \begin{cases}
c_\lambda \|u^\dag\|_{X^0}, & \lambda>0, \\
c(\|\omega^\dag\|_0+\|\xi_1^\dag\|_1 + \|\xi_3^\dag\|_1), & \lambda=0.
\end{cases}
\end{align*}
\end{list}
\end{proposition}

Finally, we write
\[
\check{\NN}(\eta_1, u^\dag)=\check{\NN}^{(2)}(\eta_1)+\check{\NN}^{(3)}(\eta_1)+\check{\NN}^\mathrm{r}(\eta_1)+\check{\NN}^\dag(u^\dag),
\]
 where
\begin{align*}
\check{\NN}^{(2)}(\eta_1)&=\tilde{\NN}^{(2)}(\eta_1),\\
\check{\NN}^{(3)}(\eta_1)&=
\tilde{\NN}^{(3)}(\eta_1)+
\tilde{\NN}^{(2)}(\check \eta_2^{(2)}(\eta_1)),\\
\check{\NN}^\mathrm{r}(\eta_1)&=\tilde{\NN}^\mathrm{r}(\eta_1+\check \eta_2^{(2)}(\eta_1)
+\check{\eta}_2^\mathrm{r}(\eta_1))
+\tilde{\NN}^{(3)}(\check \eta_2^{(2)}(\eta_1)+\check \eta_2^\mathrm{r}(\eta_1))
+\tilde{\NN}^{(2)}(\check \eta_2^\mathrm{r}(\eta_1)),\\
 \check{\NN}^\dag(u^\dag)&=\tilde{\NN}(\check \eta_2^\dag(u^\dag), u^\dag). 
\end{align*}

\begin{proposition}
\label{prop:check N expansion}
The functions $\check{\NN}^{(2)}$, $\check{\NN}^{(3)}$, $\check{\NN}^\mathrm{r}$ and $\check{\NN}^{\dag}$ satisfy the estimates
\begin{align*}
\| \check{\NN}^{(j)}(\eta_1)\|_0 &\le c\ve^{j-1} \|\eta_1\|_{0},\qquad j=2,3,\\
\|\check{\NN}^\mathrm{r}(\eta_1)\|_{0} &\le c\ve^3\|\eta_1\|_{0},\\
\| \check{\NN}^{\dag}(u^\dag)\|_0 &\le \begin{cases}
c_\lambda \|u^\dag\|_{X^0}, & \lambda>0, \\
c(\|\omega^\dag\|_0+\|\xi_1^\dag\|_1 + \|\xi_3^\dag\|_1), & \lambda=0.
\end{cases}
\end{align*}
\end{proposition}

Observing that $\chi(D)\check{N}^{(2)}(\eta)=0$ for all $\eta_1 \in \chi(D)L^2(\R)$, we find that the reduced equation \eqn{eq:Red eqn 3} may be written as
\begin{equation}
g_\ve(\mu,\lambda) \hat{\eta}_1 =\chi \FF[\check \NN^{(3)}(\eta_1) + \check \NN^\mathrm{r}(\eta_1)+\check{\NN}^\dag(u^\dag)].
\label{eq:Red eqn 4}
\end{equation}

\begin{lemma}
There exists a constant $k_\mathrm{max}>0$ such that the linear operator
$$g_\ve(D,\lambda) - \chi(D)\left(\check \NN^{(3)}(\cdot) + \check \NN^\mathrm{r}(\cdot)\right): \chi(D)L^2(\R) \rightarrow \chi(D)L^2(\R)$$ is invertible
for each $\lambda \in [\ve k_\mathrm{max}, \lambda^\star]$. 
\end{lemma}
{\bf Proof.} Write $\lambda = \ve k$ and observe that
$$g_\ve(\mu,\lambda) \geq c\ve^2 (1+k_\mathrm{max}^2)$$
(see Proposition \ref{prop:g}), while
$$\|\check{\NN}^{(3)}(\eta_1) + \check{\NN}^\mathrm{r}(\eta_1)\|_0 \leq c \ve^2 \|\eta_1\|_2.$$
It follows that the given operator is invertible for sufficiently large values of $k$.\qed

\begin{corollary}

For each $\lambda \in [\ve k_\mathrm{max},\lambda^\star]$
and each $u \in X^0$, equation \eqn{eq:Red eqn 4} has a unique
solution $\eta_1 \in \chi(D)L^2(\R)$ which satisfies the estimate
$$\|\eta_1\|_0 \leq c_\lambda \|u^\dag\|_{X^0}.$$
\end{corollary}

In view of the previous corollary we now write $\lambda = \ve k$ and
henceforth suppose that $k \leq k_\mathrm{max}$.

\subsection{Calculation of the leading-order terms in the reduced equation}
\label{Derivation step 3}

In this section we calculate the leading-order terms in the quantity
$\check\NN^{(3)}(\eta_1)$ appearing in the reduced equation \eqn{eq:Red eqn 4},
writing
$$
\eta_1(x)=\frac{1}{2}\ve \zeta(\ve x) \ee^{\ii \mu_0 x}
+\frac{1}{2}\ve \overline{\zeta(\ve x)} \ee^{-\ii \mu_0 x},
$$
where
$\zeta \in \chi_0(\ve D)L^2(\R)$.
For this purpose we first perform the corresponding calculation for
$\check\NN^{(2)}(\eta_1)$ using the formula
\begin{align}
\check{\NN}^{(2)}(\eta_1)&=\tilde{\Gamma}_x^{(2)}(\eta_1)\Big|_{y=1}+k \ve
\int_0^1 y\Phi_{1y}^\star k\ve\check \Gamma^{(1)}\dy \nonumber \\
&\qquad \mbox{}-\int_0^1\left\{ \Phi_{1x}^\star \check \Gamma_x^{(1)}-\Phi_{1y}^\star \check \Gamma_y^{(1)}+\left[
y\Phi_{1y}^\star\check \Gamma_x^{(1)}+y\Phi_{1x}^\star\check \Gamma_y^{(1)}\right]_x \right\} \dy \label{Eqn for checkNN2}
\end{align}
(see equation \eqn{Eqn for tildeNN2}) together with
\begin{equation}
\check \Gamma^{(1)}= \FF^{-1}[-\ii\mu G(y,1)\hat{\eta}_1] \label{Eqn for Phi1}
\end{equation}
and
\begin{equation}
\tilde{\Gamma}^{(2)}(\eta_1) = \FF^{-1}\left[\int_0^1 \left\{ G(y,\tilde{y}) (-\ii\mu \FF[\check{F}_1^{(2)}] - \ii k \ve\FF[\check{F}_2^{(2)}])
+ G_{\tilde{y}}(y,\tilde{y}) \FF[\check{F}_3^{(2)}]\right\}\dty\right],
\label{Eqn for Phi2}
\end{equation}
where
\begin{align}
\check F_1^{(2)}&=-\eta^\star_1\check \Gamma_x^{(1)}-\Phi_{1x}^\star\eta_1 +y\Phi_{1y}^\star\eta_{1x}+y\eta_{1x}^\star \check \Gamma_y^{(1)}, \label{Eqn for F1,2}\\
\check F_2^{(2)}&=-\ii k \ve \eta^\star_1 \check\Phi^{(1)}+\ii k \ve y \Phi_{1y}^\star \eta_1, \label{Eqn for F2,2} \\
\check F_3^{(2)}&=
y\eta_{1x}^\star\check \Gamma_x^{(1)}+y\Phi_{1x}^\star\eta_{1x}
+\eta^\star_1 \check \Gamma_y^{(1)}+\Phi^\star_{1y} \eta_1. \label{Eqn for F3,2}
\end{align}
At various stages of the calculation terms arise of the form $G(y,1)\hat{p}$ and $G(y,\tilde{y})\hat{P}$, where $\hat{p}$ and $\hat{P}$ have
compact support. We decompose such terms into a leading-order part and a higher-order remainder term using the following lemma.

\begin{lemma} $ $
\label{lemma:Green expansion 1}
\begin{list}{(\roman{count})}{\usecounter{count}}
\item
Suppose that $\underline{\mu}>0$. The Green's function $G$ admits the decomposition
\[
G(y, \tilde y)=G(y,\tilde y;\underline \mu)+R(y,\tilde y),
\]
where
\[
G(y, \tilde y;\underline \mu)=\begin{cases}
\displaystyle{\frac{\cosh(\underline \mu y)\cosh(\underline \mu(1-\tilde y))}{\underline \mu \sinh \underline \mu}}, & 0\le y\le \tilde y \le 1,\vspace{5mm}\\
\displaystyle{\frac{\cosh(\underline \mu \tilde y)\cosh(\underline \mu(1-y))}{\underline \mu \sinh \underline \mu}}, &0\le \tilde y\le y\le 1,
\end{cases}.
\]
and $R$ satisfies the estimates
\begin{align*}
|\partial_{y}^{m_1} \partial_{\tilde y}^{m_2} R(y,\tilde y)|&\le c\left(|\mu-\underline{\mu}|+\ve^2k^2\right), \qquad
m_1,m_2=0,1,2,\ldots,
\end{align*}
uniformly over $y, \tilde y$ in $[0,1]$ and $q$ in each fixed interval $[q_\mathrm{min}, q_\mathrm{max}]$
containing $\underline{\mu}$.
\item
The Green's function $G$ admits the expansion
\[
G(y, \tilde y)=G(y,\tilde y;0)+R(y,\tilde y),
\]
where
$$G(y, \tilde y;0) = \frac{1}{q^2}$$
and $R$ satisfies the estimate
$$
|R(y,\tilde y)| \le c
$$
uniformly over $y, \tilde y$ in $[0,1]$ and $q$ in each fixed interval $[0,q_\mathrm{max}]$.
\end{list}
\end{lemma}

The next proposition is proved by applying Lemma \ref{lemma:Green expansion 1}(i)
with $\underline{\mu}=\pm\mu_0$ (and $q_\mathrm{min}=\mu_0-\delta$, $q_\mathrm{max}=((\mu_0+\delta)^2+\lambda^{\star 2})^{1/2}$) to equation \eqn{Eqn for Phi1}.

\begin{proposition} \label{Expansion for Phi1}
The function $\check \Gamma^{(1)}$ is given by the formula
$$
\check \Gamma^{(1)}=- \ii \frac{\cosh(\mu_0 y)}{2\sinh \mu_0} (\ve \zeta(\ve x)\ee^{\ii\mu_0 x}- \ve \overline{\zeta(\ve x)}\ee^{-\ii\mu_0 x})+R_2,
$$
where the symbol $R_j$ denotes a quantity which satisfies the estimate
$$\|R_j\|_j \leq c\ve^{3/2}\|\zeta\|_1.$$
\end{proposition}

Substituting this formula into equations \eqn{Eqn for F1,2}--\eqn{Eqn for F3,2}
and approximating $\zeta^\star$ with
$$
\zeta^\star_\delta:=\chi_0(\ve D)\zeta^\star
$$
by means of the estimate
$$\ve \|\zeta(\ve \cdot)-\zeta_\delta(\ve \cdot)\|_\infty \leq \int_{-\infty}^\infty \left| \hat{\zeta}^\star\left(\frac{\mu}{\ve}\right)
-\hat{\zeta}^\star_\delta\left(\frac{\mu}{\ve}\right)\right| \dmu = \int_{|\mu| \geq \delta} \left| \hat{\zeta}^\star\left(\frac{\mu}{\ve}\right) \right| \dmu
=O(\ve^n), \qquad n>1,
$$
($\zeta$ is a Schwartz-class function)
yields the following result.

\begin{corollary}
The functions $\check{F}_1^{(2)}$, $\check{F}_2^{(2)}$ and $\check{F}_3^{(2)}$ are given by the formulae
\begin{align*}
\check{F}_1^{(2)}
&=
\ve^2 \frac{\mu_0(y\mu_0\sinh(\mu_0)-\cosh(\mu_0 y))}{\sinh(\mu_0)}
\zeta^\star_\delta(\ve x)
\left(\frac12 \zeta(\ve x)e^{2\ii\mu_0 x}+\frac12 \overline{\zeta(\ve x)}e^{-2\ii\mu_0 x}\right)\\
&\qquad\mbox{}-\ve^2\zeta^\star_\delta(\ve x)\frac{\mu_0(\mu_0 y\sinh(\mu_0y)+\cosh(\mu_0 y))}{\sinh(\mu_0)}
\frac1{2} ( \zeta(\ve x)+\overline{\zeta(\ve x)})+\ve R_1, \\
\\
\check{F}_2^{(2)}
&=
\ve^3 k \zeta^\star_\delta(\ve x)
\frac{y\mu_0\sinh(\mu_0 y)-\cosh(\mu_0 y)}{\sinh \mu_0}\left(\frac1{4}  \zeta(\ve x)e^{2\ii\mu_0 x}-\frac1{4} \overline{\zeta(\ve x)}e^{-2\ii\mu_0 x}\right)\\
&\qquad\mbox{}-
\ve^3 k \zeta^\star_\delta(\ve x)
\frac{y\mu_0\sinh(\mu_0 y)+\cosh(\mu_0 y)}{\sinh \mu_0}\frac1{4} ( \zeta(\ve x)-\overline{\zeta(\ve x)})+k\ve^2 R_1,\\
\\
\check{F}_3^{(2)}
&=-\ii
\ve^2  \zeta^\star_\delta(\ve x)\frac{\mu_0(\sinh(\mu_0y)-\mu_0 y \cosh(\mu_0 y))}{\sinh(\mu_0)}
\left(\frac1{2}  \zeta(\ve x)e^{2\ii\mu_0 x}-\frac1{2} \overline{\zeta(\ve x)}e^{-2\ii\mu_0 x}\right)+\ve R_1.
\end{align*}
\end{corollary}

Similarly, the next proposition is obtained from equation \eqn{Eqn for Phi2} using Lemma \ref{lemma:Green expansion 1}(i) with\linebreak
$\underline{\mu}=\pm2\mu_0$ (and $q_\mathrm{min}=2\mu_0-2\delta$, $q_\mathrm{max}=((2\mu_0+2\delta)^2+\lambda^{\star 2})^{1/2}$),
Lemma \ref{lemma:Green expansion 1}(ii) (with $q_\mathrm{max}=(4\delta^2+\lambda^{\star 2})^{1/2}$)
and Proposition \ref{prop:Green mapping}.

\begin{proposition} \label{Expansion for Phi2}
The function $\tilde{\Gamma}^{(2)}(\eta_1)$ is given by the formula
\begin{align*}
\tilde{\Gamma}^{(2)}(\eta_1)
&=-2\ii \ve^2  \zeta^\star_\delta(\ve x)\left(\frac1{2}  \zeta(\ve x)e^{2\ii\mu_0 x}-\frac1{2} \overline{\zeta(\ve x)}e^{-2\ii\mu_0 x}\right) W(y)\\
&\qquad\mbox{}+\ve \frac{\mu_0 \cosh \mu_0}{\sinh \mu_0}\FF^{-1}\left[\frac{\ii \mu}{\mu^2+k^2}\FF[ \re(\zeta^\star\zeta)](\mu)\right](\ve x)+\ve S,
\end{align*}
where 
\[
W(y)=
-\frac{\mu_0\cosh(2\mu_0 y)}{4\sinh^2(\mu_0)}
+\frac{\mu_0y\sinh(\mu_0 y)}{2\sinh(\mu_0)}
\]
and the symbol $S$ denotes a quantity which satisfies the estimate
$$\|\nabla S\|_1 + k\ve\|S\|_0 \leq c\ve^{3/2}\|\zeta\|_1.$$
\end{proposition}

Finally, we combine equation \eqn{Eqn for checkNN2} with the representations of
$\check \Gamma^{(1)}$ and $\tilde{\Gamma}^{(2)}(\eta_1)$ in
Propositions \ref{Expansion for Phi1} and \ref{Expansion for Phi2}.

\begin{lemma}$ $
\label{lemma:N^2}
The quantity $\check{\NN}^{(2)}(\eta_1)$ is given by the formula
\begin{align*}
\check{\NN}^{(2)}(\eta_1)
=\, &-\ve^2 g_0(2\mu_0,0) B_1 \zeta_\delta^\star(\ve x)\left(\frac1{2} \zeta(\ve x)e^{2\ii\mu_0 x}+\frac1{2} \overline{\zeta(\ve x)}e^{-2\ii\mu_0 x}\right)\\
&\qquad\mbox{}-\ve^2\frac{\mu_0 \cosh(\mu_0)}{\sinh(\mu_0)}\FF^{-1}\left[\frac{\mu^2}{\mu^2+k^2}\FF[\re( \zeta_\delta^\star \zeta)]\right](\ve x)
\\
&\qquad\mbox{}
-\ve^2 \frac{\mu_0^2}{2\sinh^2(\mu_0)}
\re(\zeta_\delta^\star(\ve x)\zeta(\ve x))
 + \ve R_0,
\end{align*}
where
\[
B_1=\frac{\mu_0^2(\cosh(2\mu_0)+2)}{2 g_0(2\mu_0,0) \sinh^2 \mu_0}.
\]
\end{lemma}

We now compute the leading-order terms in the quantity
\begin{align}
\check{\NN}^{(3)}(\eta_1)
&=
\check \Gamma^{(3)}_x(\eta_1)+\tilde{\Gamma}^{(2)}_x(\check \eta_2^{(2)})\Big|_{y=1}
-\frac{3\beta_0}{2} \left[(\eta_{1x}^\star)^2 \eta_{1x}\right]_x
+\frac{1}{2}k^2 \ve^2\beta_0 (\eta_{1x}^\star)^2\eta_1
\nonumber \\
&\qquad
\mbox{} -\ii k\ve\int_0^1\left\{-\ii k\ve y^2(\Phi_{1y}^\star)^2\eta_1
 + \ii k\ve y\Phi_{1y}^\star  \check \Gamma^{(2)}
 - \ii k\ve y\Phi_{2y}^\star \check \Gamma^{(1)}
  \right\}\dy \nonumber \\
&\qquad
\mbox{} -\int_0^1\Bigl\{ \Phi_{1x}^\star  \check \Gamma_x^{(2)}-\Phi_{1y}^\star \check \Gamma_y^{(2)}
+2\Phi_{1y}^\star \eta^\star_1 \check \Gamma_y^{(1)}
+(\Phi_{1y}^\star)^2 \eta_1 \} \dy \nonumber \\
&\qquad
\mbox{}-\int_0^1
\left[
y\Phi_{1y}^\star \check \Gamma_x^{(2)}+y\Phi_{1x}^\star \check \Gamma_y^{(2)}
-2y^2\eta_{1x}^\star\Phi_{1y}^\star \check \Gamma_y^{(1)}
-y^2(\Phi_{1y}^\star)^2\eta_{1x}
\right]_x \dy \nonumber \\
&\qquad
\mbox{} -\int_0^1\left\{ \Phi_{2x}^\star \check \Gamma_x^{(1)}(\eta_1)-\Phi_{2y}^\star \check \Gamma_y^{(1)}+\left[
y\Phi_{2y}^\star\check \Gamma_x^{(1)}+y\Phi_{2x}^\star\check \Gamma_y^{(1)}\right]_x \right\} \dy,
\label{Eqn for checkNN3}
\end{align}
which formula follows from equations \eqn{Eqn for tildeNN2}, \eqn{Eqn for tildeNN3} and the
relations
$$\check{\NN}^{(3)}(\eta_1)=
\tilde{\NN}^{(3)}(\eta_1)+\tilde{\NN}^{(2)}(\check \eta_2^{(2)}(\eta_1)), \qquad
\check \Gamma^{(2)} = \tilde{\Gamma}^{(2)}(\eta_1) + \tilde{\Gamma}^{(1)}(\check \eta_2^{(2)}).$$
For this purpose we also use the formulae
$$\tilde{\Gamma}^{(1)}(\check \eta_2^{(2)})= \FF^{-1}[-\ii\mu G(y,1)\FF[\check\eta_2^{(2)}]],
\qquad
\check \eta_2^{(2)}=
g_\ve(D,\lambda)^{-1}(1-\chi(D))\tilde{\NN}^{(2)}(\eta_1),$$
and
\begin{equation}
\label{eq:Final term}
\tilde{\Gamma}^{(3)}(\eta_1)+ \tilde{\Gamma}^{(2)}(\check\eta_2^{(2)})=
\FF^{-1}\left[\int_0^1 (G(y, \tilde y) (-\ii \mu \FF[\check F_1^{(3)}]-\ii\lambda \FF[\check F_2^{(3)}])
+G_{\tilde y}(y,\tilde{y})\FF[\check F_3^{(3)}])\dty\right],
\end{equation}
where
\begin{align*}
\check F_1^{(3)}&=-\eta^\star_1\check \Gamma_x^{(2)}-\eta^\star_2 \check  \Phi_x^{(1)}
-\Phi_{1x}^\star\check \eta_2^{(2)}-\Phi_{2x}^\star\eta_1+y\Phi_{1y}^\star\check \eta_{2x}^{(2)}+y\Phi_{2y}^\star\eta_{1x}
\\
&\qquad
+y\eta_{1x}^\star\check  \Phi_y^{(2)}+y\eta_{2x}^\star \check \Gamma_y^{(1)},
\\
\check F_2^{(3)}&=-\ii k\ve \eta_1^\star\check \Gamma^{(2)}
-\ii k\ve \eta^\star_2 \check \Gamma^{(1)}+\ii k\ve y \Phi_{2y}^\star \eta_1
+\ii k\ve y \Phi_{1y}^\star\check \eta_2^{(2)},
\\
\check F_3^{(3)}&=
y\eta_{1x}^\star \check \Gamma_x^{(2)}
+
y\eta_{2x}^\star\check \Gamma_x^{(1)}
+y\Phi_{1x}^\star\check \eta_{2x}^{(2)}
+y\Phi_{2x}^\star\eta_{1x}
+\eta_1^\star\check \Gamma_y^{(2)}
+\eta^\star_2 \check \Gamma_y^{(1)}
\\
&\qquad
+\Phi^\star_{1y}\check \eta_2^{(2)}
+\Phi^\star_{2y} \eta_1
-(\eta_1^\star)^2\check \Gamma_y^{(1)}
-2\Phi^\star_{1y}\eta_1^\star \eta_1
-y^2(\eta_{1x}^\star)^2\check \Gamma_y^{(1)}-2y^2\eta_{1x}^\star\Phi_{1y}^\star\eta_{1x}.
\end{align*}

\begin{proposition}
\label{prop:check eta_22}
The quantity $\check \eta_2^{(2)}$ satisfies the estimate
\begin{align*}
\check \eta_2^{(2)}
=\, &-\ve^2 B_1 \zeta_\delta^\star(\ve x)\left(\frac1{2} \zeta(\ve x)e^{2\ii\mu_0 x}+\frac1{2} \overline{\zeta(\ve x)}e^{-2\ii\mu_0 x}\right)
-\ve^2 B_2\re(\zeta_\delta^\star(\ve x)\zeta(\ve x))
\\
&
- \ve^2  B_3 \FF^{-1}\left[\frac{\mu^2}{(1-\alpha_0^{-1})\mu^2+k^2 }\FF[2\re(\zeta_\delta^\star\zeta)](\mu)\right](\ve x)
 +\ve R_2,
\end{align*}
where
\[
B_2=\frac{\mu_0^2}{2\alpha_0 \sinh^2 \mu_0}, \qquad
B_3=\alpha_0^{-1}A_4.
\]
\end{proposition}
{\bf Proof.} We apply the operator $g_\ve(D,\lambda)^{-1}(1-\chi(D))$ to the expression for $\check \NN^{(2)}(\eta_1)$ given in Lemma \ref{lemma:N^2}. The first term on the right-hand side of this expression
is handled by the approximation
$$g_\ve(\mu,k\ve)^{-1} = g_0(2\mu_0,0)^{-1} + T_1(\mu,\ve),$$
where
$$|T_1(\mu,\ve)|
\le c(|\mu- 2\mu_0|+\ve^2), \qquad |\mu-2\mu_0| \leq 2 \delta,
$$
while the second and third are treated with the approximation
\[
g_\ve(\mu,k \ve)^{-1}=\frac{q^2}{(\alpha_0-1)\mu^2+\alpha_0 k^2\ve^2}+T_2(\mu,\ve),
\]
where
$$|T_2(q,\ve)| \leq c(\mu^2+\ve^2), \qquad |\mu|<2\delta.$$
Turning to the fourth term, note that
$$g_\ve(D,\lambda)^{-1}(1-\chi(D))R_0=R_2$$
because $g_\ve(\mu,k\ve) > c(1+q^2)$ for $\mu \in \supp(1-\chi)$ (see Proposition \ref{prop:g}).\qed

The function $\tilde{\Gamma}^{(1)}(\check{\eta}_2^{(2)})$ may now be computed using
Lemma \ref{lemma:Green expansion 1}; combining the result with Proposition
\ref{Expansion for Phi2} yields the following proposition.

\begin{proposition} \label{prop:checkPhi2}
The function $\check \Gamma^{(2)}$ is given by the formula
$$
\check \Gamma^{(2)}
=-\ii
\ve^2 \zeta^\star_\delta(\ve x)\left(\frac1{2} \zeta(\ve x)e^{2\ii\mu_0 x}-\frac1{2} \overline{\zeta(\ve x)}e^{-2\ii\mu_0 x}\right)M_1(y)
+\ve\vartheta_\delta(\ve x)+\ve S,
$$
where
\begin{align*}
\vartheta_\delta(x)&=A_4
\FF^{-1}\left[\frac{\ii \mu}{(1-\alpha_0^{-1})\mu^2+k^2}\FF[ 2\re(\zeta^\star_\delta \zeta)]\right](x), \\
M_1(y)&= 
M_{11} \frac{\cosh(2\mu_0 y)}{\sinh(2\mu_0)}
+\frac{\mu_0y\sinh(\mu_0 y)}{\sinh(\mu_0)}, \\
M_{11}&= -\frac{\mu_0^2 (\cosh(2\mu_0)+2)+\mu_0g_0(2\mu_0,0)\sinh(2\mu_0)}{2\sinh^2(\mu_0) g_0(2\mu_0,0)}.
\end{align*}
\end{proposition}
\begin{corollary}
\label{lemma:check F3}
The functions $\check{F}_1^{(3)}$, $\check{F}_2^{(3)}$ and $\check{F}_3^{(3)}$ are given by the
formulae
\begin{align*}
\check F_1^{(3)}&=\ve^3(\zeta^\star_\delta(\ve x))^2\left(\frac1{2} (\zeta(\ve x)U_{11}(y)+\overline{\zeta(\ve x)}U_{12}(y))\ee^{\ii\mu_0 x}+\frac1{2}(\zeta(\ve x)U_{12}(y)+\overline{\zeta(\ve x)}U_{11}(y))\ee^{-\ii\mu_0 x}\right)\\
&\qquad\mbox{}+\ve^3\zeta^\star_\delta(\ve x) \vartheta_{\delta x}(\ve x)U_2(y)
\left(\frac{1}{2}\ee^{\ii\mu_0 x} + \frac{1}{2}\ee^{-\ii \mu_0 x}\right)
+\ve^2 T,\\
\\
\check F_2^{(3)}&=\ve T,\\
\\
\check F_3^{(3)}&=-\ii\ve^3(\zeta^\star_\delta(\ve x))^2\left(\frac1{2}(\zeta(\ve x)U_{31}(y)+\overline{\zeta(\ve x)}U_{32}(y))\ee^{\ii\mu_0 x}-\frac1{2}(\zeta(\ve x)U_{32}(y)+\overline{\zeta(\ve x)}U_{31}(y))\ee^{-\ii\mu_0 x}\right)
\\
&\qquad\mbox{}-\ii\ve^3\zeta^\star_\delta(\ve x) \vartheta_{\delta x}(\ve x)U_4(y)
\left(\frac{1}{2}\ee^{\ii\mu_0 x} - \frac{1}{2}\ee^{-\ii \mu_0 x}\right)
+\ve^2 T,
\end{align*}
where
\begin{align*}
U_{11}(y)&=\frac{\mu_0(\cosh(\mu_0y)+2\mu_0y\sinh(\mu_0 y))B_1}{2\sinh(\mu_0)} +\frac{\mu_0\cosh(\mu_0 y)(C_2+\frac12 B_2)}{\sinh(\mu_0)}\\
&\qquad\mbox{}-C_4-\frac{\mu_0yM_1^\prime(y)}{2}-\mu_0 M_1(y), \\
U_{12}(y)&=\frac{\mu_0(\cosh(\mu_0y)+2\mu_0y\sinh(\mu_0 y))C_1}{2\sinh(\mu_0)} +\frac{\mu_0\cosh(\mu_0 y)B_2}{2\sinh(\mu_0)}\\
&\qquad\mbox{}-\frac{\mu_0C_5(\cosh(2\mu_0 y)+\mu_0 y\sinh(2\mu_0 y))}{\sinh(2\mu_0)}
 -\frac{\mu_0^2 y(y\mu_0\cosh(\mu_0y)+3 \sinh(\mu_0 y))}{4\sinh(\mu_0)}, \\
U_2(y)&=-1-\frac{\mu_0 \cosh(\mu_0 y)}{\alpha_0\sinh(\mu_0)}, \\
U_{31}(y)&=\frac{\mu_0(\sinh(\mu_0y)+2\mu_0y\cosh(\mu_0 y))B_1}{2\sinh(\mu_0)} -\frac{\mu_0\sinh(\mu_0 y)(C_2+\frac12 B_2)}{\sinh(\mu_0)}\\
&\qquad\mbox{}-C_4\mu_0y+\frac{M_1^\prime(y)}{2}+\mu_0^2 y M_1(y)
-\frac{\mu_0 \sinh(\mu_0y)(1+3\mu_0^2 y^2)}{2\sinh(\mu_0)}, \\
U_{32}(y)&=\frac{\mu_0(\sinh(\mu_0y)+2\mu_0y\cosh(\mu_0 y))C_1}{2\sinh(\mu_0)} -\frac{\mu_0\sinh(\mu_0 y) B_2}{2\sinh(\mu_0)}\\
&\qquad\mbox{}+\frac{\mu_0C_5(\mu_0y \cosh(2\mu_0y)+\sinh(2\mu_0y))}{\sinh(2\mu_0)}
+\frac{y\mu_0^2(\cosh(\mu_0y)-\mu_0y \sinh(\mu_0y)}{4\sinh(\mu_0)}, \\
U_4(y)&=-\mu_0 y+\frac{\mu_0 \sinh(\mu_0 y)}{\alpha_0\sinh(\mu_0)}
\end{align*}
and the symbol $T$ denotes a quantity which satisfies the estimate
$$\|\chi(D)T\|_0 \leq c\ve^{3/2}\|\zeta\|_1.$$
\end{corollary}

The next proposition is obtained from equation \eqn{eq:Final term} using Lemma \ref{lemma:Green expansion 1} and Proposition \ref{prop:Green mapping}.

\begin{proposition}
\label{prop:Final term}
The quantity $\tilde{\Gamma}_x^{(3)}(\eta_1) + \tilde{\Gamma}_x^{(2)}(\check \eta_2^{(2)})\Big|_{y=1}$ is given by the formula
\begin{align*}
\tilde{\Gamma}_x^{(3)}(\eta_1) &+ \tilde{\Gamma}_x^{(2)}(\check \eta_2^{(2)})\Big|_{y=1} \\
&=
\mu_0\ve^2 (\zeta_\delta(\ve x))^2
\left(\frac{1}{2} (\zeta(\ve x)N_1+\overline{\zeta(\ve x)}N_2)e^{\ii\mu_0 x}
+\frac{1}{2} (\zeta(\ve x)N_2+\overline{\zeta(\ve x)}N_1)e^{-\ii\mu_0 x}\right)\\
&\qquad\mbox{}+\mu_0\ve^3\zeta_\delta^\star(\ve x) \vartheta_{\delta x}(\ve x)N_3\left(\frac{1}{2}\ee^{\ii\mu_0x}+\frac{1}{2}\ee^{-\ii \mu_0 x}\right) + \ve^2 T,
\end{align*}
where
\begin{align*}
N_1&=
\frac{\mu_0\cosh(2\mu_0)B_1}{2\sinh^2(\mu_0)}
+\frac{\mu_0(C_2+\frac12 B_2)}{\sinh^2(\mu_0)}-\frac{\cosh(\mu_0)C_4}{\sinh(\mu_0)}
-\frac{M_{11}\mu_0}{2\sinh^2(\mu_0)}
-\frac{\mu_0^2\cosh(\mu_0)}{2\sinh(\mu_0)},\\
N_2&=
\frac{\mu_0\cosh(2\mu_0)C_1}{2\sinh^2(\mu_0)}
+\frac{\mu_0B_2}{2\sinh^2(\mu_0)}
-\frac{\mu_0 C_5}{2\sinh^2(\mu_0)}
-\frac{\mu_0^2\cosh(\mu_0)}{4\sinh(\mu_0)},
\\
N_3&=\frac{-2\mu_0\alpha_0^{-1}-\sinh(2\mu_0)}{2\sinh^2(\mu_0)}.
\end{align*}
\end{proposition}

Finally, we combine equation \eqn{Eqn for checkNN3} with
Propositions \ref{prop:checkPhi2} and \ref{prop:Final term}.

\begin{lemma}
\label{lemma:N^3}
The quantity $\check{\NN}^{(3)}(\eta_1)$ is given by the formula
\begin{align*}
\check{\NN}^{(3)}(\eta_1)&=
\ve^3(\zeta^\star_\delta(\ve x))^2 \left( (A_3+A_5) \zeta(\ve x)+A_3\overline{\zeta(\ve x)}\right)\frac{\ee^{\ii \mu_0 x}}{2}\\
&\qquad\mbox{}+\ve^3(\zeta^\star_\delta(\ve x))^2 \left( A_3\zeta(\ve x)+(A_3+A_5)\overline{\zeta(\ve x)}\right)\frac{\ee^{-\ii \mu_0 x}}{2}\\
&\qquad
\mbox{}-4\ve^3 A_4 \zeta^\star_\delta(\ve x) \vartheta_{\delta x}(\ve x)
\left(\frac{1}{2}\ee^{\ii\mu_0 x} + \frac{1}{2}\ee^{-\ii \mu_0 x}\right)+\ve^2T.
\end{align*}
\end{lemma}

\subsection{Derivation of the reduced equation in its final form}
\label{Derivation step 4}

In this section we derive the reduced equation in its final form, thus completing the proof
of Theorem \ref{thm:reduction}(ii) and (iii).
For $k>0$ we define $\check{u}_{\ve,k}:\chi_0(\ve D) L^2(\R) \times X^0 \rightarrow Y^1$ by the formula
\[
\check u_{\ve, k}(\zeta, u^\dag)=(\eta_1+\check \eta_2(\eta_1, u^\dag), \check \omega(\eta_1, u^\dag), 
\check \Gamma(\eta_1, u^\dag), \check \xi(\eta_1, u^\dag)),
\]
while $\check{u}_{\ve,0}:\chi_0(\ve D) L_\mathrm{c}^2(\R) \times Z_\mathrm{r}^0 \rightarrow Y_\star^1$ is defined by
\[
\check u_{\ve, 0}(\zeta, \eta^\dag,\omega^\dag, \Gamma^\dag,\xi_1^\dag,\xi_3^\dag)=(\eta_1+\check \eta_2(\eta_1, u^\dag), \check \omega(\eta_1, u^\dag), 
\check \Gamma(\eta_1, u^\dag), \check \xi(\eta_1, u^\dag))
\]
with $\xi^\dag=-(\xi_1^\dag)_x-(\xi_3^\dag)_y$;
here
\[
\eta_1=\ve\zeta(\ve x) \frac{\ee^{\ii\mu_0 x}}{2}+\ve\overline{\zeta(\ve x)} \frac{\ee^{-\ii\mu_0 x}}{2}
\]
and $\check \omega$, $\check \xi$ are given by equations \eqn{eq:omega}, \eqn{eq:xi} with $\Phi=\check \Gamma(\eta_1, u^\dag)$, $\eta=\eta_1+\check \eta_2(\eta_1, u^\dag)$.
The theory in Sections \ref{Derivation step 1}--\ref{Derivation step 4} above shows that
$$\eta_1^+(x) = \ve\zeta(\ve x) \frac{\ee^{\ii\mu_0 x}}{2}, \qquad
\zeta \in 
\begin{cases}
\chi_0(\ve D)L^2(\R), & k>0, \\
\chi_0(\ve D)L_\mathrm{c}^2(\R), & k=0, \\
\end{cases}
$$
solves the reduced equation
$$g_\ve(\mu,k\ve)\hat{\eta}_1^+ = \chi_+\FF[
\check \NN^{(3)}(\eta_1) + \check \NN^\mathrm{r}(\eta_1)+\check{\NN}^\dag(u^\dag)],
$$
where $\eta_1=\eta_1^+ + \overline{\eta_1^+}$,
if and only if $u=\check u_{\ve, k}$ solves the resolvent equations
\eqn{eq:spectral 1}--\eqn{eq:spectral BC2}.
It follows that
$$\eta_1^+(x) = \ve\zeta(\ve x) \frac{\ee^{\ii\mu_0 x}}{2}, \qquad
\zeta \in \begin{cases}
L^2(\R), & k>0, \\
L_\mathrm{c}^2(\R), & k=0, \\
\end{cases}
$$
solves the equation
\begin{equation}
g_\ve(\mu,k\ve)\hat{\eta}_1^+ = \chi_+\FF[
\check \NN^{(3)}(\chi(D)\eta_1) + \check \NN^\mathrm{r}(\chi(D)\eta_1)+\check{\NN}^\dag(u^\dag)],
\label{Final reduced}
\end{equation}
where $\eta_1=\eta_1^+ + \overline{\eta_1^+}$,
if and only if $\zeta \in \chi_0(\ve D)L^2(\R)$ (for $k>0$) or $\zeta \in \chi_0(\ve D)L_\mathrm{c}^2(\R)$
(for $k=0$) and
$u=\check u_{\ve, k}$ solves the resolvent equations
\eqn{eq:spectral 1}--\eqn{eq:spectral BC2}.

Since all solutions of \eqn{Final reduced} have support in the interval $[\mu_0-\delta,\mu_0+\delta]$,
its solution set in $L^2(\R)$ coincides with its solution set
in $H^s(\R)$ for any $s \geq 0$, and we henceforth work in $H^2(\R)$ by
writing it as a second-order pseudodifferential
equation. Let
$$
 \tilde g_\ve(\mu, \lambda)=\ve^2 +A_1 (\mu-\mu_0)^2 +A_2 \lambda^2
$$
be the second-order Taylor polynomial of $g_\ve$ at the point $(\mu_0,0)$
and write \eqn{Final reduced} as
\[
\tilde g_\ve(\mu, \ve k)\hat{\eta}_1^+=\frac{\tilde g_\ve(\mu, \ve k)}{g_\ve(\mu, \ve k)} 
\chi_+\FF[
\check \NN^{(3)}(\chi(D)\eta_1) + \check \NN^\mathrm{r}(\chi(D)\eta_1)+\check{\NN}^\dag(u^\dag)]
\]
or equivalently as
\[
\tilde g_\ve(\mu_0+\ve \tilde \mu, \ve k)\hat \zeta(\tilde  \mu)=2\frac{\tilde g_\ve(\mu_0+\ve \tilde \mu, \ve k)}{g_\ve(\mu_0+\ve \tilde \mu, \ve k)} 
\chi_0(\ve \tilde \mu)
\FF[\check \NN^{(3)}(\chi(D)\eta_1) + \check \NN^\mathrm{r}(\chi(D)\eta_1)+\check{\NN}^\dag(u^\dag)]
 (\mu_0+\ve \tilde \mu).
\]
Taking the inverse Fourier transform with respect to $\tilde{\mu}$, we find that
\begin{align}
A_2^{-1}&\zeta - A_2^{-1}A_1\zeta_{xx}+k^2\zeta \nonumber \\
&= \FF^{-1}\bigg[
2\frac{\tilde g_\ve(\mu_0+\ve \tilde \mu, \ve k)}{g_\ve(\mu_0+\ve \tilde \mu, \ve k)} 
\chi_0(\ve \tilde \mu)\ve^{-2}A_2^{-1} \nonumber \\
&
\qquad\qquad\quad\times\FF[\check\NN^{(3)}(\chi(D)\eta_1)+\check\NN^\mathrm{r}(\chi(D)\eta_1)](\mu_0+\ve \tilde \mu)+\hat{\zeta}_{\ve, k}^\dag(u^\dag)(\tilde{\mu})
\bigg], \label{Final reduced 2}
\end{align}
where
\[
\FF[\zeta_{\ve, k}^\dag(u^\dag)](\tilde{\mu})=2A_2^{-1}\frac{\tilde g_\ve(\mu_0+\ve \tilde \mu, \ve k)}{g_\ve(\mu_0+ \ve \tilde  \mu,\ve k)}
\chi_0(\ve \tilde \mu)\ve^{-2} \FF[\check{\NN}^\dag(u^\dag)](\mu_0+\ve \tilde \mu).
\]

Observe that
\begin{align*}
\check{\NN}^{(3)}&(\chi(D)\eta_1) \\
&=
\ve^3(\zeta^\star_\delta(\ve x))^2 \left( (A_3+A_5)(\chi_0(\ve D) \zeta)(\ve x)+A_3\overline{(\chi_0(\ve D) \zeta)(\ve x)}\right)\frac{\ee^{\ii \mu_0 x}}{2}\\
&\qquad\mbox{}+\ve^3(\zeta^\star_\delta(\ve x))^2 \left( A_3(\chi_0(\ve D) \zeta)(\ve x)+(A_3+A_5)\overline{(\chi_0(\ve D) \zeta)(\ve x)}\right)\frac{\ee^{-\ii \mu_0 x}}{2}\\
&\qquad
\mbox{}+\ve^3 A_4 \zeta^\star_\delta(\ve x)(\chi_0(\ve D)\vartheta_{\delta})_x(\ve x)
\left(\frac{1}{2}\ee^{\ii\mu_0 x} + \frac{1}{2}\ee^{-\ii \mu_0 x}\right)+\ve^2T,
\end{align*}
and using the estimate
\begin{equation}
\|(\chi_0(\ve \tilde{\mu})-1)\hat{f}(\tilde{\mu})\|_0^2
=\int_{|\tilde{\mu}| \geq \delta/\ve} |\hat{f}(\tilde{\mu})|^2 \dtmu
\leq \delta^{-2}\ve^2 \int_{|\tilde{\mu}| \geq \delta/\ve} |\tilde{\mu}|^2|\hat{f}(\tilde{\mu})|^2 \dtmu
\leq \delta^{-2}\ve^2 \|f\|_1, \label{eq:Integral trick}
\end{equation}
we find that
$$
\|((\chi_0(\ve D)-1)\zeta)(\ve x)\ee^{\ii \mu_0 x} \|_0^2
= \ve^{-1} \|(\chi_0(\ve \tilde{\mu})-1)\hat{\zeta}(\tilde{\mu})\|_0^2
\leq \delta^{-2}\ve \|\zeta\|_1
$$
and
\begin{align*}
\|((\chi_0(\ve D)-1)\vartheta_\delta)_x(\ve x)\ee^{\ii \mu_0 x} \|_0^2
&= \ve^{-1}A_4
\left\| \frac{(\chi_0(\ve \tilde{\mu})-1)\tilde{\mu}^2}{(1-\alpha_0^{-1})\tilde{\mu}^2+k^2}
\FF[2\re(\zeta_\delta^\star \zeta)](\tilde{\mu})\right\|_0 \\
&\leq c\ve^{-1}
\| (\chi_0(\ve \tilde{\mu})-1)
\FF[2\re(\zeta_\delta^\star \zeta)](\tilde{\mu})\|_0 \\
&\leq c\ve \|\zeta\|_1,
\end{align*}
so that
\begin{equation}
\check{\NN}^{(3)}(\chi(D)\eta_1)
=\frac{\ve^3}{2}E(\ve x)\ee^{\ii\mu_0 x}+\frac{\ve^3}{2}F(\ve x)\ee^{-\ii\mu_0 x}+\ve^2T, \label{eq:checkNN3 estimate to E}
\end{equation}
where
\begin{align*}
E&=(A_3+A_5)(\zeta_\delta^\star)^2\zeta
+A_3(\zeta_\delta^\star)^2\overline{\zeta}
-4A_4\zeta_\delta^\star \vartheta_{\delta x}, \\
F&=A_3(\zeta_\delta^\star)^2\chi_0(\ve D)\zeta
+(A_3+A_5)(\zeta_\delta^\star)^2\chi_0(\ve D)\overline{\zeta}
-4A_4\zeta_\delta^\star(\chi_0(\ve D)\vartheta_\delta)_x.
\end{align*}
Furthermore
$$\chi_0(\ve \tilde{\mu})\FF[\ve F(\ve x)\ee^{-\ii \mu_0 x}](\mu_0+\ve\tilde{\mu})
=\chi_0(\ve \tilde{\mu}) \hat{F}(\tilde{\mu}+2\ve^{-1}\mu_0) = 0$$
because $\supp \hat{F} \subseteq [-3\ve^{-1}\delta,3\ve^{-1}\delta]$.

Using the estimate
\begin{equation}
\left|\frac{\tilde{g}_\ve(\mu,\lambda)}{g_\ve(\mu,\lambda)}-1\right|
\leq c|(\mu-\mu_0,\lambda)|^2, \qquad
|\mu-\mu_0| \leq \delta,\ \lambda \leq \lambda^\star,
\label{eq:Quotient minus one}
\end{equation}
we find in particular that
$$\left|\frac{\tilde g_\ve(\mu_0+\ve \tilde \mu, \ve k)}{g_\ve(\mu_0+\ve \tilde \mu, \ve k)} 
\chi_0(\ve \tilde \mu)\right| \leq c.$$
It follows that
\begin{eqnarray*}
\left\|\FF^{-1}\left[
\frac{\tilde g_\ve(\mu_0+\ve \tilde \mu, \ve k)}{g_\ve(\mu_0+\ve \tilde \mu, \ve k)} 
\chi_0(\ve \tilde \mu)\ve^{-2}
 \FF[\check\NN^\mathrm{r}(\chi(D)\eta_1)](\mu_0+\ve \tilde \mu)\right]\right\|_0
 & \leq & c \ve^{-5/2} \|\check\NN^\mathrm{r}(\chi(D)\eta_1)\|_0 \\
 & \leq & c\ve^{-1/2} \|\chi(D) \eta_1\|_0 \\
 & \leq & c\ve \|\zeta\|_0
 \end{eqnarray*}
 and
$$
\left\|\FF^{-1}\left[
\frac{\tilde g_\ve(\mu_0+\ve \tilde \mu, \ve k)}{g_\ve(\mu_0+\ve \tilde \mu, \ve k)} 
\chi_0(\ve \tilde \mu)\ve^{-2}
 \FF[T](\mu_0+\ve \tilde \mu)\right]\right\|_0
\leq c\ve^{-1/2} \|\chi_+(D) T\|_0 \\
\leq c\ve \|\zeta\|_1.
$$

Finally, we write
\begin{eqnarray*}
\lefteqn{\FF^{-1}\left[2
\frac{\tilde g_\ve(\mu_0+\ve \tilde \mu, \ve k)}{g_\ve(\mu_0+\ve \tilde \mu, \ve k)} 
\chi_0(\ve \tilde \mu)\ve^{-2}
 \FF\left[\frac{\ve^3}{2}E(\ve x)\ee^{\ii \mu_0 x}\right]
(\mu_0+\ve \tilde \mu)\right]} \\
& = & \FF^{-1}\left[\left(
\frac{\tilde g_\ve(\mu_0+\ve \tilde \mu, \ve k)}{g_\ve(\mu_0+\ve \tilde \mu, \ve k)} -1 \right)
\chi_0(\ve \tilde \mu)\hat{E}(\tilde{\mu})\right]
+\FF^{-1}[(\chi_0(\ve \tilde{\mu})-1)\hat{E}(\tilde{\mu})]+E(x)
 \end{eqnarray*}
and note that
$$\left\|
\FF^{-1}\left[\left(
\frac{\tilde g_\ve(\mu_0+\ve \tilde \mu, \ve k)}{g_\ve(\mu_0+\ve \tilde \mu, \ve k)} -1 \right)
\chi_0(\ve \tilde \mu)\hat{E}(\tilde{\mu})\right]
\right\|_0
\leq c (\ve \delta \|E^\prime\|_0 + k^2 \ve^2 \|E\|_0)
\leq c \ve \|\zeta\|_1,$$
(using estimate \eqn{eq:Quotient minus one}), and
$$
\|\FF^{-1}[(\chi_0(\ve \tilde{\mu})-1)\hat{E}(\tilde{\mu})]\|_0 \leq c \ve \|E\|_1
\leq c \ve\|\zeta\|_1$$
(using estimate \eqn{eq:Integral trick}).

According to the above calculations equation \eqn{Final reduced 2} can be written as
\begin{eqnarray*}
\lefteqn{A_2^{-1}\zeta - A_1A_2^{-1}\zeta_{xx}+k^2\zeta} \\
& & = A_2^{-1}
(A_3+A_5)(\zeta_\delta^\star)^2\zeta
+A_2^{-1}A_3(\zeta_\delta^\star)^2\overline{\zeta}
-4A_2^{-1}A_4\zeta_\delta^\star \vartheta_{\delta x}
+\zeta^\dag_{\ve, k}(u^\dag) + \tilde{\RR}_{\ve, k}(\zeta) \\
& & = A_2^{-1}
(A_3+A_5)(\zeta^\star)^2\zeta
+A_2^{-1}A_3(\zeta^\star)^2\overline{\zeta}
-4A_2^{-1}A_4\zeta^\star \vartheta_x
+\zeta^\dag_{\ve, k}(u^\dag) + \RR_{\ve, k}(\zeta),
\end{eqnarray*}
where
$$\vartheta = A_4 \FF^{-1}\left[\frac{\ii \mu}{(1-\alpha_0^{-1})\mu^2+k^2} \FF[2\re(\zeta^\star \zeta)]\right]$$
and
$$\|\tilde{\RR}_{\ve, k}(\zeta)\|_0,\ \|\RR_{\ve, k}(\zeta)\|_0 \leq c \ve \|\zeta\|_1.$$
In view of this calculation we now turn our attention to the solvability of the reduced equation
\begin{eqnarray}
\lefteqn{A_2^{-1}\zeta-A_2^{-1}A_1\zeta_{xx}} \nonumber \\
& & \mbox{}-A_2^{-1}
(A_3+A_5)(\zeta^\star)^2\zeta
-A_2^{-1}A_3(\zeta^\star)^2\overline{\zeta}
+4A_2^{-1}A_4\zeta^\star \vartheta_x
+k^2 \zeta - \RR_{\ve, k}(\zeta) = \zeta^\dag\qquad
\label{Final reduced 3}
\end{eqnarray}
for arbitrary $\zeta^\dag \in L^2(\R)$ (for $k>0$) or $\zeta^\dag \in L_\mathrm{c}^2(\R)$ (for $k=0$).
For $k>0$ we also note that
each step in the reduction procedure preserves the invariance of the resolvent
equations under the reflection $R$; in the present coordinates the action of this
symmetry is given by $\zeta(x) \mapsto \overline{\zeta(-x)}$, so that \eqn{Final reduced 3}
is invariant under this transformation.

\section{Spectral theory for the reduced equation}
\label{sec:analysis}

\subsection{Invertibility of the reduced operator}
\label{Analysis step 1}

In this section we reformulate the reduced equation \eqn{Final reduced 3} once more,
reducing the question of its solvability to determining the location of an eigenvalue
of an unbounded linear operator.
We first present an auxiliary spectral result, one consequence of which is the solvability
of \eqn{Final reduced 3} for $k=0$.

\begin{proposition} \label{prop:Scaling of standard spectral results}
The formulae
\begin{align*}
\CC_{0,1}(\zeta_1)&=A_2^{-1}\big(\zeta_1-A_1 \zeta_{1xx}-6\sech^2(A_1^{-1/2} x) \zeta_1\big), \\
\CC_{0,2}(\zeta_2)&=A_2^{-1}\big(\zeta_2-A_1 \zeta_{2xx}-2\sech^2(A_1^{-1/2} x) \zeta_2\big)
\end{align*}
define
\begin{list}{(\roman{count})}{\usecounter{count}}
\item
self-adjoint operators $\CC_{0,1}$, $\CC_{0,2}: H^2(\R) \subseteq L^2(\R) \rightarrow L^2(\R)$
whose spectrum consists of essential spectrum $[A_2^{-1},\infty)$ and respectively
two simple eigenvalues at $-3A_2^{-1}$, $0$ and a simple eigenvalue at $0$;
\item
self-adjoint operators $\CC_{0,1}: H^2_\mathrm{e}(\R) \subseteq L^2_\mathrm{e}(\R)
\rightarrow L^2_\mathrm{e}(\R)$,
$\CC_{0,2}: H^2_\mathrm{o}(\R) \subseteq L^2_\mathrm{o}(\R)
\rightarrow$~$L^2_\mathrm{o}(\R)$ whose spectrum consists of essential spectrum
$[A_2^{-1},\infty)$ and, in the case of $\CC_{0,1}$, a simple eigenvalue at $-3A_2^{-1}$.
\end{list}
\end{proposition}
{\bf Proof.} The spectra of the operators
$1-\partial_x^2-6 \sech^2(x): H^2(\R)\subseteq L^2(\R)\to L^2(\R)$ and
$1-\partial_x^2-2 \sech^2(x) : H^2(\R)\subseteq L^2(\R)\to L^2(\R)$
consists of essential spectrum $[1,\infty)$ and respectively
two simple eigenvalues at $-3$ and $0$ (with corresponding eigenvectors
$\sech^2(x)$ and $\sech^\prime(x)$) and a simple eigenvalue at $0$
(with corresponding eigenvector $\sech(x)$) (see Drazin \cite[Chapter 4.11]{Drazin}). The assertions follow by a scaling argument
and restricting to respectively even and odd functions.\qed
\begin{lemma}
Suppose that $k=0$.
The reduced equation \eqn{Final reduced 3} has a unique solution $\zeta \in H_\mathrm{c}^2(\R)$
for each $\zeta^\dag \in L_\mathrm{c}^2(\R)$.
\end{lemma}
{\bf Proof.} Write equation \eqn{Final reduced 3} as
$$\CC_\ve(\zeta) = \zeta^\dag,$$
where
$$\CC_\ve = \CC_0 + \RR_\ve, \qquad \CC_0 (\zeta) = (\CC_{0,1} (\zeta_1), \CC_{0,2} (\zeta_2))$$
and $\zeta_1 = \re \zeta$, $\zeta_2 = \im \zeta$, and observe that
$\CC_0: H_\mathrm{c}^2(\R) \rightarrow L_\mathrm{c}^2(\R)$ is invertible
(see Proposition \ref{prop:Scaling of standard spectral results}(ii))
and $\|\RR_\ve(\zeta)\|_0 \leq c\ve \|\zeta\|_2$.\qed

Let $W$ be the real Hilbert space  $L^2(\R)\times L^2(\R)\times L^2(\R)$ 
equipped with the inner product
\[
\llangle (\zeta_1, \zeta_2,\psi), (\tilde \zeta_1, \tilde \zeta_2, \tilde \psi)\rrangle = \langle \zeta_1, \tilde \zeta_1\rangle_0
+\langle \zeta_2, \tilde \zeta_2\rangle_0 + 2A_2^{-1}\langle \psi, \tilde \psi\rangle_0
\]
and define unbounded operators $\BB_{\ve,k}$, $\tilde{\BB}_{\ve,k}: \DD_\BB \subseteq W \rightarrow W$
by the formulae
\[
\tilde\BB_{\ve, k} \begin{pmatrix} \zeta_1 \\ \zeta_2 \\ \psi \end{pmatrix}=
\begin{pmatrix}
A_2^{-1}\zeta_1-A_2^{-1}A_1\zeta_{1xx}-A_2^{-1}(2A_3+A_5)(\zeta^\star)^2\zeta_1
+4A_2^{-1}A_4 \zeta^\star \psi_x-\re \RR_{\ve, k}(\zeta)\\
A_2^{-1}\zeta_2-A_2^{-1}A_1\zeta_{2xx}-A_2^{-1}A_5(\zeta^\star)^2\zeta_2-\im \RR_{\ve, k}(\zeta)\\
-(1-\alpha_0^{-1})\psi_{xx}-2A_4 (\zeta^\star \zeta_1)_x
\end{pmatrix}
\]
and
\[
\BB_{\ve ,k} =\tilde\BB_{\ve, k}+k^2I,
\]
where $\DD_\BB:=H^2(\R)\times H^2(\R)\times H^2(\R)$.
In this framework equation \eqn{Final reduced 3} may be written as
$$\BB_{\ve,k}(\zeta_1,\zeta_2,\psi) = (\re \zeta^\dag,\im \zeta^\dag,0),$$
where $\zeta_1 = \re \zeta$, $\zeta_2=\im \zeta$; recall that this equation
is invariant under the transformation $\tilde{R}:(\zeta_1(x),\zeta_2(x),\psi(x))=(\zeta_1(-x),-\zeta_2(-x),-\psi(-x))$
(the action of the reflection $R$).

The spectra of the self-adjoint
operators $\BB_{0,k}$ and $\tilde{\BB}_{0,k}$ (the second of which does not depend upon $k$)
can be determined precisely.

\begin{lemma} \label{lemma:spec B0}
The spectrum of the operator $\tilde{\BB}_{0,k}$ consists of essential spectrum $[0,\infty)$ and a simple
negative eigenvalue $-k_0^2$ whose eigenspace lies in $\mathrm{Fix}\, \tilde{R}=
L_\mathrm{e}^2(\R)\times L_\mathrm{o}^2(\R)\times L_\mathrm{o}^2(\R)$.
\end{lemma}
{\bf Proof.} First note that $\tilde\BB_{0,k}$ is a compact perturbation of the constant-coefficient operator $\DD_\BB \subseteq W \rightarrow W$
defined by
\[
(\zeta_1,\zeta_2,\psi)\mapsto (A_2^{-1} \zeta_1-A_2^{-1}A_2 \zeta_{1xx},A_2^{-1}\zeta_2-A_2^{-1}A_1\zeta_{2xx},-(1-\alpha_0^{-1}) \psi_{xx}),
\]
whose essential spectrum is clearly $[0,\infty)$; it follows that $\sigma_\mathrm{ess}(\tilde{\BB}_{0,k}) = [0,\infty)$
(see Kato, \cite[Chapter IV, Theorem 5.26]{Kato}). Because $\tilde{\BB}_{0,k}$ is self-adjoint the remainder of its spectrum
consists of negative real eigenvalues with finite multiplicity. Write
$W^1=L^2(\R)\times L^2(\R)$,
$\DD_B^1 = H^2(\R)\times H^2(\R)$
and $W^2=L^2(\R)$, $\DD_B^2=H^2(\R)$
and observe that (with a slight abuse of notation)
$$\tilde\BB_{0,k}(\zeta_1,\zeta_2,\psi)=(\tilde\BB_{0,k,1}(\zeta_1,\psi), \tilde\BB_{0,k,2}\zeta_2),$$
where $\tilde{\BB}_{0,  k,1}: \DD_\BB^1 \subseteq W^1 \rightarrow W^1$ and
$\tilde{\BB}_{0,  k,2}: \DD_\BB^2 \subseteq W^2 \rightarrow W^2$
are given by
\[
\tilde{\BB}_{0,  k,1} \begin{pmatrix} \zeta_1 \\ \psi \end{pmatrix} =
\begin{pmatrix}
A_2^{-1}\zeta_1-A_2^{-1}A_1\zeta_{1xx}-A_2^{-1}(2A_3+A_5)(\zeta^\star)^2\zeta_1
+4A_2^{-1}A_4 \zeta^\star \psi_x\\
-(1-\alpha_0^{-1})\psi_{xx}-2A_4(\zeta^\star \zeta_1)_x
\end{pmatrix}
\]
and
$$
\tilde{\BB}_{0, k,2}\zeta_2=\CC_{0,2}\zeta_2;
$$
the eigenvalues of $\tilde{\BB}_{0,k}$ are therefore precisely the eigenvalues of $\tilde{\BB}_{0,k,1}$ and $\tilde{\BB}_{0,k,2}$.

According to Proposition \ref{prop:Scaling of standard spectral results} the operator
$\tilde{\BB}_{0,k,2}$ has no negative eigenvalues.
Turning to the spectrum of $\tilde\BB_{0, k, 1}$, one finds by
an explicit calculation that
$$
\langle \tilde\BB_{0,  k,1}(\zeta_1,\psi), (\zeta_1,\psi)\rangle_{W^1}
=
\langle \CC_{0,1}\zeta_1,\zeta_1 \rangle_0
+2A_2^{-1}(1-\alpha_0^{-1})
\int_{\R} \left(\psi_x+\frac{2A_4}{1-\alpha_0^{-1}}  \zeta^\star \zeta_1\right)^2\dx,
$$
which quantity is positive for $(\zeta_1,\psi) \in W^1_+,$ where
\[W^1_+=\{(\zeta_1,\psi)\in W^1 : \langle (\zeta_1,\psi), ((\sech^2(A_1^{-1/2}x),0)\rangle_{W^1}=0\}.\] 
It follows that any subspace of $W^1$ upon which $\tilde{\BB}_{0,k,1}$
is strictly negative definite is one-dimensional.
The calculation
\[
\langle \tilde{\BB}_{0, k,1} ((\sech(A_1^{-1/2} x),0) ,(\sech(A_1^{-1/2} x),0)\rangle_{W^1}=-\frac{16\sqrt{A_1}A_3}{3A_2A_5}<0
\]
shows that $\inf \sigma(\tilde{\BB}_{0, k,1})<0$, so that the spectral subspace of $W^1$
corresponding to the part of the spectrum of $\tilde{\BB}_{0, k,1}$ in $(-\infty,-\ve)$
is nontrivial and hence one-dimensional for every sufficiently small value of $\ve>0$.
We conclude that $\tilde{\BB}_{0,k,1}$ has precisely one simple negative eigenvalue
$-k_0^2$.

Finally, the same argument shows that $\tilde{\BB}_{0,k}|_{\mathrm{Fix}\, \tilde{R}}$ also has precisely
one simple negative eigenvalue. It follows that this eigenvalue is $-k_0^2$, whose
eigenspace therefore lies in $\mathrm{Fix}\, \tilde{R}$.\qed

Noting that
\begin{equation}
\|\tilde\BB_{\ve,k}-\tilde\BB_{0,k}\|_{\LL(\DD_\BB,W)} \leq c\ve,
\label{BB estimate}
\end{equation}
we now use a perturbation argument to obtain a qualitative description of
a portion of  the spectrum of $\tilde{\BB}_{\ve,k}$.

\begin{lemma}$ $
\label{lemma:spectrum of tilde B}
Let $m$ and $M$ be positive real numbers with $m<k_0^2<M$ and $\gamma$ 
be an ellipse in the complex plane with major axis $[-M, -m]$. For each 
$k\in [0, k_\mathrm{max}]$ the portion of the spectrum of $\tilde\BB_{\ve, k}$ within $\gamma$ 
consists of precisely one simple real eigenvalue $\kappa_{\ve, k}$ with
$\kappa_{0, k}=-k_0^2$; its eigenspace lies in $\mathrm{Fix}\, \tilde{R}$. In particular, $\tilde\BB_{\ve, k}$ is closed
and $\tilde\BB_{\ve,k}-\kappa I:\DD_\BB \rightarrow W$ is Fredholm with index $0$ for $k\in [0, k_\mathrm{max}]$ and $\kappa \in [-M, -m]$.
\end{lemma}
{\bf Proof.}
The contour $\gamma$ defines a separation of $\sigma(\tilde{\BB}_{0,k})$: the portion of
$\sigma(\tilde{\BB}_{0,k})$ in the interior of $\gamma$ consists of precisely one simple eigenvalue
$-k_0^2$ and the remainder of $\sigma(\tilde{\BB}_{0,k})$ lies in the exterior of $\gamma$.
In view of inequality \eqn{BB estimate} and the fact that
$\tilde\BB_{0,k}$ does not depend upon $k$, we may apply a standard argument
in spectral perturbation theory (see Kato \cite[Theorem 3.16]{Kato}), which asserts that
$\gamma$ defines the same separation of $\tilde\BB_{\ve,k}$. Since $\tilde{\BB}_{\ve, k}$
is a real operator its eigenvalues arise in complex-conjugate pairs; its simple eigenvalue
$\kappa_{\ve,k}$ in the interior of $\gamma$ is therefore real. The same argument applies to
$\tilde{\BB}_{\ve, k}|_{\mathrm{Fix}\, \tilde{R}}$, whence the eigenspace corresponding to $\kappa_{\ve,k}$
lies in $\mathrm{Fix}\, \tilde{R}$.

This argument implies in particular that $\tilde{\BB}_{\ve, k}$ is closed
and that $\tilde\BB_{\ve,k}-\kappa I: \DD_\BB \rightarrow W$ is Fredholm with index $0$ for
$\kappa \in [-M,-m]$ (because $[-M,-m] \setminus \{\kappa_{\ve,k}\} \subseteq
\rho(\tilde{\BB}_{\ve,k})$ and $\kappa_{\ve,k}$ is a simple eigenvalue of $\tilde{\BB}_{\ve,k}$).
\qed

Our next result is deduced from inequality \eqn{BB estimate}
using the method given by Groves, Haragus \& Sun \cite[Lemma 3.15]{GrovesHaragusSun02}.

\begin{lemma}
The spectral projection $P_{\ve,k} \in \LL(W,\DD_\BB)$ defined by
\[
 P_{\ve,k} w=\frac{1}{2\pi \ii}\int_{\gamma}(\kappa I -\tilde{\BB}_{\ve, k})^{-1} w\dkappa
\]
satisfies the estimate
\[
 \|P_{\ve,k}-P_{0,k}\|_{\LL(W, \DD_\BB)} \le c\ve
\]
uniformly in $k\in [0, k_\mathrm{max}]$.
\end{lemma}

\begin{figure}[ht]
\hspace{3cm}\includegraphics[width=10cm]{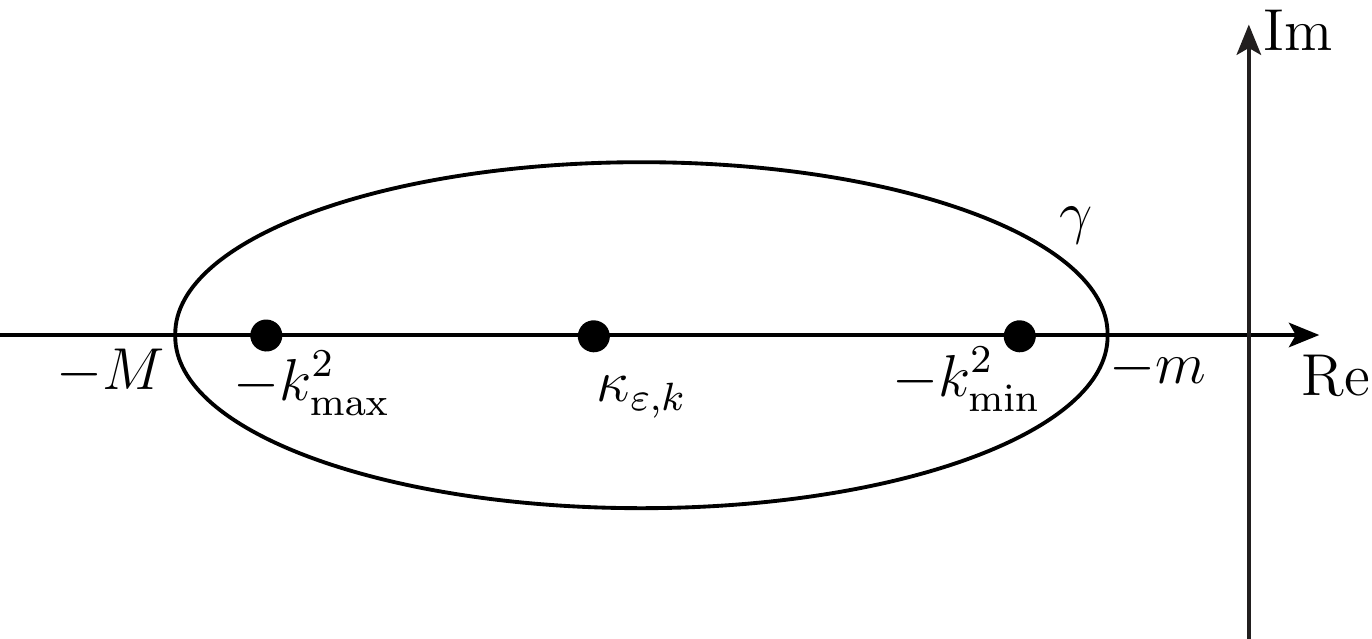}
{\it
\caption{The portion of $\sigma(\tilde\BB_{\ve, k})$ within $\gamma$ 
consists of precisely one simple real eigenvalue
$\kappa_{\ve, k} \in [-k_\mathrm{max}^2,-k_\mathrm{min}^2]$ with
$\kappa_{0, k}=-k_0^2$.
\label{separation}}}
\end{figure}

Finally, choose $k_\mathrm{max} >k_0$ and  $k_\mathrm{min} \in (0,k_0)$ and 
apply Lemma \ref{lemma:spectrum of tilde B}, choosing $m$, $M$ such that 
$M>k_\mathrm{max}^2$ and $m<k_\mathrm{min}^2$, so that the point $-k^2$ 
lies in the interior of $\gamma$ for all $k\in [k_\mathrm{min}, k_\mathrm{max}]$
(see Figure \ref{separation}). The following lemma follows from the relation
$\BB_{\ve, k}=\tilde{\BB}_{\ve, k}+k^2 I$.

\begin{lemma} \label{lemma:Evaleq}
Fix $k\in [k_\mathrm{min}, k_\mathrm{max}]$.
\begin{list}{(\roman{count})}{\usecounter{count}}
\item The operator $\BB_{\ve, k} : \DD_\BB \subseteq W \to W$ 
is invertible provided that $\kappa_{\ve, k}+ k^2\ne 0$.
\item The operator $\BB_{\ve, k} : \DD_\BB \subseteq W \to W$  has a simple zero eigenvalue whenever $\kappa_{\ve,  k}+ k^2=0$. 
\end{list}
\end{lemma}

In Section \ref{Analysis step 3} below we show that the equation $\kappa_{\ve,k}+k^2=0$
has precisely one solution $k_\ve$ in the interval $[k_\mathrm{min}, k_\mathrm{max}]$
(see Figure \ref{separation}). 

\subsection{Small spectral values} \label{Analysis step 2}

In this section we prove that $\BB_{\ve,k}|_{\mathrm{Fix}\, \tilde{R}}$ is semi-Fredholm and injective
for $k \in (0,k_\mathrm{min})$. Because this operator is closed (see Lemma \ref{lemma:spectrum of tilde B})
it suffices to establish the \emph{a priori} estimate
\[
\|w\|_W\le ck^{-2}\|\BB_{\ve,k}w\|_W, \qquad k \in (0, k_\mathrm{min}),
\]
for sufficiently small values of $k_\mathrm{min}$. For notational simplicity we here abbreviate
$\mathrm{Fix}\, \tilde{R}$ and $\DD_\BB \cap \mathrm{Fix}\, \tilde{R}$ to respectively
$W$ and $\DD_\BB$, so that $W=L_\mathrm{e}^2(\R) \times L_\mathrm{o}^2(\R) \times L_\mathrm{o}^2(\R)$ and
$\DD_\BB= H_\mathrm{e}^2(\R) \times H_\mathrm{o}^2(\R) \times H_\mathrm{o}^2(\R)$.

The first step is to choose an eigenvector corresponding to the eigenvalue $\kappa_{\ve,k}$ of $\BB_{\ve,k}$ in a perturbative fashion.

\begin{lemma}
There exists an eigenvector $w_{\ve, k}=(\zeta_{\ve,k},\psi_{\ve,k})$ of $\tilde{\BB}_{\ve, k}$ corresponding to the eigenvalue $\kappa_{\ve, k}$ such that $\|w_{\ve, k}\|_W=1$ and $\|w_{\ve,k}- w_{0,k}\|_{\DD_\BB}\le c\ve$. Furthermore, there exists a primitive
$\Psi_{\ve, k} \in H_\mathrm{e}^3(\R)$
of $\psi_{\ve, k}$ with
$\left\| \Psi_{\ve, k}- \Psi_{0, k} \right\|_0\le c\ve$.
\end{lemma}
{\bf Proof.}
Choose the eigenvector $w_{\ve,k}$ using the formula
\[
w_{\ve, k}=\|P_{\ve, k} w_0\|_W^{-1} P_{\ve, k} w_0,
\]
where $w_0$ is a unit eigenvector of $\BB_{k,0}$ corresponding to the eigenvalue
$-k_0^2$. (The calculation
\begin{align*}
 \|P_{\ve, k} w_0\|_W&=\|P_{0, k} w_0+(P_{\ve, k} -P_{0,k})w_0\|_W\\
&\geq\|w_0\|_W-\|(P_{\ve, k} -P_{0,k})w_0\|_W\\
&\geq 1-c\ve
\end{align*}
shows that $P_{\ve, k}w_0 \neq 0$.)
It follows that
\begin{align*}
\|w_{\ve, k}-w_{0, k}\|_{\DD_\BB}&\le 
\Big\|\|P_{\ve, k} w_0\|_W^{-1} P_{\ve, k} -P_{0,k}\Big\|_{\LL(W, \DD_\BB)} \|w_{0,k}\|_W\\
&\le 
\Big|\|P_{\ve, k} w_0\|_W^{-1}-1\Big|\| P_{\ve, k}\|_{\LL(W, \DD_\BB)}
+\|P_{\ve, k}- P_{0,k}\|_{\LL(W, \DD_\BB)}\\
&\leq c\ve. 
\end{align*}

Define
\[
\Psi_{\ve, k}=-\kappa_{\ve,k}^{-1} (1-\alpha_0^{-1})\psi_{\ve,k,x}-2\kappa_{\ve,k}^{-1} A_4 \zeta^\star \zeta_{\ve, k,1},
\]
so that
\begin{align*}
\Psi_{\ve,k,x} &= \kappa_{\ve,k}^{-1} \left(-(1-\alpha_0^{-1})\psi_{\ve,k,xx}-2A_4
(\zeta^\star \zeta_{\ve, k,1})_x\right) \\
&=\psi_{\ve,k}
\end{align*}
and
\begin{align*}
\|\Psi_{\ve,k}-\Psi_{\ve,0}\|_0
&= \left\|- (1-\alpha_0^{-1})\left(\frac{\psi_{\ve,k,x}}{\kappa_{\ve,k}}-\frac{\psi_{0,k,x}}{\kappa_{0,k}}\right)
-2A_4\left(\frac{\zeta_{\ve,k,1}}{\kappa_{\ve,k}}-\frac{\zeta_{0,k,1}}{\kappa_{0,k}}\right)\right\|_0 \\
& \leq c \left\| \frac{w_{\ve,k}}{\kappa_{\ve,k}}-\frac{w_{0,k}}{\kappa_{0,k}}\right\|_{\DD_\BB} \\
& \leq c\ve.
\end{align*}
\qed

We now define the auxiliary operator
$\AA_{\ve, k} : \DD_\AA \subseteq V \rightarrow V$ by the formula
\[
\AA_{\ve, k} \begin{pmatrix} \zeta_1 \\ \zeta_2 \\ \phi \end{pmatrix}=
\begin{pmatrix}
A_2^{-1}\zeta_1-A_2^{-1}A_1\zeta_{1xx}-A_2^{-1}(2A_3+A_5)(\zeta^\star)^2\zeta_1
+4A_2^{-1}A_4 \zeta^\star \phi-\re \RR_{\ve, k}(\zeta)\\
A_2^{-1}\zeta_2-A_2^{-1}A_1\zeta_{2xx}-A_2^{-1}A_5(\zeta^\star)^2\zeta_2-\im \RR_{\ve, k}(\zeta)\\
(1-\alpha_0^{-1})\phi+2A_4 \zeta^\star \zeta_1,
\end{pmatrix}
\]
where $V=L_\mathrm{e}^2(\R) \times L_\mathrm{o}^2(\R) \times L_\mathrm{e}^2(\R)$ and
$\DD_\AA=H_\mathrm{e}^2(\R) \times H_\mathrm{o}^2(\R) \times L_\mathrm{e}^2(\R)$.
Observe that $\AA_{0,k}$ is self-adjoint with respect to the inner product
$\llangle \cdot\,,\cdot\rrangle$ for $V$, and repeating the arguments used in the
proof of Lemma \ref{lemma:spec B0}, one finds that its spectrum consists of
essential spectrum $\{1-\alpha_0^{-1}\} \cup [A_2^{-1}, \infty)$ and a simple negative eigenvalue $-\omega_1^2$
(the point $1-\alpha_0^{-1}$ is an eigenvalue with infinite-dimensional eigenspace $\{(0,0)\} \times L_\mathrm{e}^2(\R)$).
Below we write $v=(\zeta,\psi_x) \in \DD_\AA$ for $w=(\zeta,\psi) \in \DD_\BB$ and in
particular set $v_{\ve,k}=(\zeta_{\ve,k},\psi_{\ve,k,x})$.

\begin{lemma} \label{lemma:A estimates}
The operator $\AA_{0,k} : \DD_\AA \subseteq V \rightarrow V$ satisfies
$$\llangle \AA_{0,k} v_{0,k},v_{0,k} \rrangle = -k_0^2, \qquad
\llangle \AA_{0,k} v_{0,k},v \rrangle = 0, \qquad
\llangle \AA_{0,k}v,v \rrangle \geq c\|v\|_V^2$$
for each $w \in \DD_\BB$ with $\llangle w, w_{0,k} \rrangle = 0$.
\end{lemma}
{\bf Proof.} Note that
\[
\AA_{0,k} (\zeta_{0,k}, \psi_{0,k,x}) = - k_0^2 (\zeta_{0,k}, -\Psi_{0,k})
\]
and therefore that
\begin{align}
\nonumber
\llangle \AA_{0, k} v_{0,k}, v_{0,k}\rrangle
&=-k_0^2 \llangle (\zeta_{0,k}, -\Psi_{0,k}), (\zeta_{0,k},\phi_{0,k})\rrangle\\
\nonumber
&=-k_0^2 \left( \|\zeta_{0,k}\|_0^2
-2A_2^{-1} \langle \Psi_{0,k},\psi_{0,k,x}\rangle_0\right)\\
\nonumber
&=-k_0^2 \|w_{0, k}\|_V^2\\
\label{eq:small spectral 1}
&=-k_0^2
\end{align}
and
\begin{align}
\nonumber
\llangle \AA_{0, k} v_{0,k}, v\rrangle&=
-k_0^2 \llangle  (\zeta_{0,k}, -\Psi_{0,k}), (\zeta,\psi_x)\rrangle\\
\nonumber
&=-k_0^2\left( \langle \zeta_{0,k}, \zeta\rangle_0
-2A_2^{-1} \langle \Psi_{0,k}, \psi_x \rangle_0 \right) \\
\nonumber
&=-k_0^2 \llangle  w_{0, k}, w\rrangle\\
\label{eq:small spectral 2}
&=0.
\end{align}

Let $\{P_\kappa\}_{\kappa \in \R}$ denote the spectral family of $\AA_{\ve,k}$ and write
\[
v_{0,k}=c_0 e +\int_{1-\alpha_0^{-1}}^\infty  \mathrm{d}(P_\kappa v_{0,k}), \qquad
v=d_0 e +\int_{1-\alpha_0^{-1}}^\infty  \mathrm{d}(P_\kappa  v).
\]
where $e$ is a unit eigenvector corresponding to the negative eigenvalue,
$c_0=\llangle v_{0,k}, e \rrangle$ and $d_0 = \llangle v,e \rrangle$,
so that
\begin{align}
\label{eq:small spectral 3}
\llangle \AA_{0, k} v_{0, k}, v_{0,k}\rrangle&=-\omega_1^2 c_0^2+\int_{1-\alpha_0^{-1}}^\infty \kappa\, 
\mathrm{d}(\llangle P_\kappa v_{0,k}, v_{0,k}\rrangle),
\\
\label{eq:small spectral 4}
\llangle \AA_{0, k} v, v\rrangle&=-\omega_1^2 d_0^2+\int_{1-\alpha_0^{-1}}^\infty \kappa\,  
\mathrm{d}(\llangle P_\kappa v, v\rrangle),
\\
\label{eq:small spectral 5}
\llangle \AA_{0, k} v_{0,k}, v\rrangle&=-\omega_1^2 c_0d_0+\int_{1-\alpha_0^{-1}}^\infty \kappa \, 
\mathrm{d}(\llangle P_\kappa v_{0,k}, v\rrangle).
\end{align}
It follows from \eqn{eq:small spectral 1} and \eqn{eq:small spectral 3} that
\[
c_0^2=\frac{1}{\omega_1^2 }\left(k_0^2+\int_{1-\alpha_0^{-1}}^\infty \kappa \,
\mathrm{d}(\llangle P_\kappa v_{0,k}, v_{0,k}\rrangle)\right) \neq 0,
\]
and equations \eqn{eq:small spectral 2} and \eqn{eq:small spectral 5} imply that
\[
\omega_1 d_0=\frac{1}{\omega_1 c_0} \int_{1-\alpha_0^{-1}}^\infty \kappa \,
\mathrm{d}(\llangle P_\kappa v_{0,k}, v\rrangle).
\]
The estimate
\begin{align*}
\llangle \AA_{0, k} v, v\rrangle
&=\frac1{c_0^2\omega_1^{2} }
\left\{-\left( \int_{1-\alpha_0^{-1}}^\infty \kappa \,
\mathrm{d}(\llangle P_\kappa v_{0,k}, v\rrangle)\right)^2\right.\\
&\qquad\qquad \qquad \left.
+\left(k_0^2+\int_{1-\alpha_0^{-1}}^\infty \kappa \,
\mathrm{d}(\llangle P_\kappa v_{0,k}, v_{0,k}\rrangle)\right)\int_{1-\alpha_0^{-1}}^\infty \kappa \,
\mathrm{d}(\llangle P_\kappa v, v\rrangle)\right\}\\
&\ge \frac{k_0^2}{c_0^2 \omega_1^2}\int_{1-\alpha_0^{-1}}^\infty \kappa \,
\mathrm{d}(\llangle P_\kappa v, v\rrangle) \\
&\geq 0,
\end{align*}
shows in particular that
\begin{align*}
d_0^2&\le \frac1{\omega_1^2}
\int_{1-\alpha_0^{-1}}^\infty \kappa\,  
\mathrm{d}(\llangle P_\kappa v, v\rrangle),
\end{align*}
(see equation \eqn{eq:small spectral 4}), whence
\begin{align*}
\|v\|_V^2&=
d_0^2 + \int_{1-\alpha_0^{-1}}^\infty 
\mathrm{d}(\llangle P_\kappa v, v\rrangle)\\
&\le \left(\frac{1}{\omega_1^2}+\frac{1}{1-\alpha_0^{-1}}\right)\int_{1-\alpha_0^{-1}}^\infty \kappa\,  
\mathrm{d}(\llangle P_\kappa v, v\rrangle)\\
&\le \frac{\omega_1^2 c_0^2}{k_0^2}\left(\frac{1}{\omega_1^2}+\frac{1}{1-\alpha_0^{-1}}\right)\llangle \AA_{0, k} v, v\rrangle.
\end{align*}
\qed

The \emph{a priori} estimate announced at the start of this section (Theorem
\ref{thm:small spectral values} below) is a corollary of the following lemma.

\begin{lemma} \label{lemma:small coercivity}
The operator $\BB_{\ve,k} : \DD_\BB\subseteq W \to W$ satisfies the 
estimate
\[
\llangle \BB_{\ve, k}w, w\rrangle\ge ck^2\|w\|_1^2, \qquad k \leq k_\mathrm{min}
\]
for each $w\in \DD_\BB$ such that $\llangle w, w_{\ve, k}\rrangle=0$
and each sufficiently small value of $k_\mathrm{min}$.
\end{lemma}
{\bf Proof.} Write
$$w=a w_{0, k}+w^\perp,$$
where $a=\llangle w, w_{0, k}\rrangle$, $\llangle w^\perp,w_{0,k} \rrangle=0$
and note that
\begin{align*}
|a|&=|\llangle w, w_{0,k}\rrangle| \\
&=|\llangle w,w_{0,k}-w_{\ve,k}\rrangle| \\
&= |\llangle v,((\zeta_{0,k}-\zeta_{\ve,k}),-(\Psi_{0,k}-\Psi_{\ve,k}))\rrangle| \\
&\leq c \ve \|v\|_V.
\end{align*}
Clearly
$$v=av_{0,k} + v^\perp,$$
so that
$$\|v^\perp\|_V^2 = \|v\|_V^2 - 2a \llangle v,v_{0,k} \rrangle +a^2,$$
and follows from Lemma \ref{lemma:A estimates} that
\begin{align*}
\llangle \AA_{0, k} v, v \rrangle&=
a^2\llangle \AA_{0, k} v_{0,k}, v_{0,k}\rrangle
+2a\llangle \AA_{0, k} v_{0,k}, v^\perp\rrangle
+\llangle \AA_{0, k} v^\perp, v^\perp\rrangle \\
& \geq c(\|v^\perp\|_V^2-a^2) \\
& \geq c(\|v\|_V^2 - 2a\llangle v, v_{0,k} \rrangle ) \\
& \geq c\|v\|_V^2
\end{align*}
because
$$|2a\llangle v,v_{0,k} \rrangle| \leq 2 |a| \|v\|_V \|v_{0,k}\|_V \leq 2c\ve \|v\|_V^2.$$

On the other hand
\begin{align*}
\llangle \AA_{0, k}v, v\rrangle
&=
\frac{1}{A_2}\int_\R (\zeta_1^2+A_1\zeta_{1x}^2-6\sech^2(A_1^{-1/2}x) \zeta_1^2)\dx \\
&\qquad\mbox{}+2A_2^{-1}(1-\alpha_0^{-1}) \int_{\R} \left(\phi+\frac{2A_4}{1-\alpha_0^{-1}}  \zeta^\star \zeta_1\right)^2\dx\\
&\qquad\mbox{}+\frac{1}{A_2}\int_\R (\zeta_1^2+A_1\zeta_{1x}^2-2\sech^2(A_1^{-1/2}x) \zeta_1^2)\dx \\
&\ge \frac{A_1}{A_2}\|\zeta_x\|_0^2- c \|\zeta\|_0^2,
\end{align*}
and combining the previous two estimates shows that
\begin{align*}
\llangle \AA_{0,k} v,v \rrangle
&\geq \frac{\delta A_1}{A_2} \|\zeta_x\|^2 - c \delta \|\zeta\|_0^2 + (1-\delta)c \|v\|^2 \\
&\geq c(\|\zeta\|_1^2 + \|\phi\|_0^2)
\end{align*}
for sufficiently small $\delta>0$. Using the estimate
$$\|\RR_{\ve,k}(\zeta)\|_0 \leq c \ve \|\zeta\|_1,$$
one finds that
\begin{align*}
\llangle \AA_{\ve,k} v,v \rrangle
&= \llangle \AA_{0,k}v, v \rrangle - \llangle (\RR_{\ve,k}(\zeta),0),(\zeta,0)\rrangle \\
&\geq c(\|\zeta\|_1^2 + \|\phi\|_0^2),
\end{align*}

We conclude that
\begin{align*}
\llangle \BB_{\ve,k}w,w \rrangle
&= \llangle \AA_{\ve,k} v,v \rrangle + k^2 \llangle w,w \rrangle \\
&\geq c (\|\zeta\|_1^2+\|\psi_x\|_0^2) + k^2(\|\zeta\|_0^2 + \|\psi\|_0^2) \\
& \geq k^2( \|\zeta\|_1^2 + \|\psi\|_1^2)
\end{align*}
for $k < k_\mathrm{min}$ and $k_\mathrm{min}$ sufficiently small.\qed

\begin{theorem}
\label{thm:small spectral values}
Choose $k \in (0, k_\mathrm{min}]$ and $w^\dag \in W$. Any solution $w\in \DD_\BB$ of the
equation
$$\BB_{\ve,k}w = w^\dag$$
satisfies the estimate
\[
\|w\|_W\le ck^{-2}\|w^\dag\|_W.
\]
\end{theorem}
{\bf Proof.}
Write
\[
w=c_{\ve, k} w_{\ve, k}+w^\perp, \quad w^\dag=d_{\ve, k} w_{\ve, k}+(w^\dag)^\perp,
\]
where $c_{\ve, k}=\llangle w, w_{\ve, k}\rrangle$, $d_{\ve, k}=\llangle w^\dag, w_{\ve, k}\rrangle$
and $\llangle w^\perp, w_{\ve,k} \rrangle =0$, $\llangle w^\dag, w_{\ve,k} \rrangle =0$, so that
\begin{equation}
c_{\ve, k}(\kappa_{\ve, k}+k^2)w_{\ve, k}+\BB_{\ve, k}w^\perp=
d_{\ve, k} w_{\ve, k}+(w^\dag)^\perp.
\label{Take inner prod}
\end{equation}
Taking the inner product of equation \eqn{Take inner prod} with $w^\perp$, we find that
\[
\llangle \BB_{\ve, k}w^\perp, w^\perp\rrangle
=\llangle (w^\dag)^\perp, w^\perp\rrangle,
\]
and it follows from Lemma \ref{lemma:small coercivity} that
\begin{equation}
\|w^\perp\|_1
\le \frac{c}{k^2}\|(w^\dag)^\perp\|_W.
\label{small coerc est 1}
\end{equation}
On the other hand, taking the inner product of \eqn{Take inner prod} with $w_{\ve,k}$ yields
\[
|c_{\ve, k}|=\frac{1}{|\kappa_{\ve, k}+k^2|}(|d_{\ve, k}|+|\llangle \BB_{0, k} w^\perp, w_{\ve, k}\rrangle|+|\llangle \RR_{\ve,k}(\zeta^\perp), \zeta_{\ve, k}\rrangle|),
\]
where $w^\perp = (\zeta^\perp, \psi^\perp)$,
because $\BB_{\ve,k}(w^\perp) = \BB_{0,k}(w^\perp) + (\RR_{\ve,k}(\zeta^\perp),0)$.
Estimating
$$|\llangle \BB_{0, k} w^\perp, w_{\ve, k}\rrangle| \leq c\|w^\perp\|_1 \|w_{\ve,k}\|_1,$$
where we have integrated by parts,
$$|\llangle \RR_{\ve,k}(\zeta^\perp), \zeta_{\ve, k}\rrangle|
\leq
\|\RR_{\ve,k}(\zeta^\perp)\|_0 \|\zeta_{\ve,k}\|_0 \leq c\|\zeta^\perp\|_1 \|\zeta_{\ve,k}\|_0,$$
and
$$|\kappa_{\ve,k}+k^2| \geq k_0^2-k_\mathrm{min}^2>0,$$
we find that
\begin{equation}
|c_{\ve,k}| \leq c (|d_{\ve,k}| + \|w^\perp\|_1).
\label{small coerc est 2}
\end{equation}

Estimates \eqn{small coerc est 1} and \eqn{small coerc est 2} show that
\begin{align*}
\|w\|_V
&=|c_{\ve,k}| + \|w^\perp\|_V \\
&\leq c\left(|d_{\ve,k}| + \frac{1}{k^2}\|(w^\dag)^\perp\|_V\right) \\
& \leq \frac{c}{k^2}\|w^\perp\|_V.
\end{align*}
\qed

\subsection{Solution of the eigenvalue equation} \label{Analysis step 3}

In this section we determine the solution set of the equation
\begin{equation}
 \label{eq:eigenvalue equation}
\kappa_{\ve, k}+k^2=0
\end{equation}
in the interval $[k_\mathrm{min}, k_\mathrm{max}]$.
The solution set is known 
explicitly for $\ve=0$, and we
approach the problem for $\ve>0$ using a perturbation argument.
Inequality \eqn{BB estimate} yields the
necessary estimate for $\kappa_{\ve,k}$,
and to obtain the corresponding estimate for
$\kappa^\prime_{\ve,k}$ we show that
the linear mapping $\RR_{\ve,k}:H^2(\R) \to L^2(\R)$ is a continuously differentiable function of $k>0$
which satisfies the estimate
$$\|\RR_{\ve, k}(\zeta)^\prime\|_0 \leq c \ve$$
uniformly over $k \in [k_\mathrm{min},k_\mathrm{max}]$ (which implies the corresponding
result for $\tilde{\BB}_{\ve,k}$).
This property of $\RR_{\ve,k}$ is established by systematically examining each step in its derivation
in Section \ref{sec:derivation} and checking the continuous differentiability of the operators appearing there.

\begin{proposition}
The linear mappings $\tilde{\Gamma}(\cdot,0)$, $\tilde{\Gamma}^{(j)}$, $j=1,2,3$, and $\tilde{\Gamma}^{\mathrm{r}}:H^2(\R) \to H^2(\Sigma)$
are continuously differentiable functions of $k>0$.
\end{proposition}
{\bf Proof.} Recall that
$$
\tilde{\Gamma}(\eta,u^\dag) = \GG^\star(\tilde{\Gamma}(\eta,u^\dag),\eta,u^\dagger),
$$
where
$$\GG^\star(\eta,\Gamma,u^\dagger)=\GG_1(F_1(\eta,\Gamma),F_2(\eta,\Gamma),F_3(\eta,\Gamma),\eta)+\GG_2(F^\dag(u^\dag)),$$
so that
\[
\tilde{\Gamma}(\eta,0)=(I-\GG^\star(0,\cdot,0))^{-1} \GG^\star(0,\eta,0).
\]
This formula and the fact that $k \mapsto \GG^\star(\cdot, \cdot, 0)$ is continuously differentiable
imply that
the map $k \mapsto \tilde{\Gamma}(\cdot,0)$ is also continuously differentiable.
The continuous differentiability of  $k\mapsto \tilde{\Gamma}^{(j)}(\cdot)$, $j=1,2,3$, follows from
that of $k \mapsto \GG_1(\cdot,\cdot,\cdot,0)$ and equations \eqn{eq:Phi1}--\eqn{eq:Phi3}, while
the continuous differentiability of $k \mapsto \tilde{\Gamma}^\mathrm{r}(\cdot)$
is a consequence of the formula
$\tilde{\Gamma}^\mathrm{r}(\eta)=\tilde{\Gamma}(\eta,0)-\tilde{\Gamma}^{(1)}(\eta)
-\tilde{\Gamma}^{(2)}(\eta)-\tilde{\Gamma}^{(3)}(\eta)$.\qed
\begin{corollary}
The linear mappings $\tilde{\NN}(\cdot,0)$, $\tilde{\NN}^{(2)}$, $\tilde{\NN}^{(3)}$ and $\tilde{\NN}^{\mathrm{r}}:H^2(\R) \to L^2(\R)$
are continuously differentiable.
\end{corollary}
\begin{lemma}
The linear mappings $\check{\eta}_2(\cdot,0)$, $\check \eta_2^{(2)}$ and $\check \eta_2^{\mathrm{r}}:\chi(D)L^2(\R) \to L^2(\R)$ are continuously differentiable
functions of $k>0$.
\end{lemma}
{\bf Proof.} Recall that
$$\check\eta_2(\eta_1,0) = \left(g_\ve(D,\ve k)-(1-\chi(D))\tilde{\NN}(\cdot)\right)^{-1}(1-\chi(D))\tilde{\NN}(\eta_1,0).$$
The map $k \mapsto g_\ve(D,\ve k)-(1-\chi(D))\tilde{\NN}(\cdot):
(1-\chi(D))H^2(\R) \to (1-\chi(D))L^2(\R)$ is continuously differentiable and
invertible (Lemma \ref{lem:invert}). It follows that its inverse and hence $\check \eta_2(\eta_1)$ 
is continuously differentiable with respect to $k$. 
The continuous differentiability of $k\mapsto \check \eta_2^{(2)}(\cdot)$ follows from that
of $k \mapsto g_\ve(D,\ve k)$, and the continuous differentiability
of $k \mapsto \check \eta_2^\mathrm{r}(\cdot)$
is a consequence of the formula
$\check \eta_2^\mathrm{r}(\cdot)(\eta_1) = \check{\eta}_2(\eta_1,0)-\check{\eta}_2^{(2)}(\eta_1)$.\qed
\begin{corollary}
$ $
\begin{list}{(\roman{count})}{\usecounter{count}}
\item
The linear mappings $\check{\Gamma}^{(j)}$, $j=1,2,3$ and
$\check{\Gamma}^\mathrm{r}:\chi(D)L^2(\R) \to H^2(\Sigma)$ are continuously differentiable functions of $k>0$.
\item
The linear mappings $\check{\NN}^{(2)}$, $\check{\NN}^{(3)}$ and $\check{\NN}^{\mathrm{r}}:\chi(D)L^2(\R) \to L^2(\R)$ are continuously differentiable
functions of $k>0$.
\end{list}
\end{corollary}
\begin{corollary}
The linear mappings $\tilde{\RR}_{\ve, k}$, $\RR_{\ve, k}: H^2(\R) \rightarrow L^2(\R)$ are
continuously differentiable functions of $k>0$.
\end{corollary}

Having established the continuous differentiability of $\tilde{\RR}_{\ve, k}(\zeta)$ and $\RR_{\ve, k}(\zeta)$, we
now turn to the task of estimating their derivatives, again by systematically examining each step in their derivation.
The following result, which is established in the same way as Proposition
\ref{prop:Green mapping}, is used in the first step.
\begin{proposition} \label{prop:tildeG}
The mapping
$(P_1,P_2,P_3,P_4,P_5) \mapsto \HH(P_1,P_2,P_3,P_4,P_5)$, where
\begin{align*}
\HH(P_1,&P_2,P_3,P_4,P_5)\\
&
=\FF^{-1}\left[\int_0^1\{G(y, \tilde y)(-2\ve^2 k \hat{P}_1-\ii \mu \hat P_2-\ii \ve \hat P_3 - \ii k \ve \hat{P}_4)+G_{\tilde y}(y,\tilde y) \hat P_5\} \dty\right]
\end{align*}
defines a linear function
$H^1(\Sigma)
\times H^1(\Sigma) \times H^1(\Sigma) \times H^1(\Sigma)
\times H^1(\Sigma) \rightarrow H^2(\Sigma)$ which satisfies the estimate
\begin{align*}
\|\nabla \HH\|_1 + \ve \|\HH\|_1 + \ve^2 \|\HH\|_0
&\leq c \left(\ve \|P_1\|_1+\sum_{j=2}^5 \|P_j\|_1
+ \ve^2 \|P_1\|_0 + \ve \sum_{j=2}^5 \|P_j\|_0 \right).
\end{align*}
\end{proposition}
\begin{lemma} \label{lemma:Phi diff}
The linear mappings $\tilde{\Gamma}(\cdot,0)$, $\tilde{\Gamma}^{(j)}$, $j=1,2,3$, and $\tilde{\Gamma}^{\mathrm{r}}:H^2(\R) \to H^2(\Sigma)$
satisfy the estimates
\begin{align*}
\|\nabla \tilde{\Gamma}(\eta,0)^\prime\|_{1}+ \ve \|\tilde{\Gamma}(\eta,0)^\prime\|_{0}
&\le c\|\eta\|_{2},\\
\|\nabla \tilde{\Gamma}^{(j)}(\eta)^\prime\|_{1}+\ve \|\tilde{\Gamma}^{(j)}(\eta)^\prime\|_{0}
&\le c\ve^{j-1}\|\eta\|_{2}, \qquad j=1,2,3,\\
\|\nabla \tilde{\Gamma}^{\mathrm{r}}(\eta)^\prime\|_{1}+\ve \|\tilde{\Gamma}^{\mathrm{r}}(\eta)^\prime\|_{0}
&\le c\ve^3\|\eta\|_{2}
\end{align*}
uniformly over $k \in [k_\mathrm{min},k_\mathrm{max}]$.
\end{lemma}
{\bf Proof.} Differentiating the boundary-value problem \eqn{eq:Phi BC1}--\eqn{eq:Phi BC3} with
respect to $k$, we find that
\begin{align*}
-\hat \Gamma_{yy}^\prime+q^2\hat \Gamma^\prime&=
-2\ve^2k\hat{\Gamma}-\ii \mu \hat F_1(\eta,\Gamma)^\prime
-\ii \ve F_2(\eta,\Gamma)
-\ii k \ve \hat F_2(\eta,\Gamma)^\prime -(\hat F_3(\eta,\Gamma))^\prime_y, && 0<y<1, \\
\hat \Gamma_y^\prime &=0 && \text{on } y=0, \\
\hat \Gamma_y^\prime &=\hat F_3(\eta,\Gamma)^\prime && \text{on } y=1,
\end{align*}
where $\Gamma^\prime$ is an abbreviation for $\tilde{\Gamma}(\eta,0)^\prime$
and $F_j(\eta,\Gamma)^\prime$ is also obtained from $F_j(\eta,\Gamma)$, $j=1,2,3$, by differentiating with
respect to $k$. It follows that
$$
\tilde{\Gamma}(\eta,0)^\prime=\HH(\tilde{\Gamma}(\eta,0),F_1(\eta,\tilde{\Gamma}(\eta,0))^\prime,F_2(\eta,\tilde{\Gamma}(\eta,0)),F_2(\eta,\tilde{\Gamma}(\eta,0))^\prime,F_3(\eta,\tilde{\Gamma}(\eta,0))^\prime)
$$
and similarly
\begin{align*}
\tilde{\Gamma}^{(1)}(\eta)^\prime &=\HH(\tilde{\Gamma}^{(1)}(\eta),0,0,0,0),\\
\tilde{\Gamma}^{(2)}(\eta)^\prime &=\HH(\tilde{\Gamma}^{(2)}(\eta),F_1^{(2)}(\eta)^\prime,F_2^{(2)}(\eta),F_2^{(2)}(\eta)^\prime,F_3^{(2)}(\eta)^\prime), \\
\tilde{\Gamma}^{(3)}(\eta)^\prime &=\HH(\tilde{\Gamma}^{(3)}(\eta),F_1^{(3)}(\eta)^\prime,F_2^{(3)}(\eta),F_2^{(3)}(\eta)^\prime,F_3^{(3)}(\eta)^\prime), \\
\tilde{\Gamma}^\mathrm{r}(\eta)^\prime &=\HH(\tilde{\Gamma}^\mathrm{r}(\eta),F_1^\mathrm{r}(\eta)^\prime,F_2^\mathrm{r}(\eta),F_2^\mathrm{r}(\eta)^\prime,F_3^\mathrm{r}(\eta)^\prime),
\end{align*}
where $F_j^{(2)}(\eta)^\prime$, $F_j^{(3)}(\eta)^\prime$, $F_j^\mathrm{r}(\eta)^\prime$
are obtained from $F_j^{(2)}(\eta)$, $F_j^{(3)}(\eta)$, $F_j^\mathrm{r}(\eta)$,
$j=1,2,3$, by differentiating with
respect to $k$. The estimates stated in the lemma are obtained by applying Proposition
\ref{prop:tildeG} to these formulae (estimates for $\tilde{\Gamma}(\eta,0)$ and $\tilde{\Gamma}^{(j)}(\eta)$,
$\tilde{\Gamma}^\mathrm{r}(\eta)$ are given in respectively 
Theorem \ref{thm:Phi} and Lemma \ref{lemma:tildePhi estimates}).\qed
\begin{corollary}
\label{cor:N derivatives}
The linear mappings $\tilde{\NN}(\cdot,0)$, $\tilde{\NN}^{(2)}$, $\tilde{\NN}^{(3)}$ and $\tilde{\NN}^{\mathrm{r}}:H^2(\R) \to L^2(\R)$ satisfy the estimates
\begin{align*}
\|\tilde{\NN}(\eta,0)^\prime\|_{0} 
&\le c\ve\|\eta\|_{2},\\
\|\tilde{\NN}^{(j)}(\eta)^\prime\|_{0}
&\le c\ve^{j-1}\|\eta\|_{2}, \qquad j=2,3,\\
\|\tilde{\NN}^{\mathrm{r}}(\eta)^\prime\|_{0} 
&\le c\ve^3 \|\eta\|_{2}
\end{align*}
uniformly over $k \in [k_\mathrm{min},k_\mathrm{max}]$.
\end{corollary}

The next estimates are obtained by differentiating the formulae defining $\check{\eta}_2(\eta_1,0)$,
$\check{\eta}_2^{(2)}(\eta_1)$ and $\check{\eta}_2^\mathrm{r}(\eta_1)$ with respect to $k$
and using Lemma \ref{lemma:tildeNestimates}, Proposition \ref{prop:checkestimates}(i), Corollary \ref{cor:N derivatives} and the fact that $|\partial_k g_\ve(\mu, \ve k)| \leq c$
for $k\in [k_\mathrm{min}, k_\mathrm{max}]$.

\begin{lemma}
The linear mappings $\check\eta_2(\cdot,0)$, $\check \eta_2^{(2)}$ and $\check \eta_2^{\mathrm{r}}:\chi(D)L^2(\R) \to L^2(\R)$ satisfy the estimates
$$
 \|(\check \eta_2(\eta_1,0))^\prime\|_2 \le c\ve \|\eta_1\|_0, \qquad
\|(\check \eta_2^{(2)}(\eta_1))^\prime\|_2 \le c\ve \|\eta_1\|_0,\qquad
\|(\check \eta_2^{\mathrm{r}}(\eta_1))^\prime\|_2 \le c\ve^2 \|\eta_1\|_0
$$
uniformly over $k \in [k_\mathrm{min},k_\mathrm{max}]$.
\end{lemma}
\begin{corollary}
$ $
\begin{list}{(\roman{count})}{\usecounter{count}}
\item
The linear mappings $\check{\Gamma}^{(j)}$, $j=1,2,3$ and
$\check{\Gamma}^\mathrm{r}:\chi(D)L^2(\R) \to H^2(\Sigma)$ satisfy the estimates
\begin{align*}
\|\nabla \check \Gamma^{(j)}(\eta_1)^\prime\|_{1}+\ve \|\check \Gamma^{(j)}(\eta_1)^\prime\|_{0}
&\le c\ve^{j-1}\|\eta_1\|_0, \qquad j=1,2,3,\\
\|\nabla \check \Gamma^{\mathrm{r}}(\eta_1)^\prime\|_{1}+\ve \|\check \Gamma^{\mathrm{r}}(\eta_1)^\prime\|_{0}
&\le c\ve^3\|\eta_1\|_0
\end{align*}
uniformly over $k \in [k_\mathrm{min},k_\mathrm{max}]$.
\item
The linear mappings $\check{\NN}^{(2)}$, $\check{\NN}^{(3)}$ and $\check{\NN}^{\mathrm{r}}:\chi(D)L^2(\R) \to L^2(\R)$ satisfy the estimates
\begin{align*}
\|\check{\NN}^{(j)}(\eta_1)^\prime\|_{0}
&\le c\ve^{j-1}\|\eta_1\|_{0}, \qquad j=2,3,\\
\|\check{\NN}^{\mathrm{r}}(\eta_1)^\prime\|_{0} 
&\le c\ve^3 \|\eta_1\|_{0}
\end{align*}
uniformly over $k \in [k_\mathrm{min},k_\mathrm{max}]$.
\end{list}
\end{corollary}

An estimate for $\check{\NN}^{(2)}(\eta_1)^\prime$ as a function of $\zeta$
is given in Lemma \ref{lemma:N^2diff}. This result is established
by differentiating \eqn{Eqn for checkNN2} and using Propositions
\ref{Expansion for Phi1}, \ref{Expansion for Phi2} and \ref{prop:checkPhi1 and tildePhi2 diff}
below, which is derived from the formulae
$$
\check{\Gamma}^{(1)\prime} = \FF^{-1}\left[\int_0^1 G(y,\tilde{y})(-2\ve^2 k \FF[\check{\Gamma}_1])\dty \right]
$$
and
\begin{align*}
\tilde{\Gamma}^{(2)}(\eta_1)^\prime
=\FF^{-1}\bigg[\int_0^1&\{G(y, \tilde y)(-2\ve^2 k \FF[\tilde{\Gamma}^{(2)}(\eta_1)]
-\ii \mu \FF[\check{F}_1^{(2)\prime}]-\ii\ve\FF[\check{F}_2^{(2)}]
-\ii k \ve \FF[\check{F}_2^{(2)\prime}]) \\
&\mbox{}+G_{\tilde y}(y,\tilde y) \FF[F_3^{(2)\prime}]\} \dty\bigg]
\end{align*}
by the methods used in Section \ref{Derivation step 3}.

\begin{proposition} \label{prop:checkPhi1 and tildePhi2 diff}
The quantities $\check\Gamma^{(1)\prime}$ and $\tilde{\Gamma}^{(2)}(\eta_1)^\prime$ satisfy the estimates
\begin{align*}
\check \Gamma^{(1)\prime} & = S, \\
\tilde{\Gamma}^{(2)}(\eta_1)^\prime &= -2\ve \FF^{-1}\left[\frac{\ii \mu k}{(\mu^2+k^2)^2}\FF[ \re(\zeta^\star\zeta)](\mu)\right](\ve x)\mu_0 \coth \mu_0+\ve S
\end{align*}
uniformly over $k \in [k_\mathrm{min},k_\mathrm{max}]$.
\end{proposition}
\begin{lemma} \label{lemma:N^2diff}
\label{lemma:check N^{(2)}'}
The quantity $\check{\NN}^{(2)}(\eta_1)^\prime$ satisfies the estimate
\begin{align*}
\check{\NN}^{(2)}(\eta_1)^\prime
= 2\ve^2\FF^{-1}\left[\frac{\mu^2k}{(\mu^2+k^2)^2}\FF[\re( \zeta_\delta^\star \zeta)](\mu)\right](\ve x)\frac{\mu_0 \cosh(\mu_0)}{\sinh(\mu_0)}
 + \ve R_0
\end{align*}
uniformly over $k \in [k_\mathrm{min},k_\mathrm{max}]$.
\end{lemma}
\begin{lemma} The quantity $\check \eta_2^{(2)\prime}$ satisfies the estimate
\begin{align*}
\check \eta_2^{(2)\prime}
= 2\ve^2  B_3\FF^{-1}\left[\frac{\mu^2k}{((1-\alpha_0^{-1})\mu^2+k^2)^2 }\FF[2\re(\zeta_\delta^\star\zeta)](\mu)\right](\ve x)
 +\ve R_2
\end{align*}
uniformly over $k \in [k_\mathrm{min},k_\mathrm{max}]$.
\end{lemma}
{\bf Proof.} Differentiating the formula
$$g_\ve(\mu,\ve k)\FF[\check{\eta}^{(2)}] = (1-\chi)\FF[\check{\NN}^{(2)}(\eta_1)]=\FF[\check{\NN}^{(2)}(\eta_1)]$$
with respect to $k$ yields
\begin{align*}
g(\mu, \ve k)\FF[\check \eta_2^{(2)^\prime}]&=\FF[\check \NN^{(2)}(\eta_1)^\prime]
-\partial_k g_\ve(\mu, \ve k)\FF[\check \eta^{(2)}(\eta_1)]\\
&=2\ve \frac{(\ve^{-1}\mu)^2k}{((\ve^{-1}\mu)^2+k^2)^2}\FF[\re( \zeta_\delta^\star \zeta)](\ve^{-1}\mu)\frac{\mu_0 \cosh(\mu_0)}{\sinh(\mu_0)}
\\
 &\qquad \mbox{}+\ve\partial_k g_\ve(\mu, \ve k)\left(\frac{1}{2}B_1 \FF[\zeta_\delta^\star\zeta](\ve^{-1}(\mu-2\mu_0))
+\frac{1}{2}B_1 \FF[\zeta_\delta^\star\overline{\zeta}](\ve^{-1}(\mu+2\mu_0))
 \right)
\\
 &\qquad \mbox{}+\ve\partial_k g_\ve(\mu, \ve k)\left(B_2
 + B_3 \frac{2(\ve^{-1}\mu)^2}{(1-\alpha_0^{-1})(\ve^{-1}\mu)^2+k^2 }\right)\FF[\re(\zeta_\delta^\star\zeta)](\ve^{-1}\mu)
\\
&\qquad +\ve(1-\chi)\hat R_0,
\end{align*}
in which we have used Proposition \ref{prop:check eta_22} and
Lemma \ref{lemma:check N^{(2)}'}. Noting that
$$\partial_k g_\ve(\mu,\ve k)= \frac{2k(\ve^{-1}\mu)^2}{((\ve^{-1}\mu)^2+k^2)^2}+O(\ve)$$
uniformly over $\{|\mu| \in [0,\mu_\mathrm{max}]\}$ for each fixed $\mu_\mathrm{max}>0$
(in particular for $\mu_\mathrm{max}=2\delta$), whence
$\partial_k g_\ve(\mu,\ve k) = O(\ve)$ uniformly over
$\{|\mu| \in [\mu_\mathrm{min},\mu_\mathrm{max}]\}$ for each fixed $\mu_\mathrm{min}$,
$\mu_\mathrm{max}>0$ (in particular for $\mu_\mathrm{min}=2\mu_0-2\delta$,
$\mu_\mathrm{max}=2\mu_0+2\delta$), we find that
\begin{align*}
& g(\mu, \ve k)\FF[\check \eta_2^{(2)^\prime}] \\
&=2\ve \frac{(\ve^{-1}\mu)^2k}{((\ve^{-1}\mu)^2+k^2)^2}
\left(\frac{\mu_0\cosh(\mu_0)}{\sinh(\mu_0)}+B_2
 + B_3 \frac{2(\ve^{-1}\mu)^2}{(\alpha_0-1)(\ve^{-1}\mu)^2+\alpha_0 k^2 }\right)\FF[\re(\zeta_\delta^\star\zeta)](\ve^{-1}\mu) \\
 &\qquad\mbox{}+\ve(1-\chi) \hat{R}_0.
 \end{align*}

Finally, we use the estimate
\[
g_\ve(\mu,k \ve)^{-1}=\frac{(\ve^{-1}\mu)^2+ k^2}{(\alpha_0-1)(\ve^{-1}\mu)^2+\alpha_0 k^2}+T_2(\mu,\ve),
\]
where
$$|T_2(\mu,\ve)| \leq c(\mu^2+\ve^2), \qquad |\mu|<2\delta,$$
and the fact that
$$g_\ve(D,\lambda)^{-1}(1-\chi(D))R_0=R_2$$
because $g_\ve(\mu,k\ve) > c(1+q^2)$ for $\mu \in \supp(1-\chi)$ (see Proposition \ref{prop:g}).\qed

Finally, we compute estimates for $\check \Gamma^{(2)\prime}$ and
$\tilde{\Gamma}^{(3)}(\eta_1)^\prime+ \tilde{\Gamma}^{(2)}(\check\eta_2^{(2)})^\prime$
using the formulae
$$
\check \Gamma^{(2)\prime} = \tilde{\Gamma}^{(2)}(\eta_1)^\prime + \tilde{\Gamma}^{(1)}(\check \eta_2^{(2)})^\prime, \qquad
\tilde{\Gamma}^{(1)}(\check \eta_2^{(2)})^\prime= \FF^{-1}\left[\int_0^1 G(y,\tilde{y})(-2\ve^2 k \FF[\tilde{\Gamma}^{(1)}(\check \eta_2^{(2)})])\dty \right]
$$
and
\begin{align*}
&\tilde{\Gamma}^{(3)}(\eta_1)^\prime+ \tilde{\Gamma}^{(2)}(\check\eta_2^{(2)})^\prime\\
&=\FF^{-1}\bigg[\int_0^1\{G(y, \tilde y)(-2\ve^2 k \FF[\tilde{\Gamma}^{(3)}(\eta_1)+ \tilde{\Gamma}^{(2)}(\check\eta_2^{(2)})]
-\ii \mu \FF[\check{F}_1^{(3)\prime}]-\ii\ve\FF[\check{F}_2^{(3)}]
-\ii k \ve \FF[\check{F}_2^{(3)\prime}]) \\
&\qquad\qquad\qquad\mbox{}+G_{\tilde y}(y,\tilde y) \FF[F_3^{(3)\prime}]\} \dty\bigg]
\end{align*}
and deduce an estimate for $\check{\NN}(\eta_1)^\prime$ by differentiating
equation \eqn{Eqn for checkNN3}.

\begin{proposition}
$ $
\begin{list}{(\roman{count})}{\usecounter{count}}
\item
The quantity $\check \Gamma^{(2)}(\eta_1)^\prime$ satisfies the estimate
\begin{align*}
\check \Gamma^{(2)}(\eta_1)^\prime
=-\ve(\vartheta_\delta(\ve x))'M_2+\ve S
\end{align*}
uniformly over $k \in [k_\mathrm{min},k_\mathrm{max}]$, where
\[
(\vartheta_\delta(x))^\prime=-2A_4
\FF^{-1}\left[\frac{\ii \mu k}{((1-\alpha_0^{-1})\mu^2+k^2)^2}\FF[ 2\re(\zeta^\star_\delta \zeta)](\mu)\right](x).
\]
\item
The quantity 
$\tilde{\Gamma}^{(3)}_x(\eta_1)^\prime + \tilde{\Gamma}^{(2)}_x(\check \eta_2^{(2)})^\prime\Big|_{y=1}$
satisfies the estimate
$$
\tilde{\Gamma}^{(3)}_x(\eta_1)^\prime + \tilde{\Gamma}^{(2)}_x(\check \eta_2^{(2)})^\prime\Big|_{y=1}
=\mu_0 \ve^3 \zeta_\delta^\star(\ve x)\vartheta_{\delta x}(\ve x)^\prime N_3 \left(\frac{1}{2}\ee^{\ii \mu_0 x}+\frac{1}{2}\ee^{-\ii \mu_0 x}\right) + \ve^2 T$$
\end{list}
uniformly over $k \in [k_\mathrm{min},k_\mathrm{max}]$.
\end{proposition}
\begin{lemma}
\label{lemma:N^3diff}
The quantity $\check{\NN}^{(3)}(\eta_1)^\prime$ is given by the formula
$$
\check{\NN}^{(3)}(\eta_1)^\prime=
-4\ve^3 A_4\zeta_\delta^\star(\ve x)\vartheta_{\delta x}(\ve x)^\prime
 \left(\frac{1}{2}\ee^{\ii \mu_0 x}+\frac{1}{2}\ee^{-\ii \mu_0 x}\right)+\ve^2 R
$$
uniformly over $k \in [k_\mathrm{min},k_\mathrm{max}]$.
\end{lemma}

We now prove our main result using Lemmata \ref{lemma:N^2}, \ref{lemma:N^3},
\ref{lemma:N^2diff} and \ref{lemma:N^3diff}.

\begin{theorem}
The linear mappings $\tilde{\RR}_{\ve, k}$, $\RR_{\ve, k}: H^2(\R) \rightarrow L^2(\R)$ satisfy the estimate
$$\|\tilde{\RR}_{\ve, k}(\zeta)^\prime\|_0,\ \|\RR_{\ve, k}(\zeta)^\prime\|_0 \leq c \ve \|\zeta\|_1$$
uniformly over $k \in [k_\mathrm{min},k_\mathrm{max}]$.
\end{theorem}
{\bf Proof.}
Recall that
\begin{align*}
\tilde\RR_{\ve, k}(\zeta)&=A_2^{-1}\FF^{-1}\left[2\frac{\tilde g_\ve(\mu_0+\ve\tilde \mu, \ve k)}{ g_\ve(\mu_0+\ve \tilde \mu,\ve k)}
\chi_0(\ve \tilde \mu)\ve^{-2} \FF[\check{\NN}^{(3)}(\chi(D)\eta_1)+\check{\NN}^\mathrm{r}(\chi(D)\eta_1)](\mu_0+\ve \tilde \mu)\right]\\
&\qquad\mbox{} - A_2^{-1}E(x),
\end{align*}
so that
\begin{align*}
\tilde\RR_{\ve, k}(\zeta)^\prime&=
A_2^{-1}\FF^{-1}\left[2\left(\frac{\tilde g_\ve(\mu_0+\ve \tilde \mu, \ve k)}{ g_\ve(\mu_0+\ve \tilde \mu,\ve k)}\right)^\prime\chi_0(\ve \tilde \mu)\ve^{-2} \FF[\check{\NN}^{(3)}(\chi(D)\eta_1)+\check{\NN}^\mathrm{r}(\chi(D)\eta_1)](\mu_0+\ve \tilde \mu)\right]\\
&\qquad\mbox{}+
A_2^{-1}\FF^{-1}\left[2\frac{\tilde g_\ve(\mu_0+\ve \tilde \mu, \ve k)}{ g_\ve(\mu_0+\ve \tilde \mu,\ve k)}
\chi_0(\ve \tilde \mu)\ve^{-2} \FF[\check{\NN}^{(3)}(\chi(D)\eta_1)^\prime+\check{\NN}^\mathrm{r}(\chi(D)\eta_1)^\prime](\mu_0+\ve \tilde \mu)\right] \\
&\qquad\mbox{}- A_2^{-1}E(x)^\prime.
\end{align*}
Noting that
$$
\check{\NN}^{(3)}(\chi(D)\eta_1)^\prime
=\frac{\ve^3}{2}E(\ve x)^\prime\ee^{\ii\mu_0 x}
+\frac{\ve^3}{2}F(\ve x)^\prime\ee^{-\ii\mu_0 x}+\ve^2T,
$$
one finds by repeating the arguments below equation \eqn{eq:checkNN3 estimate to E} that
\begin{align*}
&\left\|\FF^{-1}\left[2\frac{\tilde g_\ve(\mu_0+\ve \tilde \mu, \ve k)}{ g_\ve(\mu_0+\ve \tilde \mu,\ve k)}
\chi_0(\ve \tilde \mu)\ve^{-2} \FF[\check{\NN}^{(3)}(\chi(D)\eta_1)^\prime+\check{\NN}^\mathrm{r}(\chi(D)\eta_1)^\prime](\mu_0+\ve \tilde \mu)\right] - A_2^{-1}E(x)^\prime \right\|_0\\
&\leq c \|\zeta\|_1.
\end{align*}
Moreover, the calculation
$$
 \left(\frac{\tilde g_\ve(\mu, \ve k)}{ g_\ve(\mu,\ve k)}\right)'
=\frac{1}{g_\ve(\mu, \ve k)}\left(\left(1-\frac{\tilde g_\ve(\mu, \ve k)}{ g_\ve(\mu,\ve k)}\right)
\partial_k g_\ve(\mu,\ve k)+\partial_k g_\ve(\mu,\ve k)-\partial_k \tilde{g}_\ve(\mu,\ve k)\right)
$$
and the estimates
$$c|(\mu-\mu_0,\ve k)|^2 \leq g_\ve(\mu,\ve k) \leq \frac{1}{c}|(\mu-\mu_0,\ve k)|^2,
\qquad
\left|\frac{\tilde{g}_\ve(\mu,\ve k)}{g_\ve(\mu,\ve k)}-1\right|
\leq c|(\mu-\mu_0,\ve k)|^2$$
$$|\partial_k g(\mu,\ve k)| \leq c\ve, \qquad
|\partial_k g_\ve(\mu,\ve k)-\partial_k \tilde{g}_\ve(\mu,\ve k)| \leq c \ve|(\mu-\mu_0,\ve k)|^3$$
for $|\mu-\mu_0|<\delta$ imply that
$$\left|\left(\frac{\tilde g_\ve(\mu_0+\ve \tilde \mu, \ve k)}{g_\ve(\mu_0+\ve \tilde \mu, \ve k)} \right)^\prime
\chi_0(\ve \tilde \mu)\right| \leq c\ve,$$
whence
\begin{align*}
&\left\|\FF^{-1}\left[2\left(\frac{\tilde g_\ve(\mu_0+\ve \tilde \mu, \ve k)}{ g_\ve(\mu_0+\ve \tilde \mu,\ve k)}\right)'\chi_0(\ve \tilde \mu)\ve^{-2} \FF[\check{\NN}^{(3)}(\eta_1)+\check{\NN}^\mathrm{r}(\eta_1)](\mu_0+\ve \tilde \mu)\right]
\right\|_0 \\
&\leq c \ve^{-3/2}(\|\check{\NN}^{(3)}(\chi(D)\eta_1)\|_0+\|\check{\NN}^\mathrm{r}(\chi(D)\eta_1)\|_0) \\
& \leq c\ve^{1/2}\|\eta_1\|_1 \\
& \leq c\ve\|\zeta\|_1.
\end{align*}

According to the above calculation
$$\|\tilde\RR_{\ve,k}(\zeta)^\prime\|_0 \leq c \ve\|\zeta\|_1,$$
and the observation
$$\RR_{\ve,k}(\zeta)^\prime-\tilde{\RR}_{\ve,k}^\prime(\zeta) = A_4(\zeta_\delta^\star\vartheta_{\delta x}^\prime
-\zeta\vartheta_x^\prime),$$
shows that
$$\|\RR_{\ve,k}(\zeta)^\prime-\tilde{\RR}_{\ve,k}^\prime(\zeta)\|_0 \leq c \ve \|\zeta\|_1.\eqno{\Box}$$

\begin{corollary}
The linear mapping $\tilde{\BB}_{\ve, k}: \DD_\BB \rightarrow V$ satisfies the estimate
$$\|\tilde{\BB}_{\ve, k}^\prime\|_{\LL(\DD_\BB,V)} \leq c \ve$$
uniformly over $k \in [k_\mathrm{min},k_\mathrm{max}]$.
\end{corollary}

The results in the next lemma are deduced from the estimates
$$\|\tilde\BB_{\ve,k}-\tilde\BB_{0,k}\|_{\LL(\DD_\BB,V)} \leq c\ve, 
\qquad 
\|\tilde{\BB}_{\ve, k}^\prime\|_{\LL(\DD_\BB,V)} \leq c \ve$$
using the method given by Groves, Haragus \& Sun
\cite[Lemma 3.15]{GrovesHaragusSun02}.

\begin{lemma}
The eigenvalue $\kappa_{\ve, k}$ of $\tilde{\BB}_{\ve,k}: \DD_\BB \subseteq V \rightarrow V$
is a differentiable function of $k$ and satisfies the estimates
\[
|\kappa_{\ve, k}-\kappa_{0, k}|\le c^\star\ve, \quad |\kappa_{\ve, k}'|\le c^\star\ve
\]
uniformly in $k\in [k_\mathrm{min}, k_\mathrm{max}]$.
\end{lemma}

Finally, the solution set of \eqn{eq:eigenvalue equation} in the interval $[k_\mathrm{min},
k_\mathrm{max}]$ can be determined from the estimates presented in the previous lemma
using a perturbation argument and the contraction-mapping principle (cf.\ 
Groves, Haragus \& Sun \cite[Theorem 3.16]{GrovesHaragusSun02}).

\begin{theorem}
Equation \eqn{eq:eigenvalue equation} has precisely one solution $k_\ve$ in the interval 
$[k_\mathrm{min}, k_\mathrm{max}]$. This solution lies in the set
\[
\mathcal S=\left\{k: |k-k_0|\le \frac{c^\star\ve}{k_0}\right\}.
\]
\end{theorem}

\appendix

\section*{Appendix: The Iooss-Kirchg\"{a}ssner line solitary waves}
An asymptotic expansion of the Iooss-Kirchg\"{a}ssner line solitary waves
in powers of $\ve$ may be computed in the framework of their existence theory
(see Dias \& Iooss \cite{DiasIooss93}, noting the omission of certain terms in
equation (3.41)). One finds that
$$
\eta^\star(x)=\eta^\star_1(x)+\eta^\star_2(x)+\eta^\star_\mathrm{r}(x),\qquad
\Phi^\star(x,y)=\Phi^\star_1(x,y)+\Phi^\star_2(x,y)+ \Phi^\star_\mathrm{r}(x,y),
$$
where
\begin{eqnarray}
\eta^\star_1(x)&=&
\ve \zeta^\star(\ve x)\cos(\mu_0 x), \label{Defn eta1star}\\
\eta^\star_2(x)&=& C_0\ve^2 \zeta^\star(\ve x) \xi^\star(\ve x)\sin (\mu_0 x)
-C_1\ve^2(\zeta^\star)^2(\ve x)\cos (2\mu_0 x) \nonumber \\
& & \qquad\mbox{}-C_2 \ve^2(\zeta^\star)^2(\ve x)+C_3
\ve^2 (\zeta^\star)'(\ve x)\sin (\mu_0 x), \label{Defn eta2star}\\
\Phi^\star_1(x,y)&=&
 \ve\zeta^\star(\ve x) \sin(\mu_0x)\frac{\cosh(\mu_0y)}{\sinh(\mu_0)}, \label{Defn Phi1star}\\
\Phi^\star_2(x,y)&=&
-C_0 \ve^2 \xi^\star(\ve x)\cos (\mu_0 x) \frac{\cosh (\mu_0 y)}{\sinh \mu_0}
+ \ve^2 (\zeta^{\star})^2(\ve x)\sin(2\mu_0x)\frac{\mu_0 y\sinh(\mu_0y)}{2\sinh(\mu_0)} \nonumber \\
& & \qquad\mbox{} - C_4 \ve \partial_x^{-1}(\zeta^{\star 2})(\ve x)
- C_5 \ve^2 (\zeta^{\star})^2(\ve x)\sin(2\mu_0x)\frac{\cosh(2\mu_0y)}{\sinh(2\mu_0)}, \label{Defn Phi2star}
\end{eqnarray}
and 
$\eta_\mathrm{r}^\star$, $\nabla \Gamma_\mathrm{r}^\star = \underline{O}(\ve^3)$,
the symbol $\underline{O}(\ve^\alpha)$ denoting a smooth quantity whose derivatives are all
$O(\ve^\alpha \ee^{-\rho \ve |x|})$  for some $\rho>0$ (uniformly over $y \in [0,1]$). 
The functions $\zeta^\star$ and $\xi^\star$ are given by the formulae
$$\zeta^\star(x) =
\left(\frac{2}{A_5}\right)^{\!\! 1/2}\sech\left(\frac{x}{A_1^{1/2}}\right), \qquad \xi^\star(x) = x \zeta^\star(x)$$
with positive constants
\[
A_1=\beta_0 +(1-\mu_0 \coth(\mu_0))\cosech^2(\mu_0),  \qquad A_5=A_3+4(1-\alpha_0^{-1})^{-1}A_4^2,
\]
where
\begin{eqnarray*}
A_3 & = & -\frac{\mu_0^3}{8\sigma^3} \left(\frac{(1-\sigma^2)(9-\sigma^2)\alpha_0+\beta_0 \mu_0^2(3-\sigma^2)(7-\sigma^2)}{\alpha_0 \sigma^2-\beta_0 \mu_0^2(3-\sigma^2)}+8\sigma^2-\frac{2\mu_0}{\alpha_0\sigma}(1-\sigma^2)^2-3\beta_0 \mu_0 \sigma^3\right),\\
A_4 & = & \frac{\mu_0(\alpha_0 \sinh(2\mu_0)+\mu_0)}{4\alpha_0\sinh^2(\mu_0)}
\end{eqnarray*}
and
\[
\sigma=\tanh(\mu_0)=\frac{\mu_0}{\alpha_0+\beta_0 \mu_0^2}.
\]
Finally, the positive coefficients $C_1$, $C_2$, $C_4$ and $C_5$ are given by
\begin{align*}
C_1&=\frac{\mu_0^2(\cosh(2\mu_0)+2)}{4\sinh^2(\mu_0)g_0(2\mu_0,0)}, &
C_2&=\frac{\mu_0(\sinh(2\mu_0)+\mu_0)}{4\sinh^2(\mu_0)(\alpha_0-1)},\\
C_4&=\frac{\mu_0(\alpha_0 \sinh(2\mu_0)+\mu_0)}{4\sinh^2(\mu_0)(\alpha_0-1)}, &
C_5&=\frac{\mu_0^2 (\cosh(2\mu_0)+2)}{4\sinh^2(\mu_0) g_0(2\mu_0,0)}+
\frac{\mu_0\sinh(2\mu_0)}{4\sinh^2(\mu_0)},
\end{align*}
while the values (and signs) of the coefficients $C_0$ and $C_3$ are unimportant.

\noindent\\
{\bf Acknowledgements.} E. Wahl\'{e}n was supported by the Swedish Research Council (grant no. 621-2012-3753). We would like to thank Mariana Haragus (Universit\'e de Franche-Comt\'e Besan\c{c}on) for many helpful discussions during the preparation of this article.

\bibliographystyle{standard}
\bibliography{mdg}

\end{document}